Ͷarxiv¼ÓÉÏÕâÒ»ÐÐ
\UseRawInputEncoding
\documentclass{amsart}

\usepackage{amssymb,amsmath,amsthm}
\usepackage{amsaddr}
\usepackage{amsfonts}
\usepackage{authblk}
\usepackage{graphicx}
\usepackage{epstopdf}
\usepackage{caption}
\usepackage{subcaption}
\usepackage{hyperref}
\usepackage{color}
\usepackage{marginnote}
\usepackage{float}
\usepackage{graphicx}
\usepackage{comment}
\usepackage{tikz-cd}
\usetikzlibrary{arrows}
\usetikzlibrary{intersections,positioning}
\tikzset{>=latex}
\usepackage{lipsum}%
\usepackage{a4wide}

\newcommand{\C}{\mathbb{C}}

\newcommand{\R}{\mathbb{R}}
\newcommand{\Z}{\mathbb{Z}}
\newcommand{\Q}{\mathbb{Q}}
\newcommand{\N}{\mathbb{N}}

\newcommand{\U}{\mathcal U}

\newcommand{\ov} {\overline}

\newcommand{\Oc} {\mathcal{O}}

\newcommand{\om}{\omega}

\newcommand{\ap} {\alpha}
\newcommand{\bt} {\beta}

\newcommand{\ep} {\epsilon}

\newcommand{\sg} {\sigma}
\newcommand{\Sg} {\Sigma}

\newcommand{\defi}{\stackrel{\rm def}{=}}

\newcounter{newcounter}[section]
\numberwithin{equation}{section}
\numberwithin{newcounter}{section}
\numberwithin{figure}{section}
\numberwithin{footnote}{section}

\newcommand{\authorfootnotes}{\renewcommand\thefootnote{\@fnsymbol\c@footnote}}%

\newtheorem{thm}[newcounter]{Theorem}
\newtheorem{prop}[newcounter]{Proposition}
\newtheorem{lem}[newcounter]{Lemma}

\newtheorem{rem}[newcounter]{Remark}

\title{Non-resonant Hopf links near a Hamiltonian equilibrium point}
\date{\today}

\begin{document}

\email{ragazzo@usp.br}

\email{liulei30@email.sdu.edu.cn}

\email{psalomao@sustech.edu.cn}

\maketitle

\begin{center}

\normalsize
\authorfootnotes
C. Grotta-Ragazzo\textsuperscript{1}, Lei Liu\textsuperscript{2} and
Pedro A. S. Salom\~ao\textsuperscript{3}  \par \bigskip

\textsuperscript{1}Instituto de Matem\'atica e Estat\'istica, University of S\~ao Paulo

\textsuperscript{2}School of Mathematics, Shandong University  \par

\textsuperscript{3} Shenzhen International Center for Mathematics, SUSTech

\end{center}

\begin{abstract}
This paper is about the existence of periodic orbits near an equilibrium point of a two-degree-of-freedom Hamiltonian system. The equilibrium is supposed to be a nondegenerate minimum of the Hamiltonian. Every sphere-like component of the energy surface sufficiently close to the equilibrium contains at least two periodic orbits forming a Hopf link (A. Weinstein \cite{weinstein}).
A theorem by Hofer, Wysocki, and Zehnder \cite{HWZ1998} implies that there are either precisely two or infinitely many periodic orbits on such a component of the energy surface. This multiplicity result follows from the existence of a disk-like global surface of section. If a certain non-resonance condition on the rotation numbers of the orbits of the Hopf link is satisfied \cite{HMS2015}, then infinitely many periodic orbits follow. This paper aims to present explicit conditions on the Birkhoff-Gustavson normal forms of the Hamiltonian function at the equilibrium point that ensure the existence of infinitely many periodic orbits on the energy surface by checking the non-resonance condition as in \cite{HMS2015} and not making use of any global surface of section. The main results focus on strongly resonant equilibrium points and apply to the Spatial Isosceles Three-Body Problem, Hill's Lunar Problem, and the H\'enon-Heiles System.
\end{abstract}

\tableofcontents

\section{Introduction}
This paper is about periodic orbits near an equilibrium point of a two-degree-of-freedom Hamiltonian system. At the equilibrium,  the Hessian of the Hamiltonian function is supposed to be positive definite. In this situation,  there exist local canonical coordinates near the equilibrium such that the Hamiltonian function becomes
\begin{equation}
H=\frac{\alpha_1}{2}(x_1^2+y_1^2)+\frac{\alpha_2}{2}(x_2^2+y_2^2)+\Oc_3,
\label{H1}
\end{equation}
where $\alpha_1$ and $\alpha_2$ are positive constants and $\Oc_3$ denotes a function with all partial derivatives at the origin up to second order equal to zero. The symplectic form is $\omega_0=dy_1\wedge dx_1+dy_2\wedge dx_2$, and the Hamiltonian vector field $X$ is determined by $-i_X \omega_0=dH$. There exists a neighborhood $\U$ of the origin where the values of $H$ are positive except at the origin and
the level set of $H$ in $\U$, denoted by  $\Sg_E=H^{-1}(E)\cap \U$, is a strictly convex sphere-like hypersurface for every $E>0$ small. We restrict our attention to this neighborhood.
A theorem due to A. Weinstein \cite{weinstein} ensures that
 $\Sg_E$ contains at least two periodic orbits forming a Hopf link for every $E>0$ sufficiently small.   There are several
extensions of this theorem showing that, under suitable assumptions,
$\Sg_E$ contains periodic orbits even for larger values of $E$. A particularly important result in this direction is due to
Hofer, Wysocki and Zehnder \cite{HWZ1998}. They proved that if $\Sg_E$ is strictly convex, then
it contains a special periodic orbit that will be denoted by $\gamma_1$, which is the boundary
of a disk-like global surface of section of the Hamiltonian
flow restricted to $\Sg_E$.  This reduces the study of the topological properties of the flow to the study of an area-preserving map of the open disk. Brouwer's translation theorem implies that this map has at least one fixed point corresponding to a periodic orbit of the flow denoted by $\gamma_2$. The link $\gamma_1 \cup \gamma_2\subset \Sigma_E$ is a Hopf link.
A theorem due to Franks \cite{franks} states that if this map has a second periodic point, then it must have infinitely many periodic points. This result implies that the Hamiltonian
flow on $\Sg_E$ has either exactly two or infinitely many periodic orbits.
In particular, this result holds near the equilibrium of the Hamiltonian system under
consideration. Notice that the case where $\Sigma_E$ has precisely two periodic orbits
is exceptional. Indeed, if the local Poincar\'e  map to the flow near the
orbits $\gamma_1$ and $\gamma_2$  verifies some ``twist'' condition, then the theory of twist-maps implies the existence of infinitely many periodic orbits.
The simplest example of a Hamiltonian system with exactly two periodic orbits in each $\Sg_E$ is that of a Hamiltonian function as in equation (\ref{H1}) with $\Oc_3\equiv 0$ and $\alpha_1/\alpha_2$ equal to an irrational number.

In \cite{HMS2015}, the authors proved a theorem that can be used to check whether $\Sg_E$ contains infinitely many periodic orbits. This theorem requires only the existence of a pair of periodic orbits $\gamma_1$ and $\gamma_2$ forming a Hopf link and a sort of non-resonance condition.
This hypothesis is formulated in terms of the rotation numbers
of the periodic orbits $\gamma_1$ and $\gamma_2$. This rotation number is defined as follows.
 At first, a global trivialization of the tangent bundle of $\Sg_E$
must be found (that means three everywhere linearly independent vector-fields $V_1, V_2, V_3$  on $\Sg_E$
must be found) such that $V_3=X$. This is always possible. If $\nabla H$ is the gradient of $H$,
then the orientation of the frame $V_1,V_2, X, \nabla H$ can be chosen such that
$
\omega_0\wedge \omega_0(\nabla H,V_1,V_2,X)=2(\iota_{\nabla H}\omega_0\wedge \omega_0)(V_1,V_2,X)>0.
$
Let $\varphi_t$ and
$\varphi^{L}_t$ denote
the  flow of $X$ and the linearized flow of $X$, respectively. Then, given a point $p$ on $\gamma_1$
 we write $\varphi^L_t(V_1(p))=\ap_1(t)V_1(\varphi_t(p))+\ap_2(t)V_2(\varphi_t(p))+\ap_3(t)X(\varphi_t(p))$, where we assume that $(\alpha_1(0),\alpha_2(0))\neq (0,0)$.
Since $X$ is invariant under the flow, $\ap_1^2(t)+\ap_2^2(t)>0$ for every $t\ge 0$, and we can define a continuous argument $t\mapsto \theta(t)$, $\theta(0)=0$, satisfying $\ap_1(t) + i \ap_2(t)\in \R^+ e^{i\theta(t)}$. The rotation number $\rho_1$ associated with $\gamma_1$ is
\begin{equation}
\rho_1=T_1\lim_{t\to \infty}\frac{\theta(t)}{2\pi t},
\label{rho}
\end{equation}
where $T_1$ is the period of $\gamma_1$.
This limit always exists and does not depend on the choices involved.
The rotation number $\rho_2$ is defined in the same way. Notice that the global trivialization of $T\Sg_E$ ensures some ``coherence'' between $\rho_1$ and $\rho_2$. See Appendix \ref{sec_quaternion} for more details.

The non-resonance condition between $\gamma_1$ and $\gamma_2$ that appears in \cite{HMS2015} is
\begin{equation}
(\rho_1-1)(\rho_2-1)\ne 1.
\label{twist}
\end{equation}
If this condition is verified, then infinitely many periodic orbits exist in $\Sg_E$.
We remark that the periodic orbits $\gamma_1$ and $\gamma_2$ do
not necessarily need to be those found by Hofer, Wysocki, and Zehnder.
It is enough that they form a Hopf link. The main theorem in
\cite{HMS2015} also explains how these infinitely many periodic orbits are linked
to $\gamma_1$ and $\gamma_2$.

If $\Oc_3$ vanishes identically, then $\gamma_1$ and $\gamma_2$ are respectively given by $\Sigma_E \cap \{x_2=y_2=0\}$ and $\Sigma_E \cap \{x_1=y_1=0\}$. Their rotation numbers
are $\rho_1=1+\alpha_2/\alpha_1$ and $\rho_2=1+\alpha_1/\alpha_2$. Therefore, the Hopf link $\gamma_1\cup \gamma_2$ is resonant in this case, and non-resonance is observed only if $\Oc_3$ does not vanish identically. We shall read the non-resonance condition from the Taylor expansion of $H$ in suitable coordinates near the origin.

\subsection{Main results} Let $H=H_2 + \mathcal{O}_3$ be the Hamiltonian \eqref{H1}. In complex variables $z_j=x_j+iy_j$, $j=1,2,$ the order-$2$ term of $H$ becomes $H_2 = \frac{\alpha_1}{2}|z_1|^2 + \frac{\alpha_2}{2}|z_2|^2$. We may assume that $0<\alpha_1 \leq \alpha_2$ and we say that $H$ is resonant at the origin if there exists a pair of relatively prime integers $(m_1,m_2)\in \Z^2$ such that $\alpha_1 m_1 + \alpha_2 m_2 =0.$ In that case, we may assume that $m_1<0$, and thus $|m_1|\geq m_2>0$.  If such a pair $(m_1,m_2)$ does not exist, we say that $H$ is non-resonant at the origin. In this case, we consider $|m_1|=m_2=+\infty.$ We say that $H$ is weakly non-resonant at the origin if either $H$ is non-resonant at the origin, or $H$ is resonant at $0$ and $m_2\geq 2$. This last condition is equivalent to $\alpha_2/\alpha_1 \not \in \N$.

For each $N\geq 3$, we shall consider the Birkhoff-Gustavson normal form
$$
H = H_N + \mathcal{O}_{N+1}= H_2 + \Gamma^{(3)} + \ldots +\Gamma^{(N)} +\mathcal{O}_{N+1},
$$
where $\Gamma^{(s)},s=3,\ldots,N,$ is a $s$-homogeneous polynomial in variables $z_1,z_2,\bar z_1,\bar z_2$ satisfying
$$
D\cdot \Gamma^{(s)}:=\frac{1}{i} \sum_{j=1,2}\alpha_j(z_j \partial_{z_j} - \bar z_j \partial_{\bar z_j}) \cdot \Gamma^{(s)}=0.
$$
This means that
$$
H_N = \frac{\alpha_1}{2}|z_1|^2 + \frac{\alpha_2}{2}|z_2|^2 + \sum_{(k,l)\in \mathcal{M}_3^N} a_{kl} z^k \bar z^l,
$$
where $\mathcal{M}_3^N$ is the set of resonant exponents $(k,l)$ given in \eqref{M3N}, and $z^k\bar z^l=z_1^{k_1}z_2^{k_2}\bar z_1^{l_1}\bar z_2^{l_2},$ $k=(k_1,k_2), l=(l_1,l_2).$ Here, the coefficients $a_{kl}=a_{k_1,k_2,l_1,l_2}\in \C$ satisfy $a_{lk} = \bar a_{kl}, \forall (k,l)\in \mathcal{M}_3^N$, and are determined by the Taylor series of $H$ at the origin, as explained in section \ref{normal1}.

We wish to find a non-resonant Hopf link for small energies of $H$. To do that, we consider some special coefficients of $H$ in normal form. For fixed $N\in \N$, let $$2\leq \nu \leq \lfloor N/2 \rfloor$$ be the smallest integer so that at least one of the coefficients
$$
a_{\nu,0,\nu,0}, \quad  a_{0,\nu,0,\nu}, \quad a_{\nu-1,1\nu-1,1}, \quad a_{1,\nu-1,1,\nu-1},
$$ is non-zero. We assume the existence of $\nu$. Let
$$
\Omega_{\nu,1}:=a_{\nu-1,1,\nu-1,1}-\nu a_{\nu,0,\nu,0}\frac{\alpha_2}{\alpha_1} \quad \mbox{ and } \quad \Omega_{\nu,2}:= a_{1,\nu-1,1,\nu-1}-\nu a_{0,\nu,0,\nu}\frac{\alpha_1}{\alpha_2}.
$$
Finally, let
$$
\Omega_\nu:=\frac{\Omega_{\nu,1}}{\alpha_2\alpha_1^{\nu-1}}+\frac{\Omega_{\nu,2}}{\alpha_1\alpha_2^{\nu-1}}.
$$
We provide conditions for $H_N$ to admit periodic orbits $\gamma_1 \subset \{z_2=0\}$ and $\gamma_2 \subset \{z_1=0\}$ forming a Hopf link. In this case, $\Omega_{\nu,1}, \Omega_{\nu,2}$ and $\Omega_\nu$ quantify the non-resonant condition in the lowest order. Our first result deals with Hamiltonians, which are weakly non-resonant at the origin.

\begin{thm}\label{thm: main theorem}
Assume that $H=H_2+\mathcal{O}_3$ is weakly non-resonant at the origin, i.e., $\alpha_2$ is not a multiple of $\alpha_1$. Assume further that for some $N\geq 3$, the Birkhoff-Gustavson normal form $H= H_{N}+ \mathcal{O}_{N+1}$  satisfies one of the following conditions:
\begin{itemize}
\item[(i)] $2<m_2\leq +\infty$ and $\Omega_\nu \neq 0;$
\item[(ii)] $m_2=2$, $|m_1|>2(\nu-1)$ and $\Omega_{\nu,1},\Omega_\nu\neq 0;$
\item[(iii)] $m_2=2$, $|m_1|=2(\nu-1)$, $\Omega_{\nu,1},\Omega_\nu\neq 0$ and $a_{0,2,|m_1|,0}=0$;
\item[(iv)] $m_2=2$, $|m_1|<2(\nu-1)$ and $\Omega_{\nu,2},a_{0,2,|m_1|,0}\neq 0.$
\end{itemize}
Then, for every $E>0$ sufficiently small, the sphere-like component $\Sigma_E \subset H^{-1}(E)$ carries a pair of periodic orbits $\gamma_{1,E},\gamma_{2,E}$ forming a non-resonant Hopf link. In particular, their respective rotation numbers $\rho_{1,E},\rho_{2,E}$ satisfy $(\rho_{1,E}-1)(\rho_{2,E}-1)\neq 1$ and $\Sigma_E$ contains infinitely many periodic orbits.
\end{thm}

For resonant Hamiltonians at the origin, which are not weakly non-resonant (i.e., $\alpha_2$ is a multiple of $\alpha_1$), we may need additional symmetries to prove the existence of a non-resonant Hopf link. We first assume that $|m_1|>m_2 =1$, i.e., $\alpha_2$ is a non-trivial multiple of $\alpha_1$. In this case, we also assume that $\partial_{y_2}H(y_1,0,x_1,0)=\partial_{x_2}H(y_1,0,x_1,0)=0
$ for every $(y_1,x_1)$ sufficiently close to $(0,0)$. This assumption automatically gives a family of periodic orbits on the plane $x_2=y_2=0$, which will be part of the desired non-resonant Hopf link on the corresponding energy surface. The Birkhoff-Gustavson normal forms constructed in section \ref{normal1} share the same symmetry.

The following theorem provides conditions for the existence of a non-resonant Hopf link in the case $\alpha_2$ is a non-trivial multiple of $\alpha_1$, i.e., $|m_1|>m_2=1$.

\begin{thm}\label{thm: main theorem 2}
Assume that $H=H_2+\mathcal{O}_3$ is such that $\alpha_2$ is a non-trivial multiple of $\alpha_1$, and that $\partial_{y_2}H(y_1,0,x_1,0)=\partial_{x_2}H(y_1,0,x_1,0)=0
$ for every $(y_1,x_1)$ sufficiently close to $(0,0)$. Assume further that for some $N\geq 3$, the Birkhoff-Gustavson normal form $H= H_{N}+\mathcal{O}_{N+1}$  satisfies one of the following conditions:
\begin{itemize}
\item[(i)] $|m_1|> 2$, $|m_1|> \nu -1$ and $\Omega_{\nu,1},\Omega_\nu\neq 0$.
\item[(ii)]  $|m_1|> 2$, $|m_1|=\nu -1$, $\Omega_{\nu,1},\Omega_{\nu}\neq 0$ and $a_{0,2,2|m_1|,0}=0$.
\item[(iii)]  $|m_1|> 2$, $|m_1|<\nu -1$, $\Omega_{\nu,2}\neq 0$ and $a_{0,2,2|m_1|,0}\neq 0$.
\item[(iv)] $|m_1|=2,$ $\nu=2$, $\Omega_{2,1}\neq 0$ and $a_{0,1,2,0} \neq 0.$
\item[(v)]  $|m_1|= 2$, $\nu=2$, $\Omega_{2,1},\Omega_2\neq 0$.
\item[(vi)]  $|m_1|= 2$, $\nu=3$, $\Omega_{3,1},a_{0,1,2,0}\neq 0$ and $a_{0,2,4,0}=0$.
\end{itemize}
Then the sphere-like component $\Sigma_E\subset H^{-1}(E)$ carries a pair of periodic orbits $\gamma_{1,E}, \gamma_{2,E}$ forming a non-resonant Hopf-link for every $E>0$ sufficiently small. In particular, their rotation numbers satisfy $(\rho_{1,E}-1)(\rho_{2,E}-1)\neq 1$ and $\Sigma_E$ contains infinitely many periodic orbits.
\end{thm}

Theorems \ref{thm: main theorem} and \ref{thm: main theorem 2} apply to classical Hamiltonian systems in Celestial Mechanics, such as the Spatial Isosceles Three-Body Problem, see section \ref{SI3BP}.

Finally, we consider the case $\alpha_1=\alpha_2$ ($|m_1|=m_2=1$). The Hamiltonian is strongly resonant at the origin, and the symmetry we consider is the following. Let $p\geq 3$ be an integer and assume that $H$ is invariant under the $\Z_p$-action on $\R^4$ generated by the unitary map $R:(y,x)\mapsto(e^{2\pi i/p}y,e^{2\pi i/p}x).$ The sphere-like component $\Sigma_E\subset H^{-1}(E)$ quotients down to a lens space $\Sigma_E/\Z_p \equiv L(p,p-1).$

For each $N\in \N$, we shall construct an $R$-invariant symplectic change of coordinates so that both the Birkhoff-Gustavson normal form $H=H_N+\mathcal{O}_{N+1}$ and $H_N$ are $R$-invariant.
Consider the canonical transformation
\begin{equation}\label{Def_Psi}
\Psi:(y_1,y_2,x_1,x_2)\mapsto 2^{-1/2}(y_1+y_2,x_1-x_2,x_1+x_2,y_2-y_1).
\end{equation}
Then $H\circ \Psi = H_N \circ \Psi + \mathcal{O}_{N+1}$ is in normal form up to order $N$.

The following theorem has applications in Hill's Lunar Problem and the H\'enon-Heiles System; see sections \ref{Hill} and \ref{Henon-Heiles}.

\begin{thm}\label{thm: main theorem 3}
Assume that $H=H_2 + \mathcal{O}_3$ satisfies $\alpha_1 = \alpha_2$, and admits the $\Z_p$-symmetry as above for some $p\geq 3$. Assume further that for some $N\geq 3$, the Birkhoff-Gustavson normal form $H\circ \Psi = H_N \circ \Psi + \mathcal{O}_{N+1}$  satisfies one of the following conditions:
\begin{itemize}
\item[(i)]  $\nu=2$, $\Omega_{\nu}\neq 0$ and $a_{0,2,2,0}=0$;
\item[(ii)] $\nu=2$, $\Omega_{2,1}=-\Omega_{2,2}\neq 0, a_{0,2,2,0}=0$ and $\beta_1+\beta_2+2\Omega_{2,1}\Omega_{2,2}\neq 0$,
where
$$
\begin{aligned}
\beta_1&=6 (2a_{2,0,2,0} - a_{1,1,1,1})a_{2,0,2,0} + \alpha_1(a_{2,1,2,1} - 3 a_{3,0,3,0}),\\
\beta_2&=6 (2a_{0,2,0,2} -
a_{1,1,1,1})a_{0,2,0,2} +\alpha_1 (a_{1,2,1,2} - 3 a_{0,3,0,3}).
\end{aligned}
$$
\end{itemize}
For every $E>0$ sufficiently small, the sphere-like component $\Sigma_E\subset H^{-1}(E)$ carries a pair of periodic orbits $\gamma_{1,E}, \gamma_{2,E}$ forming a non-resonant Hopf-link. In particular, their rotation numbers satisfy $(\rho_{1,E}-1)(\rho_{2,E}-1)\neq 1$ and $\Sigma_E$ contains infinitely many periodic orbits.
\end{thm}

\begin{rem}
    The Birkhoff-Gustavson normal forms $H=H_N + \mathcal{O}_{N+1}, N\geq 3,$ in Theorems \ref{thm: main theorem}, \ref{thm: main theorem 2} and \ref{thm: main theorem 3} are canonically constructed in section \ref{normal1} via generating functions. Each $H_N$ is obtained from $H_{N-1}$ by adding a normalized $N$-homogeneous polynomial in the kernel of $D$.
\end{rem}


The paper is organized as follows. In section \ref{normal1}, we introduce the Birkhoff-Gustavson normal form $H = H_N + \mathcal{O}_{N+1}$ for each $N\geq 3.$ Symmetric versions of the normal form are also established. In section \ref{sec normal form}, we study the existence of a non-resonant Hopf link for $H_N$ under some general conditions, including various cases of resonant equilibrium points. In section \ref{sec 4}, we establish the passage from a non-resonant Hopf link for $H_N$ to a non-resonant Hopf link for $H$. Theorems \ref{thm: main theorem}, \ref{thm: main theorem 2} and \ref{thm: main theorem 3} are proved in sections \ref{sec: weakly nonresonant Hamiltonians}, \ref{sec: Hamiltonians with symmetries} and \ref{sec: zp symmetry}, respectively. In section \ref{SI3BP}, we apply our results to the Spatial Isosceles Three-Body Problem; see Theorem \ref{thm: isosceles 3bp}. Finally, we deal with Hill's Lunar Problem and the H\'enon-Heiles potential in sections \ref{Hill} and \ref{Henon-Heiles}, see Theorems \ref{thm: lunar problem} and \ref{thm: Henon-Heiles}, respectively.

\section{Birkhoff-Gustavson normal forms}\label{normal1}

Let $H$ be a Hamiltonian function as in equation (\ref{H1}) with $\alpha_1,\alpha_2>0$. This Hamiltonian is said to be resonant at the origin if there exists $0\neq m=(m_1,m_2)\in\Z^2$ such that
$\alpha \cdot m=\alpha_1m_1+\alpha_2m_2=0$, where $\alpha=(\alpha_1,\alpha_2)$.
Also, $H$ is said to be non-resonant at the origin if no such $m$ exists, namely $\alpha \cdot m=0$ if and only if
$m=(0,0)$.
We can always arrange $\alpha_1$ and $\alpha_2$ such that the following relation is verified
$$
0<\alpha_1\le\alpha_2.
$$
If  $H$ is resonant we  can always choose $(m_1,m_2)\in \Z^2$ such that
\begin{equation}\label{generator}
\gcd(m_1,m_2)=1, \qquad m_1<0, \qquad \mbox{ and } \qquad  |m_1|\ge m_2.
\end{equation}
In the non-resonant case, we take $m=(0,0)$.
The vector $\alpha=(\alpha_1,\alpha_2)$ defines a module $\mathcal{M}$ over $\Z$ given by
\begin{equation}
\mathcal{M}:=\{l=(l_1,l_2)\in \Z^2:l \cdot \alpha=0\}=\Z\cdot (m_1,m_2).
\label{M}
\end{equation}

A result due to Birkhoff \cite{Birkhoff1927} in the non-resonant case and due to Gustavson in the resonant case \cite{Gustavson1966} (also due to Moser \cite{moser58}, \cite{moser68}) state the following: given any integer $N>2$, there exists a canonical change of coordinates
$$
\Phi:(\eta_1,\eta_2,\xi_1,\xi_2) \mapsto (y_1,y_2,x_1,x_2),
$$
such that in the new coordinates $(\eta,\xi)$, the Hamiltonian
function \eqref{H1}  becomes
$$
H\circ \Phi=H_N+\mathcal{O}_{N+1},
$$
where $H_N$ is an order-$N$ polynomial in normal form, as we explain below, and $\mathcal{O}_{N+1}$ is a smooth function whose partial derivatives up to order $N$ vanish at the origin. The quadratic part $H_2$ of $H_N$ coincides with the quadratic part of $H$ as in \eqref{H1},
namely
$$
H_2=\frac{\alpha_1}{2}(\eta_1^2+\xi_1^2)+\frac{\alpha_2}{2}(\eta_2^2+\xi_2^2).
$$
Moreover, consider the partial differential operator
 $$
 D:= \sum_{j=1,2}\alpha_j(\eta_j \partial_{\xi_j} - \xi_j \partial_{\eta_j}).
 $$
Notice that $D$ is the differentiation along the Hamiltonian
vector field associated with $H_2$.
Then, $H_N$ is characterized by the equivalent properties
$$
D\cdot H_N=0 \Leftrightarrow \{H_2,H_N\}=0,
$$
where $\{\cdot,\cdot\}$ is the Poisson bracket associated with the canonical symplectic form.

For computational reasons, it is convenient to write $H_N$ in complex coordinates $$z_j:=\xi_j+i\eta_j, \qquad \bar z_j:=\xi_j - i \eta_j, \qquad j=1,2.$$
This definition implies that
$$
H_2 = \frac{\alpha_1}{2}|z_1|^2 + \frac{\alpha_2}{2}|z_2|^2 \qquad \mbox{ and } \qquad D=\frac{1}{i}\sum_{j=1,2}\alpha_j(z_j\partial_{z_j}-\bar z_j \partial_{\bar z_j}),
$$
where $\partial_{z_j}=\frac{1}{2}(\partial_{\xi_j}-i\partial_{\eta_j})$. The advantage of using complex coordinates is that, given pairs of non-negative integers
$k=(k_1,k_2)\in \Z_+^2$ and $l=(l_1,l_2)\in \Z_+^2$, the monomials
$$
z^k \bar z^l:=z_1^{k_1}z_2^{k_2}\bar z_1^{l_1}\bar z_2^{l_2}
$$
are eigenfunctions of $D$, namely
\begin{equation}\label{D on monomials}
\begin{aligned}
D\cdot  z^k\bar z^l & =
\frac{1}{i}\left[\alpha_1(k_1-l_1)+\alpha_2(k_2-l_2)\right]z^k\bar z^l\\
& =
\frac{1}{i}\left[\alpha \cdot (k-l) \right]z^k\bar z^l.
\end{aligned}
\end{equation}
This implies that $D\cdot z^k\bar z^l =0$ if and only if
$(k-l)\in \mathcal{M}$, where $\mathcal{M}$ is the $\Z$-module defined in \eqref{M}.
To simplify the notation, we define
\begin{equation}
\mathcal{M}_3^N:=\{(k,l)\in \Z_+^2\times \Z_+^2:(k-l)\in \mathcal{M}, \ 3\leq k_1+k_2+l_1+l_2\leq N\},
\label{M3N}
\end{equation}
and write
$$
H_N=\frac{\alpha_1}{2}|z_1|^2+\frac{\alpha_2}{2}|z_2|^2+
\sum_{(k,l)\in \mathcal{M}_3^N} a_{kl}z^k\bar z^l,
$$
where the coefficients $a_{kl}:=a_{k_1,k_2,l_1,l_2}\in\C$ are obtained from the order-$N$ Taylor expansion
of $H$ around the origin, as we shall explain below. The fact that $H$ is a real-valued function implies that
$a_{kl}=\bar a_{lk}.$
If $H$ is non-resonant, i.e.,  $\mathcal{M}=\{0\}$,
 then the coefficients $a_{kl}$ are uniquely
determined and all terms of $H_N$ have the form $a_{k_1,k_2,k_1,k_2}|z_1|^{2k_1}|z_2|^{2k_2}, k\geq 1$. If $H$ is resonant, i.e., $\mathcal{M}\neq \{0\}$,  then the coefficients $a_{kl}$ are not uniquely
determined (there exist canonical changes of variables that preserve the form of $H$ but change
the coefficients $a_{kl}$, see \cite{moser68}). In the resonant case, it is convenient to
rewrite $H_N$ differently. In this case, the module $\mathcal{M}$ is generated by the
vector $m=(m_1,m_2)\neq 0$, as shown in equations \eqref{generator} and \eqref{M}. We define the special monomial associated with $m$
\begin{equation}
\sigma(z_1,z_2,\bar z_1,\bar z_2):= z_2^{m_2}\bar z_1^{|m_1|}\qquad \Rightarrow \qquad \bar\sigma(z_1,z_2,\bar z_1,\bar z_2) =  z_1^{|m_1|}\bar z_2^{m_2}.
\label{sg}
\end{equation}
Notice that $\sigma=z^k\bar z^l$ for $k=(0,m_2)$ and $l=(|m_1|,0)$. It satisfies $D\cdot \sigma=0$.

We write the formal sum $\sum_{k-l\in \mathcal{M}} a_{kl}z^k\bar z^l$ as
$$
\begin{aligned}
\sum_{(k,l)\in \mathcal{M}} a_{kl}z^k\bar z^l &  = \sum_{(k,l)\in \mathcal{M}:k-l=0} a_{kl}z^k\bar z^l\\
& +\sum_{n\ge 1} \left[ \sum_{(k,l)\in \mathcal{M}:k-l=nm} a_{kl}z^k\bar z^l+\sum_{(k,l)\in \mathcal{M}:k-l=-nm} a_{kl}z^k\bar z^l\right].
\end{aligned}
$$
Using the notation $|z|^{2k}:=|z_1|^{2k_1}|z_2|^{2k_2}$, we write
$$
\begin{aligned}
\sum_{(k,l)\in \mathcal{M}:k-l=0} a_{kl}z^k\bar z^l& = \sum_{k\in \Z^2_+} a_{kk}|z|^{2k},\\
\sum_{(k,l)\in \mathcal{M}:k-l=nm} a_{kl}z^k\bar z^l & = \sg^n\sum_{k\in \Z^2_+} a_{k_1,k_2+nm_2,k_1+n|m_1|,k_2}|z|^{2k},\\
\sum_{(k,l)\in \mathcal{M}:k-l=-nm} a_{kl}z^k\bar z^l&= \bar \sg^n\sum_{k\in \Z^2_+}
\underbrace{a_{k_1+n|m_1|,k_2,k_1, k_2+nm_2}}_{=\bar a_{k_1,k_2+nm_2,k_1+n|m_1|,k_2}}|z|^{2k}.
\end{aligned}
$$
Notice that we have used $(k_1,l_2)\in \Z^2_+$ and $(l_1,k_2)\in \Z^2_+$ to parametrize the right-hand side sum in the second and third equations above, respectively, denoting them again by $(k_1,k_2)\in \Z^2_+$.  Taking  into account that $3\leq k_1+k_2+l_1+l_2\leq N$,  we define
\begin{equation}\label{An}
\begin{aligned}
A_0(|z_1|^2,|z_2|^2)& : = \sum_{k\in \Z^2_+:3\le 2k_1+2k_2\le N } a_{kk}|z_1|^{2k_1}|z_2|^{2k_2},\\
A_n(|z_1|^2,|z_2|^2)& : =
\sum_{k\in \Z^2_+:3\le 2k_1+2k_2+n(|m_1|+m_2)\le N } a_{k_1,k_2+nm_2,k_1+n|m_1|,k_2}|z_1|^{2k_1}|z_2|^{2k_2},\\
\bar A_n(|z_1|^2,|z_2|^2)& :=
\sum_{k\in \Z^2_+:3\le 2k_1+2k_2+n(|m_1|+m_2)\le N } \bar a_{k_1,k_2+nm_2,k_1+n|m_1|,k_2}|z_1|^{2k_1}|z_2|^{2k_2}.
\end{aligned}
\end{equation}
Finally,  we find an expression for $H_N$
\begin{equation}\label{hatH}
\begin{aligned}
H_N & = \frac{\alpha_1}{2}|z_1|^2+\frac{\alpha_2}{2}|z_2|^2+A_0(|z_1|^2,|z_2|^2) \\
& +\sum_{1\leq n\leq \lfloor N/(|m_1|+m_2) \rfloor}\left[\sg^n A_n(|z_1|^2,|z_2|^2)+
       \bar \sg^n\bar A_n(|z_1|^2,|z_2|^2)\right],
\end{aligned}
\end{equation}
where $\lfloor a \rfloor$ denotes the integer part of $a\in\mathbb{R}$. In the non-resonant case, defining $|m_1|=m_2=\infty$ is convenient. In this case, the sum over the index $n$ in the above equation vanishes for every $N$.
In this way, the expression \eqref{hatH}
also applies to the non-resonant case. We refer to $H_N$ as resonant (non-resonant) if $H_N$ is resonant (non-resonant) at $0$.

\subsection{The construction of the normal form}\label{sec: construction}
We explain the construction of the Birkhoff-Gustavson normal form introduced in the previous section. Let $H = H_2 + \mathcal{O}_3,$ where we assume that $H_2$ is already in normal form, that is $H_2 = \frac{\alpha_1}{2} (y_1^2+x_1^2) + \frac{\alpha_2}{2}(y_2^2+x_2^2)$. In particular, $D\cdot H_2=0$. Following Gustavson \cite{Gustavson1966} and Moser \cite{moser68}, we aim to find for each $N\geq 3$ a canonical change of coordinates $\Phi^{(N)}:(\eta_1, \eta_2, \xi_1,\xi_2)\mapsto (y_1,y_2,x_1,x_2)$ near $0$ so that the new Hamiltonian $H\circ \Phi^{(N)}$ has the form $\Gamma = \Gamma^{(2)} + \Gamma^{(3)} + \ldots + \Gamma^{(N)} +O_{N+1},$ where $\Gamma^{(2)}$ has the same form as $H_2$ after replacing $(y_1,y_2,x_1,x_2)$ with $(\eta_1,\eta_2,\xi_1,\xi_2)$, and $\Gamma^{(i)}$ is a $i$-homogeneous polynomial satisfying $D \cdot \Gamma^{(i)}=0$ for every $2\leq i\leq N$.

We consider a finite sequence of coordinate changes $\Phi^{(N)} = \Phi_3 \circ \ldots \circ \Phi_N$, where $H\circ \Phi_3$ is in normal form up to order $3$, $H \circ \Phi_3 \circ \Phi_4$ is in normal form up to order $4$ and so on up to order $N$.  Before explaining the construction of each $\Phi_i$, we observe that the operator $$
D^i: \mathcal{P}^i \to \mathcal{P}^i, \quad D^i=\sum_{j=1}^2 \alpha_j(y_j\partial_{x_j} - x_j\partial_{y_j})\big|_{\mathcal{P}^i},
$$
acting on the space $\mathcal{P}^i$ of $i$-homogeneous polynomials in variables $(y_1,y_2,x_1,x_2)$,  satisfies the following crucial properties
\begin{equation}\label{property_Di}
\ker D^i \oplus {\rm im} D^i= \mathcal{P}^i \quad \mbox{ and } \quad \ker D^i \cap {\rm im}D^i = 0,
\end{equation}
where ${\rm im} D^i$ is the range of $D^i$.

Let us assume by induction that  $H$, given in coordinates $(y_1,y_2,x_1,x_2)$, is in normal form up to order $s-1$,  that is,
$$
\begin{aligned}
\hat H_{s-1}& :=H \circ \Phi_{3}\circ \ldots \circ \Phi_{s-1}\\
& = \Gamma^{(2)} + \ldots +\Gamma^{(s-1)} + \hat H_{s-1}^{(s)} +\mathcal{O}_{s+1},
\end{aligned}
$$
where $D \cdot \Gamma^{(i)}=0$ for every $2\leq i \leq s-1$, and $\hat H^{(s)}_{s-1}$ is a $s$-homogeneous polynomial. We shall construct a canonical change of coordinates $\Phi_s:(\eta_1,\eta_2,\xi_1,\xi_2)\mapsto (y_1,y_2,x_1,x_2)$ such that the new Hamiltonian in coordinates $(\eta_1,\eta_2,\xi_1,\xi_2)$ writes as
$$
\begin{aligned}
\hat H_s & :=H \circ \Phi_3 \circ \ldots \circ \Phi_s \\
& = \Gamma^{(2)} + \ldots +\Gamma^{(s)} + \mathcal{O}_{s+1},
\end{aligned}
$$
where $D \cdot \Gamma^{(i)}=0,$ for every $2\leq i\leq s$. Here, $D=\sum_{j=1}^2\alpha_j(\eta_j\partial_{\xi_j}-\xi_j\partial_{\eta_j})$, and the functions $\Gamma^{(i)}(\eta,\xi),\ 2\leq i\leq s-1,$ coincide with $\Gamma^{(i)}(y,x)$ by formally replacing $(y,x)$ with $(\eta,\xi)$.

To find $\Phi_s$, we search for a generating function of the form $W_s(\eta_1,\eta_2,x_1,x_2)=x_1\eta_1+x_2\eta_2 + G_s(\eta_1,\eta_2,x_1,x_2)$ that satisfies $\xi d\eta+ydx = dW_s$. This leads to the relations
\begin{equation}\label{equ: relation of Ws}
\xi = \partial_\eta W_s= x+ \partial_\eta G_s, \quad y = \partial_x W_s = \eta + \partial_x G_s.
\end{equation}

We choose $G_s$ as a $s$-homogeneous polynomial in ${\rm im}D^s \subset \mathcal{P}^s$. Notice that $W_s$ induces a canonical transformation $\Phi_s:(\eta_1,\eta_2,\xi_1,\xi_2)\mapsto(y_1,y_2,x_1,x_2)$ satisfying $d\Phi_s(0)=\mathrm{Id}$ and
\begin{equation}\label{equation of W}
\hat H_{s-1}(\partial_x W_s,x)=\hat H_s(\eta,\partial_\eta W_s),\quad \mbox{ where} \quad \hat H_s:=\hat H_{s-1}\circ \Phi_s.
\end{equation}
In view of \eqref{property_Di}, we can write $\hat H_{s-1}^{(s)} = K + L,$ where $K\in \ker D^s$ and $L\in {\rm im} D^s.$ We choose $G_s(\eta,x)$ such that $D^s \cdot G_s = -L$, where we formally replace $(y,x)$ with $(\eta,x)$ and rewrite $D^s=\sum_{j=1}^2\alpha_j(\eta_j\partial_{x_j}-x_j\partial_{\eta_j})$. In this way, we obtain
$$
\hat H^{(s)}_{s-1}\circ\Phi_s(\eta,\xi)=-D^s\cdot G_s(\eta,\xi)+K(\eta,\xi)+\mathcal{O}_{s+1}.
$$
Here, the functions $-D^s \cdot G_s(\eta,\xi)$ and $K(\eta,\xi)$ are $s$-homogeneous polynomials, formally obtained from the respective functions on variables $(\eta,x)$ after replacing $x$ with $\xi$.

A direct computation gives
$$
\begin{aligned}
\Gamma^{(2)}\circ\Phi_s=\Gamma^{(2)}+D^s\cdot G_s+\mathcal{O}_{s+1},\quad
\Gamma^{(i)}\circ\Phi_s=\Gamma^{(i)}+\mathcal{O}_{s+1},\quad  i=3,\ldots, s-1.
\end{aligned}
$$
Hence, we obtain
$$
\hat H_s:=\hat H_{s-1}\circ \Phi_s=\Gamma^{(2)}+\cdots +\Gamma^{(s-1)}+\Gamma^{(s)}+\mathcal{O}_{s+1},
$$
where $\Gamma^{(s)}:=K$ as desired.

\begin{rem} Notice that in constructing the polynomials $H_N, N\geq 3,$ in normal form as above, we have used in each step a canonical transformation $\Phi=\Phi^{(N)}$ satisfying $d\Phi(0) = {\rm Id}$ and canonically determined by the splitting \eqref{property_Di} and the relations \eqref{equ: relation of Ws}. Hence, once $H_2$ is in normal form and $G_s$ is chosen in ${\rm im}D^s$, the procedure above gives a unique sequence of Birkhoff-Gustavson normal forms $H=H_N + \mathcal{O}_{N+1}, N\geq 3$, where $H_{N+1} = H_N + \Gamma^{(N+1)}$ and $\Gamma^{(N+1)}\in \ker D$.
\end{rem}

The discussion above is summarized in the following theorem.

\begin{thm}[Birkhoff-Gustavson Normal Form]\label{thm_BK_normal_form} Let the smooth Hamiltonian $H=\frac{\alpha_1}{2}(y_1^2+x_1^2) + \frac{\alpha_2}{2}(y_2^2+x_2^2)+ \mathcal{O}_3$ be defined near the origin. For every $N\geq 3$, there exists a canonical change of coordinates $\Phi^{(N)}:(\eta,\xi) \mapsto (y,x)$, defined near the origin, so that $H \circ \Phi^{(N)} = H_N + \mathcal{O}_{N+1},$ where $H_N$ is a $N$-order polynomial satisfying $D \cdot H_N =0$. Moreover, one can arrange the sequence $\Phi^{(N)}, N\geq 3,$ so that $d\Phi^{(N)}(0) = \mathrm{Id}$ and $H_{N+1}=H_N+\Gamma^{(N+1)}$, where $\Gamma^{(N+1)}$ is a $(N+1)$-homogeneous polynomial in $(\eta,\xi)$ satisfying $D \cdot \Gamma^{(N+1)}=0$ and depending only on the $(N+1)$-homogeneous part of $H\circ \Phi^{(N)}$ at the origin.
\end{thm}

In the following, we discuss particular Hamiltonians that satisfy certain symmetries. We show that the associated normal forms $H_N$ constructed as in Theorem \ref{thm_BK_normal_form} share the same symmetry.

\begin{prop}\label{lem: normal form 2}
Let $H$ be as in Theorem \ref{thm_BK_normal_form}. Assume further that
\begin{equation}\label{cond for y2,x2}
\partial_{y_2}H(y_1,0,x_1,0) = \partial_{x_2}H(y_1,0,x_1,0)=0,
\end{equation}
for every $(y_1,x_1)$ near $(0,0)$. Then, for every $N\geq 3$, the Birkhoff-Gustavson normal form $H=H_N + \mathcal{O}_{N+1}$ given in Theorem \ref{thm_BK_normal_form} satisfies
\begin{equation}\label{cond for eta2,xi2}
\partial_{\eta_2}H_N(\eta_1,0,\xi_1,0)=\partial_{\xi_2}H_N(\eta_1,0,\xi_1,0)=0,
\end{equation}
for every $(\eta_1,\xi_1)$ near $(0,0)$.
\end{prop}

\begin{proof}
Following the proof of Theorem \ref{thm_BK_normal_form}, consider the generating function $W_3(\eta,x) = x_1\eta_1 + x_2 \eta_2 + G_3$ and the induced canonical transformation $\Phi_3$. From the identity $H(\eta+\partial_xG_3,\xi-\partial_\eta G_3)=\hat H_3(\eta,\xi),$  we obtain the Hamiltonian
\begin{equation}\label{expending0}
\begin{aligned}
\hat H_3(\eta,\xi) =\Gamma^{(2)}(\eta,\xi)+\Gamma^{(3)}(\eta,\xi)+\mathcal{O}_4.
\end{aligned}
\end{equation}
Here, $-D\cdot G_3+\Gamma^{(3)}=H^{(3)}$ and $D\cdot \Gamma^{(3)}=0$, where $H^{(3)}$ is the $3$-homogeneous part of $H$.
Since $H$ satisfies \eqref{cond for y2,x2}, its $3$-homogeneous part $H^{(3)}$ also satisfies \eqref{cond for y2,x2}. In complex variables $z_j=x_j+iy_j, j=1,2$, we have $\partial_{\bar z_2}H^{(3)}(z_1,0)=0$ for every $z_1$ sufficiently close to $0$.

In view of \eqref{D on monomials}, the monomials $Q=c z^k\bar z^l$ are eigenfunctions of $D$. Hence, if a monomial $Q$ of $H^{(3)}$ satisfies $D\cdot Q =0$, then $Q$ must be part of $\Gamma^{(3)}$. In the same way, if a monomial $Q$ of $H^{(3)}$ satisfies $D \cdot Q = \beta Q$, where $\beta=-i\alpha\cdot(k-l) \neq0$, then $G_3$ must contain $-Q/\beta$ so that $-D (-Q/\beta) = Q$. We conclude that $\Gamma^{(3)}$ coincides with the part of $H^{(3)}$ containing the monomials $Q$ satisfying $D \cdot Q=0$, and $G_3$ contains precisely the same monomials $Q$ of $H^{(3)}$ satisfying $D\cdot Q \neq 0$, perhaps with distinct coefficients.   Since $\partial_{z_2} H^{(3)}(z_1,0)=0$ for every $z_1$ sufficiently close to $0$, $H^{(3)}$ does not have monomials of the form $\beta z_1^{k_1}\bar z_1^{l_1} z_2$ or $\beta z_1^{k_1}\bar z_1^{l_1} \bar z_2$ where $k_1,l_1\in \{0,1,2\}$ satisfy $k_1+l_1=2$. Hence, both $G_3$ and $\Gamma^{(3)}$ do not have such monomials, and thus   $\partial_{z_2}G_3(z_1,0)=\partial_{z_2}\Gamma^{(3)}(z_1,0)=0$ for every $z_1$ sufficiently close to $0$. Here, the complex coordinates are formally replaced with $z_j = x_j+i \eta_j$ and $z_j=\xi_j + i \eta_j, j=1,2$.  Since $G_3$ and $\Gamma^{(3)}$ are real-valued functions, this implies that $\partial_{x_2}G_3(\eta_1,0,x_1,0)=\partial_{\eta_2}G_3(\eta_1,0,x_1,0) = \partial_{\xi_2}\Gamma^{(3)}(\eta_1,0,\xi_1,0)=\partial_{\eta_2}\Gamma^{(3)}(\eta_1,0,\xi_1,0)=0$ for every $(\eta_1,x_1)$ and $(\eta_1,\xi_1)$ sufficiently close to $(0,0)$. In particular, $G_3$ and $\Gamma^{(3)}$ satisfy \eqref{cond for y2,x2} in the corresponding coordinates.

Recall that $\Phi_3$ is a local diffeomorphism induced by the relations $y_i=\eta_i+\partial_{x_i}G_3$ and $x_i=\xi_i-\partial_{\eta_i}G_3, i=1,2,$ where $G_3$ is seen as a function of $(\eta,x)$ after formally replacing $(y,x)$ with $(\eta,x)$. It follows that if $(\eta_2,x_2)=(0,0)$, then $(\eta_2,\xi_2)=(y_2,x_2)=(0,0)$. In particular, $\Phi_3$ sends the coordinate plane $\{\eta_2=\xi_2=0\}$ to $\{y_2=x_2=0\}$.

From $\hat H_3(\eta,\xi) = H(\eta+\partial_xG_3,\xi-\partial_\eta G_3)$, we obtain
\begin{equation}\label{equ: d of H3}
\begin{aligned}
\partial_{\eta_2}\hat H_3&=\partial_{y_1} H \partial_{\eta_2} y_1+\partial_{y_2} H \partial_{\eta_2} y_2+\partial_{x_1} H \partial_{\eta_2} x_1+\partial_{x_2} H \partial_{\eta_2} x_2,\\
\partial_{\xi_2}\hat H_3&=\partial_{y_1} H \partial_{\xi_2} y_1+\partial_{y_2} H \partial_{\xi_2} y_2+\partial_{x_1} H \partial_{\xi_2} x_1+\partial_{x_2} H \partial_{\xi_2} x_2.
\end{aligned}
\end{equation}
We also have
$$
\begin{aligned}
\partial_{\eta_2} y_1&= \partial^2_{x_1\eta_2}G_3+\partial^2_{x_1x_1}G_3\partial_{\eta_2}x_1+\partial^2_{x_1x_2}G_3\partial_{\eta_2}x_2\\
\partial_{\xi_2} y_1&=\partial^2_{x_1x_1}G_3\partial_{\xi_2}x_1+\partial^2_{x_1x_2}G_3\partial_{\xi_2}x_2\\
\partial_{\eta_2} x_1&= -\partial^2_{\eta_1\eta_2}G_3-\partial^2_{\eta_1x_1}G_3\partial_{\eta_2}x_1-\partial^2_{\eta_1x_2}G_3\partial_{\eta_2}x_2\\
\partial_{\xi_2} x_1&=-\partial^2_{\eta_1x_1}G_3\partial_{\xi_2}x_1-\partial^2_{\eta_1x_2}G_3\partial_{\xi_2}x_2.
\end{aligned}
$$

Since $G_3$ satisfies \eqref{cond for y2,x2} in coordinates $(\eta,x)$, the monomials of $G_3$ containing $\eta_2^{k_2}x_2^{l_2}$ must satisfy $k_2+l_2\neq 1$.  Assuming that $\eta_2=x_2=0$, it follows from the relations above that
$$
\begin{aligned}
\partial_{\eta_2} y_1=\partial^2_{x_1x_1}G_3\partial_{\eta_2}x_1,\quad &\partial_{\eta_2}x_1(1+\partial^2_{\eta_1x_1}G_3)=0,\\
\partial_{\xi_2} y_1=\partial^2_{x_1x_1}G_3\partial_{\xi_2}x_1,\quad &\partial_{\xi_2}x_1(1+\partial^2_{\eta_1x_1}G_3)=0.
\end{aligned}$$
Since $\partial_{\eta_1x_1}^2G_3$ is a first order homogeneous polynomial,  $1+\partial_{\eta_1x_1}^2G_3\neq 0$ near the origin. This implies
$\partial_{\eta_2}x_1=\partial_{\xi_2}x_1=0=\partial_{\eta_2}y_1=\partial_{\xi_2}y_1=0$.  From \eqref{equ: d of H3}, together with conditions \eqref{cond for y2,x2} for $H$, we conclude that $\hat H_3(\eta_1,\eta_2,\xi_1,\xi_2)$ satisfies \eqref{cond for y2,x2} for every $(\eta_1,\xi_1)$ close to $(0,0)$.

By induction, we assume that  \eqref{cond for y2,x2} holds for $G_i, \hat H_i, i=3,\ldots, s-1$. We aim to prove that \eqref{cond for y2,x2} also holds for $G_s, \hat H_s$, where $G_s$ determines the generating function $W_s$ and the induced canonical transformation $\Phi_s$. We see $\hat H_{s-1}$ as a function of $(y,x)$ and $\hat H_s= \hat H_{s-1} \circ \Phi_s$ as a function of $(\eta,\xi)$. We have
\begin{equation}\label{expending s}
\begin{aligned}
\hat H_{s-1}&=\Gamma^{(2)}+\ldots+\Gamma^{(s-1)}+\mathcal{O}_{s}\\
\hat H_{s}&=\Gamma^{(2)}+\ldots+\Gamma^{(s-1)}+\Gamma^{(s)}+\mathcal{O}_{s+1}.
\end{aligned}
\end{equation}
Using that $y=\eta+\partial_xG_s$, $x=\xi-\partial_\eta G_s$, we choose $G_s$ so that
$-D\cdot G_s+\Gamma^{(s)}= \hat H_{s-1}^{(s)},$ where  $\hat H^{(s)}_{s-1}$ is the $s$-homogeneous part of $\hat H_{s-1}$ and $D\cdot \Gamma^{(s)}=0$. 
Using that the monomials are eigenfunctions of $D$, we see as before that both $G_s$ and $\Gamma^{(s)}$ satisfy \eqref{cond for y2,x2} in the corresponding coordinates. By a similar discussion, we conclude that $\Psi_s$ sends the coordinate plane $\{\eta_2=\xi_2=0\}$ to $\{y_2=x_2=0\}$. We also compute
$$
\begin{aligned}
\partial_{\eta_2}\hat H_s&=\partial_{y_1} \hat H_{s-1} \partial_{\eta_2} y_1+\partial_{y_2} \hat H_{s-1}
\partial_{\eta_2} y_2+\partial_{x_1} \hat H_{s-1}
\partial_{\eta_2} x_1+\partial_{x_2} \hat H_{s-1} \partial_{\eta_2} x_2,\\
\partial_{\xi_2}\hat H_s&=\partial_{y_1} \hat H_{s-1} \partial_{\xi_2} y_1+\partial_{y_2} \hat H_{s-1} \partial_{\xi_2} y_2+\partial_{x_1} \hat H_{s-1} \partial_{\xi_2} x_1+\partial_{x_2} \hat H_{s-1} \partial_{\xi_2} x_2,
\end{aligned}
$$
and
$$
\begin{aligned}
\partial_{\eta_2} y_1&= \partial^2_{x_1\eta_2}G_s+\partial^2_{x_1x_1}G_s\partial_{\eta_2}x_1+\partial^2_{x_1x_2}G_s\partial_{\eta_2}x_2\\
\partial_{\xi_2} y_1&=\partial^2_{x_1x_1}G_s\partial_{\xi_2}x_1+\partial^2_{x_1x_2}G_s\partial_{\xi_2}x_2\\
\partial_{\eta_2} x_1&= -\partial^2_{\eta_1\eta_2}G_s-\partial^2_{\eta_1x_1}G_s\partial_{\eta_2}x_1-\partial^2_{\eta_1x_2}G_s\partial_{\eta_2}x_2\\
\partial_{\xi_2} x_1&=-\partial^2_{\eta_1x_1}G_s\partial_{\xi_2}x_1-\partial^2_{\eta_1x_2}G_s\partial_{\xi_2}x_2.
\end{aligned}
$$
Since $G_s$ satisfies \eqref{cond for y2,x2} for every $(\eta_1,x_1)$ near $(0,0)$, we can follow the same reasoning as before to conclude that \eqref{cond for y2,x2} also holds for $\hat H_s$ for every $(\eta_1,\xi_1)$ near $(0,0)$. This finishes the proof.
\end{proof}

\section{Non-resonant Hopf links for \texorpdfstring{$H_N$}{Lg}}\label{sec normal form}

In this section, we study the dynamics of the Hamiltonian $H_N$ as given in \eqref{An}-\eqref{hatH}. More precisely, we present sufficient conditions for the existence of a non-resonant Hopf link on the sphere-like component near the origin.

For $E>0$ small, denote by $\Sigma^N_E$ the sphere-like component of the energy surface  $H_N^{-1}(E)$. In complex coordinates $z_j=\xi_j + i\eta_j$, $\bar z_j=\xi_j-i\eta_j$, $j=1,2,$ the corresponding Hamiltonian system becomes
\begin{equation}\label{H system 1}
 \dot{\bar z}_1  =2i \partial_{z_1}H_N, \qquad
 \dot{\bar z}_2  =2i\partial_{z_2}H_N,
\end{equation}
where $\partial_{z_j}=\frac{1}{2}(\partial_{\xi_j}-i\partial_{\eta_j})$. We look for solutions $\gamma(t)$ to these equations of the following form
$$
(z_1(t),z_2(t))=(c_1e^{i\beta \alpha_1 t},c_2e^{i\beta \alpha_2t}), \ \ \ c_1,c_2\in \C \ \ \mbox{ and } \ \ \beta \in \R \setminus \{0\}.
$$
If $H_N$ is resonant, then these solutions are periodic. If $H_N$ is non-resonant,
these solutions are periodic only if $c_1$ or $c_2$ is zero. Notice that, for $\bt=1$,  these are precisely the solutions to the Hamiltonian equations
associated with the quadratic part $H_2$ of $H_N$. Moreover, from the definition of $\sigma$ in equation
\eqref{sg} we obtain the important fact along $\gamma(t)$
$$
\sg|_{\gamma} = c_2^{m_2}\bar c_1^{|m_1|} \qquad \Rightarrow \qquad \bar \sigma|_\gamma = c_1^{|m_1|}\bar c_2^{m_2},
$$
i.e., $\sigma|_\gamma$ does not depend on $t$. This fact and the expression for $H_N$ imply that
$H_N(\gamma(t))= H_N (\gamma(0))$, namely
for any $c_1,c_2$ and $\bt$, the curve $\gamma(t)$ is contained in a
level set of $H_N$.

Equations \eqref{An}, \eqref{hatH} and \eqref{H system 1} give
\begin{equation}\label{H system 2}
\begin{aligned}
 -i\beta\alpha_1 \bar z_1 & = \dot{\bar z}_1= 2i\bar z_1\big\{
 \frac{\alpha_1}{2}+\partial_{|z_1|^2}A_0(|c_1|^2,|c_2|^2)
+\sum_{1\leq n\le \lfloor N/(|m_1|+m_2)\rfloor}[\sigma^n\partial_{|z_1|^2}A_n(|c_1|^2,|c_2|^2)\\
&  +
\bar \sigma^n\partial_{|z_1|^2}\bar A_n(|c_1|^2,|c_2|^2)]\big\}
 +2i\sum_{1\leq n\le \lfloor N/(|m_1|+m_2)\rfloor}\partial_{z_1} \bar \sigma^n \bar A_n(|c_1|^2,|c_2|^2) \\
 -i\beta\alpha_2 \bar z_2 & = \dot{\bar z}_2 = 2i\bar z_2\big\{
 \frac{\alpha_2}{2}+\partial_{|z_2|^2}A_0(|c_1|^2,|c_2|^2)
+\sum_{1\leq n\le \lfloor N/(|m_1|+m_2)\rfloor}[\sigma^n\partial_{|z_2|^2}A_n(|c_1|^2,|c_2|^2)\\
&  +
\bar \sigma^n\partial_{|z_2|^2}\bar A_n(|c_1|^2,|c_2|^2)]\big\}
 +2i\sum_{1\leq n\le \lfloor N/(|m_1|+m_2)\rfloor}\partial_{z_2}\sigma^n A_n(|c_1|^2,|c_2|^2),
\end{aligned}
\end{equation}
where $\sigma=z_2^{m_2}\bar z_1^{|m_1|}$ and the expressions above are evaluated at $\gamma(t)$. Equations \eqref{H system 2} imply that $\gamma(t)$ is a solution to the Hamiltonian
differential equations  if
and only if the constants $c_1$, $c_2$ and $\bt$ satisfy a set of algebraic
equations that do not involve the parameter $t$. Moreover, equations \eqref{H system 2} can be written as
\begin{equation}\label{variational}
\partial_{z_j} H_2+\beta^{-1}\partial_{z_j}H_N=0, \qquad j=1,2,
\end{equation}
and we are led to the following variational problem for the determination of $c_1, c_2$ and $\beta$:
\begin{itemize}
\item[(P)] {\it Find critical points of $H_2 = \frac{\alpha_1}{2} |z_1|^2 + \frac{\alpha_2}{2}|z_2|^2$ restricted to $\Sigma^N_E\subset H_N^{-1}(E)$, where $E>0$ is small.}
\end{itemize}
Notice that $\beta^{-1}$ corresponds to the Lagrange multiplier of the associated critical point. Since $\Sigma_E$ is a regular sphere-like hypersurface, problem (P) has at least two distinct solutions that realize extreme values of $H_2$ restricted to $\Sigma^N_E$. Such a solution determines the constants $c_1, c_2$ and $\beta$.
We have proved the following theorem.

\begin{prop}
\label{existence}
For every $E>0$ sufficiently small,  the Hamiltonian differential equations
associated with $H_N$ admit two geometrically distinct solutions $\gamma_1,\gamma_2 \subset \Sigma^N_E$ both of the form
$$
(z_1(t),z_2(t))=(c_1e^{i\beta \alpha_1 t},c_2e^{i\beta \alpha_2t}), \quad \forall t.
$$
The constants $c_1,c_2\in \C$ and $\beta \in \R \setminus \{0\}$ can be found by solving the variational
problem (P) above.
If $H_N$ is resonant, then these solutions are periodic. If $H_N$ is non-resonant, these solutions are periodic only if $c_1$ or $c_2$ is zero.
\end{prop}

There are several remarks concerning problem (P).
At first, the function $H_N$ in equation \eqref{hatH} can be decomposed as
$H_N=H_2+R_N$. This means that $H_2$ restricted to $\Sigma^N_E$ coincides with the function $E-R_N$ restricted to $\Sigma^N_E$. Therefore, we can change the variational
problem (P) by that of finding the critical points of $R_N$ restricted to the
set $\Sigma^N_E$. Second, all  functions
 $H_2$, $R_N$, and $H_N$ are invariant under the action of $\R$ given by
$t\cdot(z_1,z_2)\to (z_1e^{i\alpha_1t}, z_2e^{i\alpha_2t})$. This action can be used
to reduce the dimension of the variational problem. For instance, in the particular
case where  $\alpha_1=\alpha_2$, this action reduces to an $S^1$-action, and the problem reduces to finding critical points of $R_N$ restricted to
 $\Sigma^N_E/ S^1 \cong S^2$.

Let us assume that  $c_1=0$ and $c_2\neq 0.$ Then  $\partial_{z_1}H_N|_\gamma=0$ and the first equation of \eqref{H system 2} is automatically satisfied, except in the case $|m_1|=1\Rightarrow m_2=1$, where it is equivalent to
\begin{equation}\label{cond: gamma2 m2=1}
0=\bar A_1(0,|c_2|^2)=\sum_{3\leq 2k_2+2\leq N}\bar a_{0,k_2+1,1,k_2}|c_2|^{2k_2}.
\end{equation}
This condition is not verified for $|c_2|>0$ arbitrarily small if some coefficient $\bar a_{0,k_2+1,1,k_2}$ in the sum above does not vanish.
The second equation in \eqref{H system 2} is equivalent to
\begin{equation}\label{c1=0}
\begin{aligned}
0 & =\frac{\alpha_2}{2}(1+\beta) + \partial_{|z_2|^2}A_0(0,|c_2|^2)\\
& = \frac{\alpha_2}{2}(1+\beta) + \sum_{3\leq 2k_2\leq N} k_2a_{0,k_2,0,k_2}|c_2|^{2k_2-2},
\end{aligned}
\end{equation}
which is solvable for $\beta$, for every $|c_2|>0$. Such solutions near the origin are in one-to-one correspondence with the energy surfaces $H_N^{-1}(E)$, where $E>0$ is small.

Now assume that $c_1\neq 0$ and $c_2=0$. Then $\partial_{z_2}H_N=0$ and the second equation of \eqref{H system 2} is automatically satisfied, except in the case $m_2=1\Rightarrow 1\leq |m_1|<+\infty$, where it is equivalent to
\begin{equation}\label{cond: gamma1 m2=1}
0=A_1(|c_1|^2,0)=\sum_{3\leq 2k_1+1+|m_1|\leq N} a_{k_1,1,k_1+|m_1|,0}|c_1|^{2k_1}.
\end{equation}
This condition is not verified for $|c_1|>0$ arbitrarily small if some coefficient $a_{k_1,1,k_1+|m_1|,0}$ above does not vanish.
The first equation in \eqref{H system 2} is equivalent to
\begin{equation}\label{c2=0}
\begin{aligned}
0 & =\frac{\alpha_1}{2}(1+\beta) + \partial_{|z_1|^2}\bar A_0(|c_1|^2,0)\\
& = \frac{\alpha_1}{2}(1+\beta) + \sum_{3\leq 2k_1\leq N} k_1a_{k_1,0,k_1,0}|c_1|^{2k_1-2},
\end{aligned}
\end{equation}
which is solvable for $\beta$, for every $|c_1|>0$. Such solutions near the origin are in one-to-one correspondence with the energy surfaces $H_N^{-1}(E)$, where $E>0$ is small.

We have proved the following proposition.

\begin{prop}
\label{prop_explicit}The following assertions hold: \begin{itemize}\item[(i)] If $\alpha_2 \neq \alpha_1$ (i.e., $|m_1|>1$),
 then, for every $E>0$ sufficiently small, the Hamiltonian differential equations associated with $H_N$ admit a periodic solution  $\gamma_2 \subset \Sigma^N_E,$ of the form
$$
(z_1(t),z_2(t))=(0,c_2e^{i\beta \alpha_2t})
$$
for every $t$. The constants $c_2>0$ and $\beta \neq 0$ satisfy \eqref{c1=0}.

\item[(ii)] If $\alpha_2$ is not a multiple of $\alpha_1$ (i.e., $m_2>1$), or if $\alpha_2$ is a multiple of $\alpha_1$ (i.e., $m_2=1$) and condition \eqref{cond: gamma1 m2=1} holds, then for every $E>0$ sufficiently small,  the Hamiltonian differential equations associated with $H_N$ admit a periodic solution  $\gamma_1\subset \Sigma^N_E,$ of the form
$$
(z_1(t),z_2(t))=(c_1e^{i\beta \alpha_1 t},0)
$$
for every $t$. The constants $c_1>0$ and $\beta\neq 0$ satisfy \eqref{c2=0}.

\item[(iii)] If $\alpha_2=\alpha_1$ (i.e., $|m_1|=m_2=1$) and conditions \eqref{cond: gamma2 m2=1}, \eqref{cond: gamma1 m2=1} hold, then for every $E>0$ sufficiently small,  the Hamiltonian differential equations associated with $H_N$ admit both periodic solutions $\gamma_2,\gamma_1\subset \Sigma^N_E$ as in (i) and (ii), respectively.
\end{itemize}
In particular, if $H_N$ satisfies
$a_{0,k_2+1,1,k_2}=a_{k_1,1,k_1+|m_1|,0}=0$, for every $3\leq 2k_2+2\leq N$ and $3\leq 2k_1+1+|m_1|\leq N$, then there always exists a pair of periodic orbits $\gamma_1,\gamma_2$ as above.
\end{prop}

If $\alpha_2$ is not a multiple of $\alpha_1$, then both conditions in Proposition \ref{prop_explicit}-(i) and (ii) hold for $H_N$. In particular, the link $\gamma_1\cup \gamma_2\subset \Sigma^N_E$ given in Proposition \ref{prop_explicit} is a Hopf link for every $E>0$ sufficiently small. If $\alpha_2$ is a multiple of $\alpha_1$, then the periodic solutions must be computed using the coefficients of the polynomial $H_N$.

\subsection{Rotation numbers} Our next task is to compute the rotation numbers of $\gamma_1$ and $\gamma_2$.

The rotation numbers of the periodic orbits $\gamma_1$ and $\gamma_2$ given in Proposition \ref{prop_explicit} can be computed using the formulas in Appendix \ref{sec_quaternion}.
To explain the detail, it is necessary to linearize the differential equations \eqref{H system 1} along a periodic orbit $\gamma$.
Using complex coordinates $w_j$ and $\bar w_j$ tangent to $z_j$ and $\bar z_j$, respectively, the linearized equations become
\begin{equation}
\begin{aligned}
 \dot {\bar w}_1&=  2i\left( \partial_{z_1z_1} H_N w_1+ \partial_{\bar z_1z_1}H_N\bar w_1 +\partial_{z_2z_1} H_N w_2+
\partial_{\bar z_2z_1} H_N\bar w_2\right),\\
 \dot {\bar w}_2&= 2i\left( \partial_{z_1z_2} H_N w_1+ \partial_{\bar z_1z_2} H_N \bar w_1 +\partial_{z_2z_2} H_N w_2+
\partial_{\bar z_2z_2} H_N\bar w_2\right),
\end{aligned}
\end{equation}
where the second derivatives of $H_N$ are evaluated at $\gamma(t)$.

If $\gamma=\gamma_2$, then for $|m_1|\geq3$, we have $z_1(t)=0$, $z_2(t)=c_2e^{i\beta \alpha_2 t}, c_2>0$,  and the expression for $H_N$ in equation \eqref{hatH} implies that
the $w_1$-component of the linearized equations (the only component used in the computation of the rotation number) satisfies
\begin{equation}\label{w1bdot}
 \dot {\bar w}_1= 2i\big\{\frac{\alpha_1}{2}  +\partial_{|z_1|^2} A_0(0,c_2^2) \big\}\bar w_1.
\end{equation}
If $|m_1|=2$ or if $|m_1|=1$ and condition \eqref{cond: gamma2 m2=1} holds, then $m_2=1$ and the linearized equation of $\gamma_2$ can be written as
\begin{equation}\label{linw2a}
\dot{\bar w}_1=2i\big\{\frac{\alpha_1}{2}  +\partial_{|z_1|^2} A_0(0,c_2^2) \big\}\bar w_1+\underbrace{4i\bar z_2^{2/|m_1|}\bar A_{2/|m_1|}(0,c_2^2)}_{=2i\partial_{z_1z_1}H_N|_{\gamma_2}}w_1.
\end{equation}

Similarly, if $\gamma=\gamma_1$, then
 $z_1(t)=c_1e^{i\beta\alpha_1 t}, c_1>0$, $z_2(t)=0$,  and the expression for $H_N$ in equation \eqref{hatH} implies that
for $m_2\ge 3$ the $w_2$-component of the linearized equations satisfies
\begin{equation}\label{linw2}
 \dot {\bar w}_2= 2i\big\{\frac{\alpha_2}{2}  +\partial_{|z_2|^2} A_0(c_1^2,0) \big\}\bar w_2.
\end{equation}
If $m_2=2$, or if $m_2=1$ and  condition \eqref{cond: gamma1 m2=1} holds, then the linearized equation becomes
\begin{equation}\label{linw2b}
 \dot {\bar w}_2= 2i\big\{\frac{\alpha_2}{2}  +\partial_{|z_2|^2} A_0(c_1^2,0)\big\} \bar w_2
+\underbrace{4i \bar z_1^{2|m_1|/m_2}(t)A_{2/m_2}(c_1^2,0)}_{=2i\partial_{z_2z_2}H_N|_{\gamma_1}} w_2.
\end{equation}

\subsection{The case \texorpdfstring{$m_2\geq 3$}{Lg}} Here, we assume that $|m_1|>m_2\ge 3$. Recall that $|m_1|=m_2=+\infty$ if $H_N$ is non-resonant at the origin. From the previous discussion, the rotation numbers of $\gamma_2$ and $\gamma_1$ are computed from the linearized equations \eqref{w1bdot} and \eqref{linw2}, respectively, regardless of whether $H_N$ is resonant or non-resonant at the origin.

We compute the rotation number $\rho_1$ of $\gamma_1\subset \Sigma^N_E$ as a function of $E$, where $c_1=c_1(E)>0$ and $\beta=\beta(E) \neq 0$ satisfy \eqref{c2=0}.
Since $\gamma_1 \subset H_N^{-1}(E)$, we obtain the implicit function
$$
\begin{aligned}
E & = H_N(\gamma_1(t))=\frac{\alpha_1}{2}c_1^2 + A_0(c_1^2,0)\\ & = \frac{\alpha_1}{2}c_1^2+ \sum_{3\leq 2k_1\leq N}a_{k_1,0,k_1,0}c_1^{2k_1},
\end{aligned}
$$
for every $E>0$ sufficiently small. Since $\partial_{|z_1|^2}A_0(0,0)=0$, there
exists a smooth function $s_1=s_1(E)\geq 0$, with $s_1(0)=0$ and $s_1'(0)=2/\alpha_1>0$, defined for $E\geq 0$ sufficiently small, such that
$c_1(E)^2=s_1(E)=\frac{2E}{\alpha_1} + O(E^2)$ for every $E>0$ sufficiently small. We thus have
$$
c_1(E)^2 = \frac{2E}{\alpha_1} -\frac{2}{\alpha_1}\sum_{3\leq 2k_1\leq N}a_{k_1,0,k_1,0}\left(\left(\frac{2E}{\alpha_1} \right)^{k_1}+O(E^{k_1+1})\right).
$$
Since $\dot {\bar z}_1 = 2i \partial_{z_1} H_N(\gamma_1) = 2i\bar z_1\{\frac{\alpha_1}{2} + \partial_{|z_1|^2}A_0(|c_1|^2,0)\}$, the frequency $\beta \alpha_1$ of $\gamma_1$ is given by
\begin{equation}\label{Omega1}
\begin{aligned}
\omega_1(E) :=-\beta \alpha_1&=\alpha_1 + 2\partial_{|z_1|^2}A_0(c_1(E)^2,0)=\alpha_1 + 2\sum_{3\leq 2k_1 \leq N}k_1 a_{k_1,0,k_1,0} c_1(E)^{2k_1-2} \\
& =\alpha_1 + 2\sum_{3\leq 2k_1 \leq N}k_1 a_{k_1,0,k_1,0}\left(\left(\frac{2E}{\alpha_1}\right)^{k_1-1}+O(E^{k_1})\right).
\end{aligned}
\end{equation}
The period of $\gamma_1$ is $T_1(E)=2\pi/\omega_1(E)$.

Along the orbit $\gamma_1$, for which $(\eta_1(t),\eta_2(t),\xi_1(t),\xi_2(t))=(c_1\sin(-\omega_1t),0,c_1\cos(-\omega_1t),0)$, $c_1>0$, we have
$$\nabla H_N(\gamma_1)=-J\dot \gamma_1=c_1\omega_1(-\sin(\omega_1 t)\partial_{\eta_1}+\cos(\omega_1 t)\partial_{\xi_1})=\omega_1(\eta_1\partial_{\eta_1}+\xi_1\partial_{\xi_1}).$$
As explained in the Appendix \ref{sec_quaternion}, the quaternion trivialization along $\gamma_1$ is induced by the vector $V_0$ normal to $\Sigma^N_E$
\begin{equation}
\begin{aligned}
V_0(\gamma_1) & = \frac{\nabla  H_N(\gamma_1)}{|\nabla H_N(\gamma_1)|}=c_1^{-1}(\eta_1\partial_{\eta_1}+\xi_1\partial_{\xi_1})=e^{i\omega_1t}\partial_{\bar z_1}+e^{-i\omega_1t}\partial_{z_1}.
\end{aligned}
\label{V0}
\end{equation}
Hence the orthonormal basis tangent to $H_N^{-1}(E)$ along $\gamma_1$ is given by
\begin{equation}\label{equ: V1,V2,V3}
\begin{aligned}
    V_1(\gamma_1) & =c_1^{-1}(-\xi_1\partial_{\eta_2}-\eta_1\partial_{\xi_2})= i(e^{-i\omega_1t}\partial_{\bar z_2}-e^{i\omega_1t}\partial_{z_2}),\\
    V_2(\gamma_1) & =c_1^{-1}(\eta_1\partial_{\eta_2}-\xi_1\partial_{\xi_2}) = -e^{-i\omega_1 t}\partial_{\bar z_2} - e^{i\omega_1 t}\partial_{z_2}, \\
    V_3(\gamma_1) & =c_1^{-1}(-\xi_1\partial_{\eta_1}+\eta_1\partial_{\xi_1})= i(e^{i\omega_1t}\partial_{\bar z_1} - e^{-i\omega_1t}\partial_{z_1}),
\end{aligned}
\end{equation}
where $\omega_1=\omega_1(E)$ satisfies \eqref{Omega1}. While the basis $\{V_1(\gamma_1),V_2(\gamma_1)\}$ rotates with  frequency $-\omega_1$ with respect to the frame $\partial_{\bar z_2},$ we see from \eqref{linw2} that the linearized flow with respect to the frame $\partial_{\bar z_2}$ has  frequency
\begin{equation}\label{hatOmega2}
\begin{aligned}
\hat \omega_2(E):&= \alpha_2 + 2\partial_{|z_2|^2}A_0(c_1(E)^2,0) = \alpha_2 +2\sum_{3\leq 2k_1\leq N}a_{k_1-1,1,k_1-1,1} c_1(E)^{2k_1-2}\\
& = \alpha_2 +2\sum_{3\leq 2k_1\leq N}a_{k_1-1,1,k_1-1,1}\left(\left(\frac{2E}{\alpha_1}\right)^{k_1-1}+O(E^{k_1})\right)
\end{aligned}
\end{equation}
Hence, the linearized flow on the frame $\{V_1(\gamma_1),V_2(\gamma_2)\}$ rotates with frequency $\hat \omega_2(E) + \omega_1(E)$, and the rotation number of $\gamma_1\subset H_N^{-1}(E)$ is then given by
\begin{equation}\label{rho1 m2>2}
\rho_1(E)  = T_1(E)\cdot \frac{\omega_1(E) + \hat \omega_2(E)}{2\pi} = 1 + \frac{\hat \omega_2(E)}{\omega_1(E)}.
\end{equation}

The computation of the rotation number of $\gamma_2\subset \Sigma^N_E$ is similar. As before, we write $|c_2|^2$ as a function of $E>0$, namely
$$
c_2(E)^2 = \frac{2E}{\alpha_2} - \frac{2}{\alpha_2}\sum_{3\leq 2k_2\leq N}a_{0,k_2,0,k_2}\left(\left(\frac{2E}{\alpha_2} \right)^{k_2}+O(E^{k_2+1})\right).
$$
We obtain
\begin{equation}\label{rho2 m2>2}
\rho_2(E) = 1 + \frac{\hat \omega_1(E)}{\omega_2(E)},
\end{equation}
where
\begin{equation}\label{Omega2}
\begin{aligned}
\omega_2(E)  :=-\beta \alpha_2&=\alpha_2 + 2\partial_{|z_2|^2}A_0(0,c_2(E)^2) = \alpha_2 + 2\sum_{3\leq 2k_2\leq N} k_2 a_{0,k_2,0,k_2}  c_2(E)^{2k_2-2}\\
& = \alpha_2 + 2\sum_{3\leq 2k_2\leq N} k_2 a_{0,k_2,0,k_2} \left(\left(\frac{2E}{\alpha_2}\right)^{k_2-1}+O(E^{k_2})\right)
\end{aligned}
\end{equation}
is the frequency of $\gamma_2$ and
\begin{equation}\label{hatOmega1}
\begin{aligned}
\hat \omega_1(E) :&= \alpha_1 + 2\partial_{|z_1|^2}A_0(0,c_2(E)^2) = \alpha_1 + 2\sum_{3\leq 2k_2\leq N} a_{1,k_2-1,1,k_2-1} c_2(E)^{2k_2-2}\\
& = \alpha_1 + 2\sum_{3\leq 2k_2\leq N} a_{1,k_2-1,1,k_2-1}\left(\left(\frac{2E}{\alpha_2}\right)^{k_2-1}+O(E^{k_2})\right)
\end{aligned}
\end{equation}
is the frequency of the linearized flow on the frame $\partial_{\bar z_1}$ along $\gamma_2$.

The twist condition \eqref{twist} becomes
$$
(\rho_1(E)-1)(\rho_2(E)-1)=\frac{\hat \omega_2(E)}{\omega_1(E)}\frac{\hat \omega_1(E)}{\omega_2(E)}\ne 1,
$$
which is equivalent to
$ \hat \omega_2(E) \hat \omega_1(E) \neq \omega_1(E) \omega_2(E).$

Let $2\leq \nu \leq \lfloor N/2 \rfloor$ be the smallest integer such that one of the following numbers does not vanish
\begin{equation}\label{condition of nv}
a_{0,\nu,0,\nu}, \quad a_{\nu,0,\nu,0}, \quad a_{1,\nu-1,1,\nu-1}, \quad a_{\nu-1,1,\nu-1,1}.
\end{equation}
If such a $\nu$ does not exist, then $
\hat \omega_2(E) \hat \omega_1(E) = \omega_1(E) \omega_2(E)=\alpha_1\alpha_2$, and the link $\gamma_1\cup \gamma_2$ is resonant.
We thus assume that such $\nu$ exists. From this,  we get the following power series expansions:
\begin{align*}
\hat \omega_2(E) & =\alpha_2+2a_{\nu-1,1,\nu-1,1}\left(\frac{2E}{\alpha_1}\right)^{\nu-1}
+O(E^\nu),\
&\omega_1(E) & = \alpha_1+2\nu a_{\nu,0,\nu,0}\left(\frac{2E}{\alpha_1}\right)^{\nu-1}+O(E^\nu)\\
\hat \omega_1(E) & = \alpha_1 +2a_{1,\nu-1,1,\nu-1}\left(\frac{2E}{\alpha_2}\right)^{\nu-1}
+O(E^\nu),\
&\omega_2(E) & = \alpha_2+2\nu a_{0,\nu,0,\nu}\left(\frac{2E}{\alpha_2}\right)^{\nu-1}+O(E^\nu).
\end{align*}
From these expansions, we obtain
\begin{equation}\label{omega/Omega}
\begin{aligned}
\frac{\hat \omega_2}{\omega_1}
&=\frac{\alpha_2}{\alpha_1}+\underbrace{(a_{\nu-1,1,\nu-1,1}-\nu a_{\nu,0,\nu,0}\frac{\alpha_2}{\alpha_1})}_{\Omega_{\nu,1}}\left(\frac{2}{\alpha_1}\right)^\nu E^{\nu-1}+O(E^\nu),\\
\frac{\hat \omega_1}{\omega_2}
&=\frac{\alpha_1}{\alpha_2}+\underbrace{(a_{1,\nu-1,1,\nu-1}-\nu a_{0,\nu,0,\nu}\frac{\alpha_1}{\alpha_2})}_{\Omega_{\nu,2}}\left(\frac{2}{\alpha_2}\right)^\nu E^{\nu-1}+O(E^\nu),
\end{aligned}
\end{equation}
and
\begin{equation}
\begin{aligned}
(\rho_1(E)-1)(\rho_2(E)-1) & = \frac{\hat \omega_2(E) \hat \omega_1(E)}{\omega_1(E) \omega_2(E)} = 1+2^\nu\Omega_\nu E^{\nu-1}+O(E^\nu),
\label{equ: nonresonancy for m2>2}
\end{aligned}
\end{equation}
where
\begin{equation}\label{equ: omega nu}
\begin{aligned}
\Omega_\nu:&=\frac{\Omega_{\nu,2}}{\alpha_1\alpha_2^{\nu-1}}+\frac{\Omega_{\nu,1}}{\alpha_2\alpha_1^{\nu-1}}\\
&=\frac{1}{\alpha_1\alpha_2}
\left(\frac{a_{\nu-1,1,\nu-1,1}}{\alpha_1^{\nu-2}}+\frac{a_{1,\nu-1,1,\nu-1}}{\alpha_2^{\nu-2}}\right)-\nu\left(\frac{a_{\nu,0,\nu,0}}{\alpha_1^{\nu}}+\frac{a_{0,\nu,0,\nu}}{\alpha_2^{\nu}}\right).
\end{aligned}
\end{equation}
Hence, the twist condition \eqref{twist} is verified by the Hopf link $\gamma_1 \cup \gamma_2$ if the coefficient $\Omega_\nu$ is different from zero. This is summarized in the following theorem.

\begin{thm}
Assume that $|m_1|>m_2 \geq 3$. Let $\gamma_1,\gamma_2 \subset H_N^{-1}(E)$ be the periodic orbits forming a Hopf link, where $E>0$ is sufficiently small. Assume that there exists a least $2\leq \nu \leq \lfloor N/2\rfloor$ so that at least one of the coefficients in \eqref{condition of nv} does not vanish. Let $\Omega_\nu$ be given in \eqref{equ: omega nu}. If $\Omega_\nu \neq 0$, then the Hopf link $L=\gamma_1 \cup \gamma_2$ is non-resonant for every $E>0$ sufficiently small. Moreover, if $\nu$  does not exist, then the Hopf link $L$ is resonant.
\end{thm}

\subsection{The general linearized equation}\label{sec linearization}
To further generalize the twist condition \eqref{equ: nonresonancy for m2>2} to the particular cases $m_2=1,2$, we need to study the linearized equations \eqref{linw2a} and \eqref{linw2b}. 

Assume that the periodic solution $\gamma_1$ of the form $(z_1(t),z_2(t))=(c_1e^{-i\omega_1 t},0)$ exists for the flow of $H_N$, where $c_1>0$ and $\beta\neq 0$. The existence of $\gamma_1$ follows from condition \eqref{cond: gamma1 m2=1} whenever $m_2=1$. Using \eqref{equ: V1,V2,V3}, we compute
\begin{equation}\label{H33}
\begin{aligned}
H_{N33}(\gamma_1)&=(V_3,\nabla^2H_N(\gamma_1)V_3)\\
&=c_1^{-2}\big(\xi_1^2\partial_{\eta_1\eta_1}H_N|_{\gamma_1}-2\eta_1\xi_1\partial_{\eta_1\xi_1}H_N|_{\gamma_1}+\eta_1^2\partial_{\xi_1\xi_1}H_N|_{\gamma_1}\big)\\
&=c_1^{-2}\big((4\xi_1^2\eta_1^2-8\eta_1^2\xi_1^2+4\eta_1^2\xi_1^2)\partial_{c_1^2}^2H_N|_{\gamma_1}+2(\xi_1^2+\eta_1^2)\partial_{c_1^2}H_N|_{\gamma_1}\big)\\&=c_1^{-2}(\eta_1^2+\xi_1^2)\cdot 2\partial_{c_1^2}H_N|_{\gamma_1}=\omega_1(E),
\end{aligned}
\end{equation}
where
$$\begin{aligned}
\partial_{\eta_1\eta_1}H_N\big|_{\gamma_1}&=
4\eta_1^2\partial_{c_1^2}^2H_N(\gamma_1)+2\partial_{c_1^2}H_N(\gamma_1),\\
\partial_{\xi_1\xi_1}H_N\big|_{\gamma_1}&=
4\xi_1^2\partial_{c_1^2}^2H_N(\gamma_1)+2\partial_{c_1^2}H_N(\gamma_1),\\
\partial_{\xi_1\eta_1}H_N\big|_{\gamma_1}&=
4\eta_1\xi_1\partial_{c_1^2}^2H_N(\gamma_1).
\end{aligned}$$
Using that $\partial_{\xi_j}=\partial_{\bar z_j}+\partial_{z_j}$ and $\partial_{\eta_j}=-i(\partial_{\bar z_j}-\partial_{z_j}), j=1,2$, we  obtain
\begin{equation}\label{H11}
\begin{aligned}
&\quad\ H_{N11}(\gamma_1)=(V_1,\nabla^2H_N(\gamma_1)V_1)\\
&=c_1^{-2}\big(\xi_1^2\partial_{\eta_2\eta_2}H_N|_{\gamma_1}+2\eta_1\xi_1\partial_{\eta_2\xi_2}H_N|_{\gamma_1}+\eta_1^2\partial_{\xi_2\xi_2}H_N|_{\gamma_1}\big)\\
&=c_1^{-2}\big(-\xi_1^2(\partial_{\bar z_2}-\partial_{z_2})^2H_N|_{\gamma_1}-2i\eta_1\xi_1(\partial_{\bar z_2}^2-\partial_{ z_2}^2)H_N|_{\gamma_1}+\eta_1^2(\partial_{\bar z_2}+\partial_{z_2})^2H_N|_{\gamma_1}\big)\\
&=c_1^{-2}\big(2c_1^2\partial_{\bar z_2 z_2}H_N|_{\gamma_1}-z_1^2\partial_{\bar z_2\bar z_2}H_N|_{\gamma_1}-\bar z_1^2\partial_{z_2z_2}H_N|_{\gamma_1}\big),
\end{aligned}
\end{equation}
Similarly, we obtain that
\begin{equation}\label{H22}
\begin{aligned}
&\quad\ H_{N22}(\gamma_1)=(V_2,\nabla^2H_N(\gamma_1)V_2)\\
&=c_1^{-2}\big(\eta_1^2\partial_{\eta_2\eta_2}H_N|_{\gamma_1}-2\eta_1\xi_1\partial_{\eta_2\xi_2}H_N|_{\gamma_1}+\xi_1^2\partial_{\xi_2\xi_2}H_N|_{\gamma_1}\big)\\
&=c_1^{-2}\big(-\eta_1^2(\partial_{\bar z_2}-\partial_{z_2})^2H_N|_{\gamma_1}+2i\eta_1\xi_1(\partial_{\bar z_2}^2-\partial_{ z_2}^2)H_N|_{\gamma_1}+\xi_1^2(\partial_{\bar z_2}+\partial_{z_2})^2H_N|_{\gamma_1}\big)\\
&=c_1^{-2}\big(2c_1^2\partial_{\bar z_2 z_2}H_N|_{\gamma_1}+z_1^2\partial_{\bar z_2\bar z_2}H_N|_{\gamma_1}+\bar z_1^2\partial_{z_2z_2}H_N|_{\gamma_1}\big),
\end{aligned}
\end{equation}
and
\begin{equation}\label{H12}
\begin{aligned}
&\ \quad H_{N12}(\gamma_1)=H_{N21}(\gamma_1)=(V_1,\nabla^2H_N(\gamma_1)V_2)\\
&=c_1^{-2}\big(-\eta_1\xi_1\partial_{\eta_2\eta_2}H_N|_{\gamma_1}+(\xi_1^2-\eta_1^2)\partial_{\eta_2\xi_2}H_N|_{\gamma_1}+\eta_1\xi_1\partial_{\xi_2\xi_2}H_N|_{\gamma_1}\big)\\
&=c_1^{-2}\big(\eta_1\xi_1(\partial_{\bar z_2}-\partial_{z_2})^2H_N|_{\gamma_1}-i(\xi_1^2-\eta_1^2)(\partial_{\bar z_2}^2-\partial_{ z_2}^2)H_N|_{\gamma_1}+\eta_1\xi_1(\partial_{\bar z_2}+\partial_{z_2})^2H_N|_{\gamma_1}\big)\\
&=ic_1^{-2}\big(-z_1^2\partial_{\bar z_2\bar z_2}H_N|_{\gamma_1}+\bar z_1^2\partial_{z_2z_2}H_N|_{\gamma_1}\big).
\end{aligned}
\end{equation}
Recall $\hat \omega_2=2\partial_{|z_2|^2}H_N|_{\gamma_1}=2\partial_{\bar z_2 z_2}H_N|_{\gamma_1}$. We denote
\begin{equation}\label{equ: def of A}
A:=e^{-2i\omega_1 t}\partial_{\bar z_2\bar z_2}H_N|_{\gamma_1},\quad \partial_{\bar z_2\bar z_2}H_N|_{\gamma_1}=\tilde c_1e^{i\tilde \theta_1},\quad \tilde c_1\geq 0.
\end{equation}
We see that $A=\partial_{\bar z_2\bar z_2}H_N|_{\gamma_1}=0$ if $m_2>2$. If $m_2=2$ or if $m_2=1$ and condition \eqref{cond: gamma1 m2=1} holds, then
\begin{equation}\label{partial z_2
z_2 partial bar z_2 z_2}
\begin{aligned}
\partial_{\bar z_2\bar z_2}H_N|_{\gamma_1}&=2z_1^{\frac{2|m_1|}{m_2}}\bar A_{2/m_2}(c_1^2,0)=\sum_{3\leq 2k_1+2+2|m_1|/m_2\leq N}2\bar a_{k_1,2,k_1+\frac{2|m_1|}{m_2},0}c_1^{2k_1}z_1^{\frac{2|m_1|}{m_2}}.
\end{aligned}
\end{equation}
In this case,
$\tilde c_1=2c_1^{\frac{2|m_1|}{m_2}}|\bar A_{2/m_2}(c_1^2,0)|$ and $\tilde \theta_1=\arg \bar A_{2/m_2}(c_1^2,0)-\frac{2|m_1|}{m_2}\omega_1t.$ In particular, $\tilde c_1=0$ if $m_2>2$.
Applying \eqref{thdot}, we obtain that
\begin{equation}\label{equation of dot theta}
\begin{aligned}
\dot\theta(t)&=\omega_1+\cos^2\theta\cdot H_{N11}(\gamma_1)+2\cos\theta\sin\theta\cdot H_{N12}(\gamma_1)+\sin^2\theta\cdot H_{N22}(\gamma_1)\\
&=\omega_1+\hat\omega_2+(\sin^2\theta-\cos^2\theta)(A+\bar A)+2i\cos\theta\sin\theta(\bar A-A)\\
&=\omega_1+\hat\omega_2-2\mathrm{Re}(e^{2\theta i}A),
\end{aligned}
\end{equation}
which coincides with $\dot \theta_1=\omega_1+\hat \omega_2$ if $m_2>2$. Defining $\bar \theta(t):=\arg(e^{2\theta i}A)=2\theta(t)-(2+\frac{2|m_1|}{m_2})\omega_1t+\arg \bar A_1(c_1^2,0)$, we conclude that
\begin{equation}\label{bartheta}
\dot{\bar\theta}(t)=2\dot \theta(t)-\left(2+\frac{2|m_1|}{m_2}\right)\omega_1=2\hat \omega_2-\frac{2|m_1|}{m_2}\omega_1-4\tilde c_1\cos\bar\theta(t).
\end{equation}

Similarly, let us assume that $\gamma_2=(0,c_2e^{-i\omega_2 t})$ exists, and denote $\partial_{\bar z_1\bar z_1}H_N|_{\gamma_2}=\tilde c_2e^{i\tilde \theta_2}$ as before, where $c_2,\tilde c_2>0$, $\beta\neq 0$. Since $\nabla H_N(\gamma_2)=c_2\omega_2(-\sin(\omega_2t)\partial_{\eta_2}+\cos(\omega_2t)\partial_{\xi_2})$, the quaternion trivalization along $\gamma_2$ is given by
\begin{align*}
V_0(\gamma_2)& =c_2^{-1}(\eta_2\partial_{\eta_2}+\xi_2\partial_{\xi_2})=e^{i\omega_2t}\partial_{\bar z_2}+e^{-i\omega_2t}\partial_{z_2},\\
V_1(\gamma_2)& =c_2^{-1}(\xi_2\partial_{\eta_1}+\eta_2\partial_{\xi_1})=i(-e^{-i\omega_2t}\partial_{\bar z_1}+e^{i\omega_2t}\partial_{z_1}),\\
V_2(\gamma_2)& =c_2^{-1}(-\eta_2\partial_{\eta_1}+\xi_2\partial_{\xi_1})=e^{-i\omega_2t}\partial_{\bar z_1}+e^{i\omega_2t}\partial_{z_1},\\
V_3(\gamma_2)& =c_2^{-1}(-\xi_2\partial_{\eta_2}+\eta_2\partial_{\xi_2})=i(e^{i\omega_2t}\partial_{\bar z_2}-e^{-i\omega_2t}\partial_{z_2}).
\end{align*}
Hence, as in \eqref{H33}-\eqref{H12}, we obtain
\begin{equation}\label{Hij gamma2}
\begin{aligned}
H_{N11}(\gamma_2)&=c_2^{-2}(2c_2^2\partial_{\bar z_1 z_1}H_N|_{\gamma_2}-z_2^2\partial_{\bar z_1\bar z_1}H_N|_{\gamma_2}-\bar z_2^2\partial_{z_1z_1}H_N|_{\gamma_2}),\\
H_{N22}(\gamma_2)&=c_2^{-2}(2c_2^2\partial_{\bar z_1 z_1}H_N|_{\gamma_2}+z_2^2\partial_{\bar z_1\bar z_1}H_N|_{\gamma_2}+\bar z_2^2\partial_{z_1z_1}H_N|_{\gamma_2},\\
H_{N12}(\gamma_2)&=ic_2^{-2}(-z_2^2\partial_{\bar z_1}^2H_N|_{\gamma_2}+\bar z_2^2\partial_{z_1}^2H_N|_{\gamma_2}),\ H_{N33}(\gamma_2)=\omega_2(E).
\end{aligned}
\end{equation}
Moreover, we see that $\partial_{z_1z_1}H_N|_{\gamma_2}=\partial_{\bar z_1\bar z_1}H_N|_{\gamma_2}=0$ and $\dot \theta(t)=\omega_2+\hat \omega_1$ if $|m_1|>2$. If $|m_1|\in\{1,2\}$, $m_2=1$ and condition \eqref{cond: gamma2 m2=1} holds, then we have
\begin{equation}\label{partial z_2
z_2 partial bar z_2 z_2b}
\begin{aligned}
\partial_{\bar z_1\bar z_1}H_N|_{\gamma_2}&=2z_2^{\frac{2}{|m_1|}} A_{2/|m_1|}(0,c_2^2)=\sum_{3\leq 2k_2+2/|m_1|+2\leq N}2a_{0,k_2+\frac{2}{|m_1|},2,k_2}c_2^{2k_2} z_2^{\frac{2}{|m_1|}}.
\end{aligned}
\end{equation}
Denoting
\begin{equation}\label{equ: def of B}
B=e^{-2i\omega_2t}\partial_{\bar z_1\bar z_1}H_N|_{\gamma_2},\quad \partial_{\bar z_1\bar z_1}H_N|_{\gamma_2}=\tilde c_2e^{i\tilde \theta_2},\quad \tilde c_2\geq 0,
\end{equation}
where $\tilde c_2=2c_2^{2/|m_1|}|A_{2/|m_1|}(0,c_2^2)|$ and $\tilde \theta_2=\mathrm{arg} A_{2/|m_1|}(0,c_2^2)-\frac{2}{|m_1|}\omega_2t$. In particular, $\tilde c_2=0$, if $|m_1|>2$. The equation \eqref{thdot} can be rewritten as
\begin{equation}\label{equation of dot theta 2}
\begin{aligned}
\dot \theta(t)&=\omega_2+\hat \omega_1-2\mathrm{Re}(e^{2\theta i}B).
\end{aligned}
\end{equation}
Defining $\bar \theta(t):=\arg(e^{2\theta i}B)=2\theta(t)-(2+\frac{2}{|m_1|})\omega_2t+\arg A_{2/m_2}(0,c_2^2)$, we obtain
\begin{equation}\label{bartheta2}
\dot{\bar\theta}(t)=2\dot \theta(t)-\left(2+\frac{2}{|m_1|}\right)\omega_2=2\hat \omega_1-\frac{2}{|m_1|}\omega_2-4\tilde c_2\cos\bar\theta(t).
\end{equation}

Another way to obtain equations \eqref{bartheta} and \eqref{bartheta2} is as follows. If $m_2=2$, then $|m_1|>m_2$. 
Let $w_2(t) = r_2(t)e^{i\theta_2(t)}$ be a linearized solution along $\gamma_1$.  From the linearized equation \eqref{linw2b}, we obtain $\dot{\bar w}_2=i\hat \omega_2\bar w_2+2i \partial_{z_2z_2}H_N|_{\gamma_1}w_2$. Write $\partial_{z_2z_2}H_N|_{\gamma_1}:=\tilde c_1e^{-i\tilde \theta_1}$ as before. Recall that $\tilde \theta_1=\mathrm{arg}\bar A_1(c_1^2,0)-|m_1|\omega_1t$. Then we have
$$
\frac{\dot r_2}{r_2} - i \dot \theta_2(t) =\frac{\dot{\bar w}_2e^{i\theta_2}}{r_2}= i\hat \omega_2 +2i \tilde c_1 e^{2i\theta_2-i\tilde \theta_1},
$$
where $\hat \omega_2=\hat \omega_2(E),\omega_1=\omega_1(E)$ are given in \eqref{hatOmega2} and \eqref{Omega1}, respectively. Hence $\theta_2(t)$ satisfies
$$
-\dot \theta_2(t) = \hat \omega_2 +2\tilde c_1 \cos (-2\theta_2(t) - |m_1|\omega_1 t+\mathrm{arg}\bar A_1(c_1^2,0)), \quad \forall t.
$$
Finally, by defining
$\bar \theta(t) := -2\theta_2(t) - |m_1|\omega_1 t+\mathrm{arg}\bar A_1(c_1^2,0)+\pi,$
we obtain the same differential equation as \eqref{bartheta} for $\bar \theta(t)$.

If $m_2=1$ and condition \eqref{cond: gamma1 m2=1} holds, then $\gamma_1$ exists and from the linearized equation \eqref{linw2b}, we obtain
$$
-\dot \theta_2(t) = \hat \omega_2 +2\tilde c_1 \cos (-2\theta_2(t) - 2|m_1|\omega_1 t+\mathrm{arg}\bar A_2(c_1^2,0)), \quad \forall t.
$$
By defining
$\bar \theta(t) := -2\theta_2(t) - 2|m_1|\omega_1 t+\mathrm{arg}\bar A_2(c_1^2,0)+\pi,$
we end up with the same differential equation \eqref{bartheta} for $\bar \theta(t)$.

If $|m_1|\in\{1,2\}$, then $m_2=1$. Let $w_1(t)=r_1(t)e^{i\theta_1(t)}$ be a linearized solution along $\gamma_2(t)$. If $|m_1|=2$, then we obtain from \eqref{linw2a} that
$\dot {\bar w}_1=i\hat \omega_1\bar w_1+2i\partial_{z_1z_1}H_N|_{\gamma_2}w_1$. We write $\partial_{z_1z_1}H_N|_{\gamma_2}=\tilde c_2e^{-i\tilde \theta_2}$ as before. Recall that $\tilde \theta_2=\mathrm{arg} A_1(0,c_2^2)-\omega_2t$. It follows that
$$
\frac{\dot r_1}{r_1} - i \dot \theta_1(t) =\frac{\dot{\bar w}_1e^{i\theta_1}}{r_1}= i\hat \omega_1 +2i \tilde c_2 e^{2i\theta_1-i\tilde \theta_2},
$$
where $\hat \omega_1=\hat \omega_1(E),\omega_2=\omega_2(E)$ are given in \ref{hatOmega1} and \eqref{Omega2}, respectively. We conclude that $\theta_1(t)$ satisfies
$$
-\dot \theta_1(t) = \hat \omega_1 +2\tilde c_2 \cos (-2\theta_1(t) - \omega_2 t+\mathrm{arg} A_1(0,c_2^2)), \quad \forall t.
$$
Defining
$\bar \theta(t) := -2\theta_1(t) - \omega_2 t+\mathrm{arg} A_1(0,c_2^2)+\pi,$
we obtain \eqref{bartheta2}.

Similarly, if $|m_1|=1$ and condition \eqref{cond: gamma2 m2=1} holds, then the periodic orbit $\gamma_2$ exists and from the linearized equation \eqref{linw2a} we obtain
$$
-\dot \theta_1(t) = \hat \omega_1 +2\tilde c_2 \cos (-2\theta_1(t) - 2\omega_2 t+\mathrm{arg}A_2(0,c_2^2)), \quad \forall t.
$$
As before, defining
$\bar \theta(t) := -2\theta_1(t) - 2\omega_2 t+\mathrm{arg}\bar A_2(0,c_2^2)+\pi,$
we obtain \eqref{bartheta2}.

Using the linearized equations \eqref{bartheta} and \eqref{bartheta2} in the following lemma, we shall compute the rotation numbers of $\gamma_1$ and $\gamma_2$.

\begin{lem}\label{lem: asymptotic limit of theta}
Let $\theta:\R \to \R$ satisfy $\dot \theta(t)=a+b\cos \theta(t), \forall t,$ where $a,b \in \R$. Then the following statements hold.
\begin{itemize}
\item[(i)] If $|a|\leq |b|$, then $\lim_{t \to +\infty} \theta(t)/t=0$.
\item[(ii)] If $|a|>|b|$, then $\lim_{t \to +\infty} \theta(t)/t =\mathrm{sign}(a)\sqrt{a^2-b^2}$.
\end{itemize}
\end{lem}

\begin{proof}
If $|a| \leq |b|$, then $\dot \theta=0$ has infinitely many solutions, and thus, every solution of $\dot \theta=a+b\cos \theta$ is bounded. This implies $\theta(t)/t \to 0$ as $t \to +\infty.$ Item (i) follows.

Now suppose that $|a|> |b|.$  Then every solution $\theta(t)$ is periodic with period
$$
\tau={\rm sign}(a)\int_0^{2\pi}\frac{d\theta}{a+b\cos \theta}=\frac{2\pi}{\sqrt{a^2-b^2}}>0.
$$
Therefore, we may write $\theta(t)= \frac{2\pi t}{\tau}{\rm sign}(a)+g(t)$ for some $\tau$-periodic function $g(t)$. This implies
$\lim_{t\to +\infty}\theta(t)/t=\mathrm{sign}(a)\sqrt{a^2-b^2}.$
\end{proof}


Let us start with the rotation number $\rho_1=\rho_1(E)$ of $\gamma_1(t)=(c_1 e^{-i\omega_1 t},0)$, with $c_1>0$. If $m_2=2$, then $\gamma_1$ exists, and we compute $\rho_1$ using the linearized equation \eqref{bartheta}. If $m_2=1$, then $\gamma_1$ exists provided that \eqref{cond: gamma1 m2=1} holds, and in that case we also use \eqref{bartheta} to compute $\rho_1$.

Recall that the period of $\gamma_1$ is $T_1 = 2\pi / \omega_1$. Using \eqref{rho}, \eqref{bartheta} and Lemma \ref{lem: asymptotic limit of theta}, we compute $\rho_1$ for $m_2=1$ or $2$
\begin{equation}\label{rot of gamma1 p=2 1}
\begin{aligned}
\rho_1 & =\lim_{t\rightarrow +\infty}T_1\frac{\theta(t)}{2\pi t}=\lim_{t\rightarrow +\infty}\frac{\theta(t)}{\omega_1 t}= 1+\frac{|m_1|}{m_2} + \frac{1}{2\omega_1}\lim_{t \to +\infty} \frac{\bar \theta(t)}{t} \\
& =
\left\{
\begin{aligned}
& 1+\frac{|m_1|}{m_2}, & \mbox{ if } & \ \left|\frac{\hat \omega_2}{\omega_1}-\frac{|m_1|}{m_2}\right| \leq \frac{2\tilde c_1}{\omega_1}\\
& 1+\frac{\hat \omega_2}{\omega_1} -\Delta_1, & \mbox{ if } &  \  \frac{\hat \omega_2}{\omega_1}-\frac{|m_1|}{m_2} > \frac{2\tilde c_1}{\omega_1}\\
& 1+\frac{\hat \omega_2}{\omega_1} +\Delta_1, & \mbox{ if } &  \ \frac{|m_1|}{m_2}-\frac{\hat \omega_2}{\omega_1} > \frac{2\tilde c_1}{\omega_1},
\end{aligned}\right.
\end{aligned}
\end{equation}
where
$$
\tilde c_1=2c_1^{\frac{2|m_1|}{m_2}}|A_{2/m_2}(c_1^2,0)|=\bigg|\sum_{3\leq 2k_1+2+\frac{2|m_1|}{m_2}\leq N}2a_{k_1,2,k_1+\frac{2|m_1|}{m_2},0}\cdot c_1^{2k_1+\frac{2|m_1|}{m_2}}\bigg|,
$$
and
$$
\Delta_1=\frac{\left(\frac{2\tilde c_1}{\omega_1}\right)^2}{\sqrt{\left(\frac{|m_1|}{m_2}-\frac{\hat\omega_2}{\omega_1}\right)^2}+\sqrt{\left(\frac{|m_1|}{m_2}-\frac{\hat\omega_2}{\omega_1}\right)^2-\left(\frac{2\tilde c_1}{\omega_1}\right)^2}}.
$$
Let
$$
C_1:=\left(\frac{\hat\omega_2}{\omega_1}-\frac{|m_1|}{m_2}\right)^2-\left(\frac{2\tilde c_1}{\omega_1}\right)^2.
$$
Notice that $\rho_1=1+\frac{|m_1|}{m_2}\pm \sqrt{C_1}$ if $C_1\geq0$ and $\rho_1=1+\frac{|m_1|}{m_2}$ if $C_1<0$.

Let us estimate $C_1$ and $\Delta_1$ as functions of $E$.
In view of \eqref{omega/Omega} we have
\begin{equation}\label{square minus square 1}
\begin{aligned}
C_1
&=\left(\Omega_{\nu,1}\left(\frac{2}{\alpha_1}\right)^\nu E^{\nu-1}+O(E^\nu)\right)^2\\
&\quad \ -\left(\frac{4\big|\sum_{3\leq 2k_1+\frac{2|m_1|}{m_2}+2\leq N}a_{k_1,2,k_1+\frac{2|m_1|}{m_2},0}\cdot(\frac{2E}{\alpha_1}+O(E^\nu))^{k_1+\frac{|m_1|}{m_2}}\big|}{\alpha_1+2\nu a_{\nu,0,\nu,0}\left(\frac{2E}{\alpha_1}\right)^{\nu-1}+O\left(E^{\nu}\right)}\right)^2\\
&=\Omega_{\nu,1}^2\left(\frac{2}{\alpha_1}\right)^{2\nu}E^{2(\nu-1)}+O(E^{2\nu-1})-\frac{4^2}{\alpha_1^2}\big|a_{0,2,\frac{2|m_1|}{m_2},0}\big|^2\left(\frac{2E}{\alpha_1}\right)^{\frac{2|m_1|}{m_2}}-O\left(E^{\frac{2|m_1|}{m_2}+1}\right).
\end{aligned}
\end{equation}
where
$$
\Omega_{\nu,1}=a_{\nu-1,1,\nu-1,1}-\nu a_{\nu,0,\nu,0}\frac{|m_1|}{m_2}.
$$
We also have
\begin{equation}\label{Delta1_estimate}
\Delta_1=\frac{\frac{4^2}{\alpha_1^2}\big|a_{0,2,\frac{2|m_1|}{m_2},0}\big|^2\left(\left(\frac{2E}{\alpha_1}\right)^{\frac{2|m_1|}{m_2}}+O(E^{\frac{2|m_1|}{m_2}+1})\right)+O\left(E^{\frac{2|m_1|}{m_2}+2}\right)}{|\Omega_{\nu,1}|\frac{2^\nu}{\alpha_1^\nu}E^{\nu-1}+O(E^\nu)+\sqrt{\Omega_{\nu,1}^2\left(\frac{2}{\alpha_1}\right)^{2\nu}E^{2(\nu-1)}-\frac{4^2}{\alpha_1^2}\big|a_{0,2,\frac{2|m_1|}{m_2},0}\big|^2\left(\frac{2E}{\alpha_1}\right)^{\frac{2|m_1|}{m_2}}+R_1}},
\end{equation}
where $R_1=\mathcal{O}_{\min\{2\nu-1,2|m_1|/m_2+1\}}$.

To estimate the non-resonance condition, recall that
\begin{equation}\label{ohatdivo}
\frac{\hat \omega_2}{\omega_1}=\frac{|m_1|}{m_2}+\Omega_{\nu,1}\frac{2^\nu}{\alpha_1^\nu}E^{\nu-1}+O(E^\nu)\quad \mbox{ and } \quad \frac{\hat \omega_1}{\omega_2}=\frac{m_2}{|m_1|}+\Omega_{\nu,2}\frac{2^\nu}{\alpha_2^\nu}E^{\nu-1}+O(E^\nu),
\end{equation}
and thus
$$
\frac{\hat \omega_2\hat \omega_1}{\omega_1\omega_2}=1+2^\nu\Omega_\nu E^{\nu-1}+O(E^\nu),
$$
where
$$
\Omega_\nu = \frac{|m_1|}{m_2}\frac{\Omega_{\nu,2}}{\alpha_2^\nu} + \frac{m_2}{|m_1|}\frac{\Omega_{\nu,1}}{\alpha_1^\nu}.
$$

If $|m_1|\geq 3,$ the rotation number $\gamma_2(t)=(0,c_2e^{-i\omega_2 t})$ is computed as in \eqref{rho2 m2>2}, that is $\rho_2 = 1 + \hat \omega_1/\omega_2=1+m_2/|m_1| + \Omega_{\nu,2}(2/\alpha_2)^\nu E^{\nu-1} + O(E^{\nu})$.

We see from \eqref{rot of gamma1 p=2 1} that the following criteria hold for $m_2\in \{1,2\}$ and $|m_1| \geq 3$:

\begin{lem}\label{lem_criteria1} Assume that $m_2\in \{1,2\}$ and $|m_1|\geq 3$. Then the following assertions hold:
\begin{itemize}
\item[(i)] if $C_1\geq 0$,
then
$$
(\rho_1-1)(\rho_2-1)=\frac{\hat \omega_2 \hat \omega_1}{\omega_1 \omega_2}\pm \Delta_1 \frac{\hat \omega_1}{\omega_2}.
$$

\item[(ii)] if $C_1 < 0$, then
$$
(\rho_1-1)(\rho_2-1)=\frac{|m_1|}{m_2}\frac{\hat \omega_1}{\omega_2}.
$$
\end{itemize}
Here, $C_1=C_1(E)$ and $\Delta_1=\Delta_1(E)$ satisfy the estimates \eqref{square minus square 1} and \eqref{Delta1_estimate}.
\end{lem}

We use Lemma \ref{lem_criteria1} to estimate $(\rho_1-1)(\rho_2-1)$ as a function of $E>0$. Let us start with the case $m_2=2$. In particular, $|m_1|>2$. The following theorem gives conditions for the existence of non-resonant Hopf links for $H_N$.

\begin{thm}\label{thm: non-resonant theorem}
Assume that $m_2=2$ (in particular, $|m_1|\geq 3$), and let $\gamma_1,\gamma_2\subset \Sigma_E^N\subset H_N^{-1}(E)$ be the pair of periodic orbits of $H_N$ forming a Hopf link as above. The following statements hold:
\begin{itemize}

\item[(i)] If $|m_1|>2(\nu-1)$ and $\Omega_\nu,\Omega_{\nu,1} \neq 0$, then $\rho_1,\rho_2\notin \mathbb{Z}$ and
$$(\rho_1 -1)(\rho_2-1) = 1+2^{\nu} \Omega_\nu E^{\nu-1} + O(E^{\nu})\neq 1,$$ for every $E>0$ sufficiently small.

\item[(ii)] If $|m_1|=2(\nu-1)$, $\Omega_\nu,\Omega_{\nu,1} \neq 0$ and $a_{0,2,|m_1|,0}=0$, then  $\rho_1,\rho_2\notin \mathbb{Z}$ and
$$(\rho_1 -1)(\rho_2-1) = 1+2^\nu \Omega_\nu E^{\nu-1} + O(E^{\nu})\neq 1,$$ for every $E>0$ sufficiently small.

\item[(iii)] If $|m_1|<2(\nu-1)$,  $\Omega_{\nu,2}\neq 0$ and $a_{0,2,|m_1|,0}\neq 0$, then $\rho_1,\rho_2\notin \mathbb{Z}$ and
$$
(\rho_1 -1)(\rho_2-1) = 1+\frac{|m_1|}{2}\left(\frac{2}{\alpha_2}\right)^\nu \Omega_{\nu,2} E^{\nu-1} + O(E^\nu) \neq 1,
$$ for every $E>0$ sufficiently small.
\end{itemize}
In particular, in all cases (i), (ii), and (iii) above, the Hopf link $L=\gamma_1 \cup \gamma_2$ is non-resonant for every $E>0$ sufficiently small.
\end{thm}

\begin{proof}
Assume that conditions in (i) hold. Then $C_1 =\Omega_{\nu,1}^2\frac{2^{2\nu}}{\alpha_1^{2\nu}}E^{2(\nu-1)} + O(E^{2\nu-1})>0$ for every $E>0$ sufficiently small and
$\Delta_1 = O(E^{|m_1|-(\nu -1)})$. Hence $\Delta_1 = O(E^\nu)$ and by Lemma \ref{lem_criteria1}-(i) and \eqref{ohatdivo}, we obtain that $\rho_1=1+\frac{|m_1|}{2}+O(E^{\nu-1})$, $\rho_2=1+\frac{2}{|m_1|} + O(E^{\nu-1})$ are not integers and $(\rho_1 -1)(\rho_2-1) = 1+2^{\nu} \Omega_\nu E^{\nu-1} + O(E^{\nu})\neq 1$ for every $E>0$ sufficiently small. 

Assuming (ii), we see as in (i) that $C_1 =\Omega_{\nu,1}^2\frac{2^{2\nu}}{\alpha_1^{2\nu}}E^{2(\nu-1)} + O(E^{2\nu-1})>0$ for every $E>0$ sufficiently small, and $\Delta_1 = O(E^{|m_1|+1-(\nu-1)})=O(E^\nu)$. Hence, by Lemma \ref{lem_criteria1}-(i) and \eqref{ohatdivo}, $\rho_1,\rho_2$ admit the same estimate as before and $(\rho_1 -1)(\rho_2-1) = 1+2^\nu \Omega_\nu E^{\nu-1} + O(E^{\nu})\neq 1$ for every $E>0$ sufficiently small.

Assuming (iii), we see that $C_1=-\frac{4^2}{\alpha_1^2}|a_{0,2,|m_1|,0}|^2(\frac{2E}{\alpha_1})^{|m_1|}+O(E^{|m_1|+1})<0$ for every $E>0$ sufficiently small. It follows by Lemma \ref{lem_criteria1}-(ii) and \eqref{ohatdivo} that $\rho_1=1+\frac{|m_1|}{2}$, $\rho_2=1+\frac{2}{|m_1|} +\Omega_{\nu,2}( \frac{2}{\alpha_2})^\nu E^{\nu-1} +  O(E^{\nu})$ and $(\rho_1 -1)(\rho_2-1) = 1+\frac{|m_1|}{2}(\frac{2}{\alpha_2})^\nu \Omega_{\nu,2} E^{\nu-1} + O(E^\nu) \neq 1$ for every $E>0$ sufficiently small.
\end{proof}

Now, let us treat the case $m_2=1$. As in the case $m_2=2$, the rotation number of $\gamma_1$ is given by \eqref{rot of gamma1 p=2 1} provided that $\gamma_1$ exists, i.e., equation \eqref{cond: gamma1 m2=1} holds.
If $|m_1|\geq 3,$ then $\gamma_2$ always exists and the rotation number $\rho_2$ of  $\gamma_2(t)=(0,c_2 e^{-i\omega_2 t})$ is computed as in \eqref{rho2 m2>2}, that is $\rho_2 = 1 + \hat \omega_1/\omega_2$ as before. We obtain the following lemma.

\begin{thm}\label{thm: non-resonant theorem_A}
Assume that $m_2=1$, $|m_1|\geq 3$, and that the periodic orbit $\gamma_1 \subset \Sigma_E^N \subset H_N^{-1}(E)$ exists, i.e., \eqref{cond: gamma1 m2=1} holds for $H_N$. The following statements hold:
\begin{itemize}

\item[(i)] If $|m_1|>\nu-1$ and $\Omega_{\nu,1},\Omega_\nu\neq 0$, then $\rho_1,\rho_2\notin \mathbb{Z}$ and $$(\rho_1 -1)(\rho_2-1) = 1+2^\nu \Omega_\nu E^{\nu-1} + O(E^{\nu})\neq 1,$$ for every $E>0$ sufficiently small.

\item[(ii)] If $|m_1|=\nu-1$, $\Omega_{\nu,1},\Omega_\nu\neq 0$ and $a_{0,2,2|m_1|,0}=0$, then $\rho_1,\rho_2\notin \mathbb{Z}$ and $$(\rho_1 -1)(\rho_2-1) = 1+2^\nu \Omega_\nu E^{\nu-1} + O(E^{\nu})\neq 1,$$ for every $E>0$ sufficiently small.

\item[(iii)] If $|m_1|<\nu-1$,  $\Omega_{\nu,2}\neq 0$ and $a_{0,2,2|m_1|,0}\neq 0$, then $\rho_1\in \mathbb{Z},\rho_2\notin \mathbb{Z}$ and $$(\rho_1 -1)(\rho_2-1) = 1+|m_1|\left(\frac{2}{\alpha_2}\right)^\nu \Omega_{\nu,2} E^{\nu-1} + O(E^\nu) \neq 1,$$ for every $E>0$ sufficiently small.
\end{itemize}
In particular, in all cases (i), (ii), and (iii) above, the Hopf link $L=\gamma_1 \cup \gamma_2$ is non-resonant for every $E>0$ sufficiently small.
\end{thm}

\begin{proof}Assume that conditions in (i) hold. Then $C_1 =\Omega_{\nu,1}^2(\frac{2}{\alpha_1})^{2\nu}E^{2(\nu-1)} + O(E^{2\nu-1})>0$ for every $E>0$ sufficiently small and
$\Delta_1 = O(E^{2|m_1|-(\nu -1)})$. Hence $\Delta_1 = O(E^\nu)$ and by the proof of Theorem \ref{thm: non-resonant theorem}-(i) and \eqref{ohatdivo}, we obtain $\rho_1=1+|m_1|+\Omega_{\nu,1}\frac{2^\nu}{\alpha_1^\nu}E^{\nu-1}+O(E^{\nu})$, $\rho_2=1+\frac{1}{|m_1|} + O(E^{\nu-1})$ and $(\rho_1 -1)(\rho_2-1) = 1+2^\nu \Omega_\nu E^{\nu-1} + O(E^{\nu})\neq 1$ for every $E>0$ sufficiently small.

Assuming (ii), we see as in (i) that $C_1 =\Omega_{\nu,1}^2(\frac{2}{\alpha_1})^{2\nu}E^{2(\nu-1)} + O(E^{2\nu-1})>0$ for every $E>0$ sufficiently small, and $\Delta_1 = O(E^{2|m_1|+1-(\nu-1)})=O(E^\nu)$. By Theorem \ref{thm: non-resonant theorem}-(i) and \eqref{ohatdivo}, $\rho_1,\rho_2$ admit the same estimate as above and $(\rho_1 -1)(\rho_2-1) = 1+2^\nu \Omega_\nu E^{\nu-1} + O(E^{\nu})\neq 1$ for every $E>0$ sufficiently small.

Assuming (iii), we see that $C_1=-\frac{4^2}{\alpha_1^2}|a_{0,2,2|m_1|,0}|^2(\frac{2E}{\alpha_1})^{|m_1|}+O(E^{|m_1|+1})<0$ for every $E>0$ sufficiently small. It follows from Lemma \ref{lem_criteria1}-(ii) and \eqref{ohatdivo} that $\rho_1=1+|m_1|,\rho_2=1+\frac{1}{|m_1|}+\Omega_{\nu,2}(\frac{2}{\alpha_2})^\nu E^{\nu-1}   + O(E^{\nu})$ and $(\rho_1 -1)(\rho_2-1) = 1+|m_1|(\frac{2}{\alpha_2})^\nu \Omega_{\nu,2} E^{\nu-1} + O(E^\nu) \neq 1$ for every $E>0$ sufficiently small. \end{proof}

We now assume that $m_2=1$ and $|m_1|\in \{1,2\}$. We keep the assumption that $\gamma_1$ exists, i.e., equation \eqref{cond: gamma1 m2=1} holds. If $|m_1|=2$, then $
\gamma_2$ exists. If $|m_1|=1$, we assume that \eqref{cond: gamma2 m2=1} also holds and $\gamma_2$ exists. We consider the linearized equation \eqref{linw2a} to compute $\rho_2$. Recall that the period of $\gamma_2$ is $T_2=2\pi/\omega_2$.

Using \eqref{rho}, \eqref{bartheta2} and Lemma \ref{lem: asymptotic limit of theta}, we obtain
\begin{equation}\label{rot of gamma2 q=2 1}
\begin{aligned}
\rho_2 & =\lim_{t\rightarrow +\infty}T_2\frac{\theta(t)}{2\pi t}=\lim_{t\rightarrow +\infty}\frac{\theta(t)}{\omega_2 t}= 1+\frac{1}{|m_1|} + \frac{1}{2\omega_2}\lim_{t \to +\infty} \frac{\bar \theta(t)}{t} \\
& =
\left\{
\begin{aligned}
& 1+\frac{1}{|m_1|}, & \mbox{ if } & \ \left|\frac{\hat \omega_1}{\omega_2}-\frac{1}{|m_1|}\right| \leq \frac{2\tilde c_2}{\omega_2}\\
& 1+\frac{\hat \omega_1}{\omega_2}-\Delta_2, & \mbox{ if } &  \ \frac{\hat \omega_1}{\omega_2}-\frac{1}{|m_1|} > \frac{2\tilde c_2}{\omega_2}\\
& 1+\frac{\hat \omega_1}{\omega_2}+\Delta_2, & \mbox{ if } &  \ \frac{1}{|m_1|}-\frac{\hat \omega_1}{\omega_2} >\frac{2\tilde c_2}{\omega_2},
\end{aligned}\right.
\end{aligned}
\end{equation}
where
$$
\tilde c_2=2c_2^{\frac{2}{|m_1|}}|A_{2/|m_1|}(0,c_2^2)|=\bigg|\sum_{3\leq 2k_2+\frac{2}{|m_1|}+2\leq N}2a_{0,k_2+\frac{2}{|m_1|},2,k_2}\cdot c_2^{2k_2+\frac{2}{|m_1|}}\bigg|
\geq 0$$  and
$$
\Delta_2=\frac{\left(\frac{2\tilde c_2}{\omega_2}\right)^2}{\sqrt{\left(\frac{1}{|m_1|}-\frac{\hat\omega_1}{\omega_2}\right)^2}+\sqrt{\left(\frac{1}{|m_1|}-\frac{\hat\omega_1}{\omega_2}\right)^2-\left(\frac{2\tilde c_2}{\omega_2}\right)^2}}.
$$
Let
$$
C_2 := \left(\frac{\hat \omega_1}{\omega_2} - \frac{1}{|m_1|}\right)^2 - \left(\frac{2\tilde c_2}{\omega_2}\right)^2.$$
Notice that if $C_2\geq 0$, then $\rho_2=1+\frac{1}{|m_1|}\pm \sqrt{C_2}$. Otherwise, $\rho_2=1+\frac{1}{|m_1|}$.

If $m_2=1$ and $|m_1|\in\{1,2\}$, then \eqref{omega/Omega} implies
\begin{equation}\label{square minus square 2}
\begin{aligned}
C_2
&=\left(\Omega_{\nu,2}\left(\frac{2}{\alpha_2}\right)^{\nu}E^{\nu-1}+O(E^\nu)\right)^2\\
&\quad\ -\left(\frac{4\big|\sum_{3\leq 2k_2+\frac{2}{|m_1|}+2\leq N}a_{0,k_2+\frac{2}{|m_1|},2,k_2}\cdot\left(\frac{2E}{\alpha_2}+O(E^\nu)\right)^{k_2+\frac{1}{|m_1|}}\big|}{\alpha_2+2\nu a_{0,\nu,0,\nu}\left(\frac{2E}{\alpha_2}\right)^{\nu-1}+O(E^{\nu})}\right)^2\\
&=\Omega_{\nu,2}^2\left(\frac{2}{\alpha_2}\right)^{2\nu}E^{2(\nu-1)}+O(E^{2\nu-1})-\frac{4^2}{\alpha_2^2}\big|a_{0,\frac{2}{|m_1|},2,0}\big|^2\left(\frac{2E}{\alpha_2}\right)^{\frac{2}{|m_1|}}-O\left(E^{\frac{2}{|m_1|}+1}\right).
\end{aligned}
\end{equation}
where $\Omega_{\nu,2}=a_{1,\nu-1,1,\nu-1}-\nu a_{0,\nu,0,\nu}\frac{1}{|m_1|}$.
From \eqref{square minus square 1} and \eqref{square minus square 2}, we conclude that
\begin{equation}\label{Delta2_estimate}
\Delta_2=\frac{\frac{4^2}{\alpha_2^2}\big|a_{0,\frac{2}{|m_1|},2,0}\big|\left(\left(\frac{2E}{\alpha_2}\right)^{\frac{2}{|m_1|}}+O\left(E^{\frac{2}{|m_1|}+1}\right)\right)+O\left(E^{\frac{2}{|m_1|}+2}\right) }{|\Omega_{\nu,2}|\frac{2^{\nu}}{\alpha_2^{\nu}}E^{\nu-1}+O(E^\nu)+\sqrt{\Omega_{\nu,2}^2\left(\frac{2}{\alpha_2}\right)^{2\nu}E^{2(\nu-1)}-\frac{4^2}{\alpha_2^2}\big|a_{0,\frac{2}{|m_1|},2,0}\big|^2\left(\frac{2E}{\alpha_2}\right)^{\frac{2}{|m_1|}}+R_2}}.
\end{equation}
where $R_2=\mathcal{O}_{\min\{2(\nu-1),2/|m_1|\}+1}$.

Hence, if $m_2=1$ and $m_1\in \{1,2\}$, we obtain from \eqref{rot of gamma1 p=2 1} and \eqref{rot of gamma2 q=2 1} the following lemma.

\begin{lem}\label{lem_criteria2} Assume that $|m_1| \in \{1,2\}$ (in particular, $m_2=1$). Then, the following assertions hold:
\begin{itemize}
\item[(i)] if $C_1\geq 0$ and $C_2\geq 0,$ then
$$
\begin{aligned}
(\rho_{1}-1)(\rho_{2}-1)&=\frac{\hat \omega_2\hat \omega_1}{\omega_1\omega_2}\pm \frac{\hat \omega_1}{\omega_2}\Delta_1\pm \frac{\hat \omega_2}{\omega_1}\Delta_2 \pm \Delta_1\Delta_2.
\end{aligned}
$$

\item[(ii)] if $C_1<0$ and $C_2 \geq 0$, then
$$
\begin{aligned}
(\rho_{1}-1)(\rho_{2}-1)&=|m_1|\cdot\frac{\hat \omega_1}{\omega_2}\pm|m_1|\Delta_2.
\end{aligned}
$$

\item[(iii)] if $C_1\geq 0$ and $C_2< 0,$ then
$$
\begin{aligned}
(\rho_{1}-1)(\rho_{2}-1)&=\frac{\hat \omega_2}{\omega_1}\cdot \frac{1}{|m_1|}\pm\frac{1}{|m_1|}\Delta_1.
\end{aligned}
$$

\item[(iv)] if $C_1< 0$ and $C_2< 0,$ then
$$
\begin{aligned}
(\rho_{1}-1)(\rho_{2}-1)&=|m_1|\cdot \frac{1}{|m_1|}=1.
\end{aligned}
$$
\end{itemize}
Here, $C_i=C_i(E), \Delta_i=\Delta_i(E), i=1,2,$ satisfy thee estimates \eqref{square minus square 1}, \eqref{Delta1_estimate}, \eqref{square minus square 2} and \eqref{Delta2_estimate}.
\end{lem}

The following theorem provides conditions for the link $\gamma_1 \cup \gamma_2$ to be non-resonant.

\begin{thm}\label{thm: non-resonant theorem B}
Assume that $m_2=1$, $|m_1|=2$, and that the periodic orbit  $\gamma_1 \subset \Sigma_E^N \subset H_N^{-1}(E)$ exists, i.e., \eqref{cond: gamma1 m2=1} holds for $H_N$. The following statements hold:
\begin{itemize}

\item[(i)] If $\nu=2$, $\Omega_{2,1}\neq 0$ and $a_{0,1,2,0}\neq 0$, then $\rho_1\notin\mathbb{Z},\rho_2=3/2$ and $$(\rho_1 -1)(\rho_2-1) =1+ \frac{2}{\alpha_1^2}\Omega_{2,1} E + O(E^2) \neq 1,$$ for every $E>0$ sufficiently small.

\item[(ii)] If $\nu=2$, $\Omega_{2,1}, \Omega_{2}\neq 0$, then $\rho_1\notin\mathbb{Z},\rho_2=3/2+O(E)$ and $$(\rho_1 -1)(\rho_2-1) =1+  cE + O(E^2)\neq 1,$$ for every $E>0$ sufficiently small. Here, $c\neq 0$ is a constant independent of $E$.

\item[(iii)] If $\nu=3$, $\Omega_{3,1}\neq 0$, $a_{0,1,2,0}\neq 0$ and $a_{0,2,4,0}=0$, then $\rho_1\notin \mathbb{Z},\rho_2=3/2$ and
$$(\rho_1-1)(\rho_2-1) = 1+ \frac{2^2}{\alpha_1^3} \Omega_{3,1} E^2 + O(E^3) \neq 1,$$ for every $E>0$ sufficiently small.

\end{itemize}
In particular, in all cases (i), (ii), and (iii) above, the Hopf link $L=\gamma_1 \cup \gamma_2$ is non-resonant for every $E>0$ sufficiently small.
\end{thm}

\begin{proof}
Assuming the conditions in (i), we see that $C_1 =\Omega_{2,1}^2\frac{2^4}{\alpha_1^4}E^2 + O(E^3)>0$, for every $E>0$ sufficiently small and
$\Delta_1 = O(E^3)$. We also see that $C_2 =-\frac{2^5}{\alpha_2^3}|a_{0,1,2,0}|^2E + O(E^2)<0$ for every $E>0$ sufficiently small. From Lemma \ref{lem_criteria2}-(iii), \eqref{rot of gamma1 p=2 1} and \eqref{rot of gamma2 q=2 1}, we obtain $\rho_1=3+\Omega_{2,1}\frac{2^2}{\alpha_1^2}E+O(E^2)$, $\rho_2=3/2$ and $(\rho_1 -1)(\rho_2-1) =1+ \frac{2}{\alpha_1^2}\Omega_{2,1} E + O(E^2) \neq 1$ for every $E>0$ sufficiently small.

Assuming the conditions in (ii), we see that $C_1 =\Omega_{2,1}^2\frac{2^4}{\alpha_1^4}E^2 + O(E^3)>0$, for every $E>0$ sufficiently small and
$\Delta_1 = O(E^3)$. We also see that $C_2 =-\frac{2^5}{\alpha_2^3}|a_{0,1,2,0}|^2E + O(E^2)$. If $C_2<0$ for every $E>0$ sufficiently small, then by Lemma \ref{lem_criteria2}-(iii), \eqref{rot of gamma1 p=2 1} and \eqref{rot of gamma2 q=2 1}, we obtain the same estimates of $\rho_1,\rho_2$ as above and $(\rho_1 -1)(\rho_2-1) =1+ \frac{2}{\alpha_1^2}\Omega_{2,1} E + O(E^2) \neq 1$ for every $E>0$ sufficiently small. If $C_2\geq0$ for every $E>0$ sufficiently small, then $a_{0,1,2,0}$ must be zero and we have $\Delta_2=O(E^2)$. Then by Lemma \ref{lem_criteria2}-(i), \eqref{rot of gamma1 p=2 1} and \eqref{rot of gamma2 q=2 1}, we obtain $\rho_1=3+\Omega_{2,1}\frac{2^2}{\alpha_1^2}E+O(E^2)$, $\rho_2=3/2+\Omega_{2,2}\frac{2^2}{\alpha_2^2}E+O(E^2)$ and $(\rho_1 -1)(\rho_2-1) =1+  4\Omega_{2} E + O(E^2)\neq 0$ for every $E>0$ sufficiently small.

Under the conditions in (iii), we have $C_1=\Omega_{3,1}^2 \frac{2^6}{\alpha_1^6}E^4 + O(E^5) >0$ for every $E>0$ sufficiently small and $\Delta_1 = O(E^4)$. Moreover, $C_2 = -\frac{2^5}{\alpha_2^3}|a_{0,1,2,0}|^2 E+ O(E^2)<0$ for every $E>0$ sufficiently small, then by Lemma \ref{lem_criteria2}-(iii), \eqref{rot of gamma1 p=2 1} and \eqref{rot of gamma2 q=2 1}, we obtain $\rho_1=3+\Omega_{3,1}\frac{2^3}{\alpha_1^3}E^2+O(E^3)$, $\rho_2=3/2$ and $(\rho_1-1)(\rho_2-1) = 1+ \frac{2^2}{\alpha_1^3} \Omega_{3,1} E^2 + O(E^3) \neq 1$ for every $E>0$ sufficiently small.
\end{proof}

Finally, we consider the most degenerate case $m_2=|m_1|=1$, that is $\alpha_1=\alpha_2$. We restrict to $\nu=2$ and consider both cases $\Omega_2\neq 0$ and $\Omega_2 =0.$ If $\Omega_2=0$, we need to consider the second order term $E^2$ in the expansion of $(\rho_1-1)(\rho_2-1)$, and some higher order coefficients $a_{kl}$ come into play. To do that, recall that $\gamma_1(t) = (c_1 e^{-i \omega_1 t},0)$ and $\gamma_2(t) = (0, c_2 e^{-i \omega_2 t}),$ where $c_1=c_1(E)>0$ and $c_2=c_2(E)>0$ satisfy
$$
\begin{aligned}
    E = \frac{\alpha_1}{2} c_1^2+\sum_{3\leq k_1\leq N} a_{k_1,0,k_1,0}c_1^{2k_1}, \quad
     E = \frac{\alpha_1}{2} c_2^2+\sum_{3\leq k_2\leq N} a_{0,k_2,0,k_2}c_2^{2k_2}.
\end{aligned}
$$
This gives
$$
\begin{aligned}
    c_1^2 = \frac{2E}{\alpha_1} - a_{2,0,2,0}\frac{2^3}{\alpha_1^3}E^2 + O(E^3), \quad
     c_2^2 = \frac{2E}{\alpha_1} - a_{0,2,0,2}\frac{2^3}{\alpha_1^3}E^2 + O(E^3).
\end{aligned}
$$
Now recall that
\begin{alignat*}{2}
\hat \omega_2 & = \alpha_2  + 2 \sum_{3\leq 2k_1\leq N}  a_{k_1-1,1,k_1-1,1} c_1^{2k_1-2}, \quad
&\hat \omega_1 & = \alpha_1 + 2 \sum_{3\leq 2k_2\leq N}  a_{1,k_2-1,1,k_2-1} c_2^{2k_2-2},
\\
\omega_1 & = \alpha_1 + 2 \sum_{3\leq 2k_1\leq N} k_1 a_{k_1,0,k_1,0} c_1^{2k_1-2},
\quad  &\omega_2 & = \alpha_2  + 2 \sum_{3\leq 2k_2\leq N} k_2 a_{0,k_2,0,k_2} c_2^{2k_2-2}.
\end{alignat*}
Hence
$$
\begin{aligned}
\frac{\hat \omega_2 }{\omega_1} =
1 + \Omega_{2,1}\frac{2^2}{\alpha_1^2}E + \frac{2^3}{\alpha_1^4}\beta_1 E^2 + O(E^3),\quad
\frac{\hat \omega_1}{\omega_2} = 1 + \Omega_{2,2}\frac{2^2}{\alpha_1^2}E + \frac{2^3}{\alpha_1^4} \beta_2 E^2 + O(E^3),
\end{aligned}
$$
where
\begin{equation}\label{beta1,2}
\begin{aligned}
\beta_1&=6 (2a_{2,0,2,0} - a_{1,1,1,1})a_{2,0,2,0} + \alpha_1(a_{2,1,2,1} - 3 a_{3,0,3,0}),\\
\beta_2&=6 (2a_{0,2,0,2} -
a_{1,1,1,1})a_{0,2,0,2} +\alpha_1 (a_{1,2,1,2} - 3 a_{0,3,0,3}).
\end{aligned}
\end{equation}
We obtain
\begin{equation}\label{o2o1/w1w2}
\frac{\hat \omega_2\hat \omega_1}{\omega_1\omega_2}=1+\frac{2^2}{\alpha_1^2}(\Omega_{2,1}+\Omega_{2,2})E+\frac{2^3}{\alpha_1^4}(\beta_1+\beta_2+2\Omega_{2,1}\Omega_{2,2})E^2+O(E^3).
\end{equation}

\begin{thm}\label{thm: non-resonant theorem C}
Assume that $\alpha_2=\alpha_1$, i.e., $m_2=|m_1|=1$, and that the periodic orbits $\gamma_1, \gamma_2 \subset \Sigma_E^N \subset H_N^{-1}(E)$ exist, i.e., \eqref{cond: gamma2 m2=1} and \eqref{cond: gamma1 m2=1}  hold. The following assertions hold:
\begin{itemize}
\item[(i)] If $\nu=2$, $\Omega_2\neq 0$ and $a_{0,2,2,0}=0$, then  $(\rho_1-1)(\rho_2-1) =1+cE + O(E^2) \neq 1,$ for every $E>0$ sufficiently small. Here, $c\neq 0$ is a constant independent of $E$. Moreover, if $\Omega_{2,1},\Omega_{2,2}\neq 0$, then $\rho_1,\rho_2\notin \mathbb{Z}$.

\item[(ii)] If $\nu=2$, $\Omega_{2,1}=-\Omega_{2,2}\neq 0$, $a_{0,2,2,0}=0$ and
$\beta_1+\beta_2+2\Omega_{2,1}\Omega_{2,2}\neq 0$, then $$
(\rho_1-1)(\rho_2-1)=1+\frac{2^3}{\alpha_1^4}(\beta_1+\beta_2+2\Omega_{2,1}\Omega_{2,2})E^2+O(E^3)\neq 1,
$$
for every $E>0$ sufficiently small.
\end{itemize}
In particular, in both cases (i) and (ii) above, the Hopf link $L=\gamma_1 \cup \gamma_2$ is non-resonant for every $E>0$ sufficiently small.
\end{thm}

\begin{proof}
Let us assume that conditions in (i) hold. Since $\Omega_2\neq0$,  at least one of $\Omega_{2,1}$ and $\Omega_{2,2}$ does not vanish. By symmetry, we may assume that $\Omega_{2,1}\neq 0$. Then  $C_1 = \Omega_{2,1}^2 \frac{2^4}{\alpha_1^4}E^2+O(E^3)>0$ for every $E>0$ sufficiently small. We also obtain $\Delta_1 = O(E^3).$ In the same way, if $\Omega_{2,2}\neq 0$, then $C_2= \Omega_{2,2}^2 \frac{2^4}{\alpha_2^4} E^2 + O(E^3)>0$ for every $E>0$ sufficiently small and $\Delta_2 = O(E^3).$ Then Lemma \ref{lem_criteria2}-(i) implies that $(\rho_1-1)(\rho_2-1) = 1 +4\Omega_2 E + O(E^2) \neq 1$ for every $E>0$ sufficiently small. If $\Omega_{2,2}=0$, then $C_2$ could be non-negative or negative, and $\hat \omega_1/\omega_2 = 1 + O(E^2)$. If $C_2\geq 0$ for every $E>0$ sufficiently small, then $(1- \hat\omega_1/\omega_2)^2\geq(2\tilde c_2/\omega_2)^2$, which implies that $\tilde c_2=O(E^2)$ and thus $|\Delta_2| \leq |2\tilde c_2/\omega_2| = O(E^2)$.
Lemma \ref{lem_criteria2}-(i) implies that $(\rho_1-1)(\rho_2-1) = 1 +4\Omega_2 E + O(E^2) \neq 1$ for every $E>0$ sufficiently small. If $C_2< 0$ for every $E>0$ sufficiently small, then we obtain from Lemma \ref{lem_criteria2}-(iii) that $(\rho_1-1)(\rho_2-1) = 1+ \Omega_{2,1}\frac{2^2}{\alpha_1^2}E + O(E^2)\neq 1$ for every $E>0$ sufficiently small.

Under the conditions in (ii), we have
$C_1 = \Omega_{2,1}^2 \frac{2^4}{\alpha_1^4}E^2+O(E^3)>0$ and $C_2= \Omega_{2,2}^2 \frac{2^4}{\alpha_2^4} E^2 + O(E^3)>0$, for every $E>0$ sufficiently small. We also obtain that $\Delta_1 = O(E^3)$ and $\Delta_2 = O(E^3)$. Therefore, we obtain
$$
(\rho_1-1)(\rho_2-1)=\frac{\hat \omega_2\tilde \omega_1}{\omega_1\omega_2}+O(E^3)=1+\frac{2^3}{\alpha_1^4}(\beta_1+\beta_2+2\Omega_{2,1}\Omega_{2,2})E^2+O(E^3)\neq 1,
$$
for every $E>0$ sufficiently small.
\end{proof}

\section{Non-resonant Hopf links for \texorpdfstring{$H$}{Lg}}\label{sec 4}
In this section, we show that if $H_N$ satisfies certain weak transversality conditions, then its non-resonant Hopf link can be continued to $H$ for every $E>0$ sufficiently small.

Write $y=(y_1,y_2)$ and $x=(x_1,x_2)$, and consider the Hamiltonian $H(y,x)=H_2+\mathcal{O}_3$ as in \eqref{H1}. Denote by $\Sigma_{E} \subset H^{-1}(E)$ the sphere-like component near the origin, for every $E>0$ small. Let
$$\tilde H(y,x):=\epsilon^{-2} H(\epsilon y,\epsilon x)=\frac{\alpha_1}{2}(x_1^2+y_1^2)+\frac{\alpha_2}{2}(x_2^2+y_2^2)+\epsilon^{-2}\mathcal{O}_{3}(\epsilon y,\epsilon x).$$
We omit $\epsilon$ in the notation of $\tilde H$ for simplicity. Since $\epsilon^{-2}\mathcal{O}_{3}(\epsilon y,\epsilon x)=O(\epsilon)$, the re-scaled sphere-like component $$\tilde \Sigma_{\epsilon^2}:=\epsilon^{-1}\Sigma_{\epsilon^2}\subset \tilde H^{-1}(1)=\epsilon^{-1}H^{-1}(\epsilon^2),$$  converges in $C^\infty$ to the ellipsoid $H_2^{-1}(1)$.
Let
$$
\gamma_{1,0}(t):=\sqrt{\frac{2}{\alpha_1}}(\sin(\alpha_1 t),0,\cos(\alpha_1 t),0),\quad \gamma_{2,0}(t):=\sqrt{\frac{2}{\alpha_2}}(0,\sin(\alpha_2 t),0,\cos(\alpha_2 t)),
$$
be the simple periodic orbits in $H_2^{-1}(1)$ forming a Hopf link. Their respective periods are denoted by $T_{1,0}$ and $T_{2,0}$. Let $P_j, j=1,2,$ be the Poincar\'e map given by the first return map to a surface of section $\chi_j\subset H_2^{-1}(1)$ transverse to $\gamma_{j,0}(0)$. The linear maps $DP_1(\gamma_{1,0}(0))$ and $DP_2(\gamma_{2,0}(0))$ possess eigenvalues $e^{\pm i\alpha_2/\alpha_1}$ and $e^{\pm i\alpha_1/\alpha_2}$, respectively. If $m_2\geq 2$ (i.e., $\alpha_2$ is not a trivial multiple of $\alpha_1$), then both $\gamma_{1,0}$ and $\gamma_{2,0}$ are non-degenerate. Hence, by the implicit function theorem, there exist families of periodic orbits $\gamma_{1,\epsilon}, \gamma_{2,\epsilon}$, where $\epsilon>0$ is small.

\subsection{Weakly non-resonant Hamiltonians}\label{sec: weakly nonresonant Hamiltonians}
Let us start considering the case $|m_1|>m_2\geq 2.$ Fix $N>0$ and denote $H=H_N+R_{N+1}$, where $H_N$ is in Birkhoff-Gustavson normal form and $R_{N+1}=\mathcal{O}_{N+1}$. Consider the $(\epsilon,\delta)$-dependent family of Hamiltonians defined by
\begin{equation}\label{equ: tilde H}
\tilde H_\delta:=\tilde H_N+\delta^{N-1}\tilde R,\quad \delta\in[0,\epsilon],
\end{equation}
where $\tilde H_N=\epsilon^{-2}H_N(\epsilon\cdot,\epsilon\cdot)$ and $\tilde R=\epsilon^{-(N+1)}R_{N+1}(\epsilon\cdot,\epsilon\cdot)$ are smooth $\epsilon$-dependent functions. Here, $(\epsilon,\delta)$ lies in a small neighborhood of $(0,0)\in \R^2$.

\begin{lem}\label{lem: perturbation}
Let $\tilde H_\delta$ be the family in \eqref{equ: tilde H}. Then there exist $M>0$ and smooth families of periodic orbits $\gamma^\delta_{j,\epsilon}\in \tilde H_{\delta}^{-1}(1), j=1,2,$ with period $T^\delta_{j,\epsilon}$, defined for every $(\epsilon,\delta)$ sufficiently close to $(0,0)$, such that
\begin{equation}\label{estimates_gamma10}
{\rm dist}(\gamma^\delta_{j,\epsilon}, \gamma^0_{j,\epsilon})\leq M \delta^{N-1}
,\quad
|T^\delta_{j,\epsilon}-T^{0}_{j,\epsilon}|\leq M\delta^{N-1},\quad |\rho^\delta_{j,\epsilon}-\rho^{0}_{j,\epsilon}|\leq M \delta^{(N-1)/2},
\end{equation}
where $\rho^\delta_{j,\epsilon}$ is the rotation number of $\gamma^\delta_{j,\epsilon}$. Here, $\gamma_{1,\epsilon}^0, \gamma_{2,\epsilon}^0\subset \tilde H_N^{-1}(1)$ are periodic orbits of $\tilde H_N$, contained in $\{z_2=0\}$ and $\{z_1=0\},$ respectively, and forming a Hopf link. Moreover, $\gamma_{j,\epsilon}^\delta \to \gamma_{j,0}$ in $C^\infty$ as $(\epsilon,\delta) \to (0,0).$
\end{lem}

\begin{proof}
Since $\gamma_{1,0},\gamma_{2,0}$ are nondegenerate,  there exist smooth families of periodic orbit $\gamma^\delta_{j,\epsilon}\in \tilde H_\delta^{-1}(1), j=1,2,$ near $\gamma_{j,0}$ with period $T^\delta_{j,\epsilon}\approx T_{j,0}$, smoothly depending on $(\epsilon,\delta^{N-1})$, for $(\epsilon,\delta^{N-1})$ sufficiently close to $(0,0)$.
 From the smooth dependence of solutions with respect to initial conditions and parameters, we obtain the first two uniform estimates in \eqref{estimates_gamma10} for some $M>0$.

The coefficients of the linearized solutions along $\gamma^\delta_{j,\epsilon}$ also smoothly depend on the initial conditions and parameters. The rotation number of $\gamma_{j,\epsilon}^\delta$ depends on the phase of the Floquet multipliers along $\gamma_{j,\epsilon}^\delta$. It is a smooth function on $(\epsilon,\delta^{N-1})$ if the multipliers are not $\pm 1$. Because we are assuming that $\alpha_2$ is not a multiple of $\alpha_1$, we know that $|m_1|> m_2\geq 2$. If $m_2=2$,  the Floquet multipliers along $\gamma_{1,0}$ are $-1$ with multiplicity $2$. In this case, the dependence of the multipliers has order $(N-1)/2$ due to the square root on the coefficients of the fundamental
 matrix of the linear system. If $m_2>2$, then this dependence has order $\delta^{N-1}$. Since $|m_1|>2$, we necessarily have the dependence of the multipliers along $\gamma_{2,\epsilon}^\delta$ has order $\delta^{N-1}$. In any case, taking $M>0$ sufficiently large, we have  $|\rho_{j,\epsilon}^\delta-\rho_{j,\epsilon}^0|\leq M \delta^{(N-1)/2}, j=1,2.$
\end{proof}

Based on the previous lemma, we obtain the following result for $m_2>2$.

\begin{thm} \label{thm: m2>2}
Assume that $|m_1|>m_2>2$ and that the Birkhoff-Gustavson normal form $H=H_N + \mathcal{O}_{N+1}$ satisfies $\Omega_\nu \neq0,$ see \eqref{equ: omega nu}, where $2\leq \nu\leq \lfloor N/2\rfloor$ is the smallest integer so that at least one of the numbers in \eqref{condition of nv} does not vanish. Then the sphere-like component $\Sigma_E\subset H^{-1}(E)$ carries a non-resonant Hopf link for every $E>0$ sufficiently small.
\end{thm}

\begin{proof}
It is enough to prove that $\tilde \Sigma_E=E^{-1}\Sigma_E \subset \tilde H^{-1}(1)$ contains a non-resonant Hopf link for every $E=\epsilon^2>0$ sufficiently small, where $\tilde H=\epsilon^{-2}H(\epsilon\cdot ,\epsilon \cdot)$ as before. The inductive construction of $H_N$ in section \ref{sec: construction} implies that $H_N$ satisfies the condition $\Omega_\nu\neq 0$ for every $N$ sufficiently large, where both $\nu\geq 2$ and $\Omega_\nu$ do not depend on $N$. 

Since $\gamma_{1,0},\gamma_{2,0}\subset H_2^{-1}(1)$ are nondegenerate, $\tilde H_N^{-1}(1)$ and $\tilde \Sigma_{E}$ admit respective pairs of periodic orbits $\gamma_1,\gamma_2$ and $\gamma_{1,E}, \gamma_{2,E}$ forming a Hopf link for every $E>0$ sufficiently small. Let $\rho_i$ and $\rho_{i,E}$ be the rotation number of $\gamma_i$ and $\gamma_{i,E}$, respectively. Since $|m_1|>2$, we obtain from \eqref{rho1 m2>2} and \eqref{rho2 m2>2}
$$
\rho_1=1+\frac{\hat \omega_2(E)}{\omega_1(E)},\quad \rho_2=1+\frac{\hat \omega_1(E)}{\omega_2(E)}.
$$
It follows from \eqref{equ: nonresonancy for m2>2} that
$$(\rho_1-1)(\rho_2-1)=\frac{\hat \omega_1\hat \omega_2}{\omega_1\omega_2}=1+2^\nu\Omega_\nu E^{\nu-1}+O(E^{\nu}).$$
By Lemma \ref{lem: perturbation}, we$|\rho_{i,E}-\rho_i|\leq M E^{(N-1)/4}$ for some constant $M$. Since $N$ can be taken arbitrarily large and $\Omega_\nu \neq 0$, we obtain $(\rho_{1,E}-1)(\rho_{2,E}-1)=1+2^\nu\Omega_\nu E^{\nu-1} + O(E^\nu) \neq 1$ for every $E>0$ sufficiently small. Hence, the Hopf link $\gamma_{1,E} \cup \gamma_{2,E}\subset \tilde \Sigma_E$ is non-resonant.
\end{proof}

\begin{thm} \label{thm: m2=2}
Assume that $|m_1|>m_2=2$ and that the Birkhoff-Gustavoson normal form $H=H_N+\mathcal{O}_{N+1}$ satisfies one of the following conditions:
\begin{itemize}
\item[(i)] $2(\nu-1)<|m_1|$ and $\Omega_{\nu,1},\Omega_\nu \neq 0$  see \eqref{omega/Omega} and \eqref{equ: omega nu}.
\item[(ii)] $2(\nu-1)=|m_1|$ and $\Omega_{\nu,1},\Omega_\nu \neq 0$ and $a_{0,2,|m_1|,0}=0$.
\item[(iii)] $2(\nu-1)>|m_1|$ and $\Omega_{\nu,2},a_{0,2,|m_1|,0}\neq 0$.
\end{itemize}
Here, $2\leq \nu\leq \lfloor N/2\rfloor$ is the smallest integer so that at least one of the numbers in \eqref{condition of nv} does not vanish. Then the sphere-like component $\Sigma_E\subset H^{-1}(E)$ carries a non-resonant Hopf link for every $E>0$ sufficiently small.
\end{thm}

\begin{proof}
It is enough to show that $\tilde \Sigma_E\subset \tilde H^{-1}(1)$ contains a non-resonant Hopf link for every $E=\epsilon^2>0$ sufficiently small, where $\tilde H=\epsilon^{-2}H(\epsilon\cdot ,\epsilon \cdot)$ as before. We may assume that $N$ is arbitrarily large.   Recall that $\nu, \Omega_\nu,\Omega_{\nu,1},\Omega_{\nu,2}$ do not depend on $N$ for $N$ sufficiently large.

Since $\gamma_{1,0},\gamma_{2,0}\subset H_2^{-1}(1)$ are nondegenerate, $\tilde H_N^{-1}(1)$ and $\tilde \Sigma_{E}$ admit respective pairs of periodic orbits $\gamma_1,\gamma_2$ and $\gamma_{1,E}, \gamma_{2,E}$ forming a Hopf link for every $\epsilon>0$ sufficiently small. Let $\rho_i$ and $\rho_{i,E}$ be the rotation number of $\gamma_i$ and $\gamma_{i,E}$, respectively. Since $|m_1|>m_2= 2$, the rotation number of $\gamma_2$ is given by \eqref{rho2 m2>2}, that is $\rho_2=1+\hat \omega_1(E)/\omega_2(E).$  Since $m_2=2$, the rotation number of $\gamma_1$ is given by \eqref{rot of gamma1 p=2 1}. Under conditions in (i) and (ii), we see from Theorem \ref{thm: non-resonant theorem_A} that
$$
(\rho_1-1)(\rho_2-1)=1+2^\nu\Omega_\nu E^{\nu-1}+O(E^{\nu}).
$$
Under the conditions in (iii), Theorem \ref{thm: non-resonant theorem_A} gives
$$
(\rho_{1}-1)(\rho_{2}-1)=1+\frac{|m_1|}{2}\cdot \left(\frac{2}{\alpha_2}\right)^\nu\Omega_{\nu,2}E^{\nu-1}+O(E^\nu).
$$
According to Lemma \ref{lem: perturbation}, $|\rho_{i,E}-\rho_i|\leq ME^{(N-1)/4}$ for some constant $M$ and every $E>0$ sufficiently small. Since $N$ can take arbitrarily large, we conclude in all cases (i), (ii), and (iii) that $(\rho_{1, E}-1)(\rho_{2, E}-1)\neq 1$ for every $E>0$ sufficiently small. This implies that the Hopf link $\gamma_{1,E}\cup \gamma_{2,E}\subset \tilde \Sigma_E$ is non-resonant.
\end{proof}

Theorem \ref{thm: main theorem} directly follows from Theorem \ref{thm: m2>2} and \ref{thm: m2=2}.

\subsection{Hamiltonians with symmetries}\label{sec: Hamiltonians with symmetries}
To deal with strongly resonant Hamiltonians at the origin, we need to consider Hamiltonians possessing some extra symmetries. Such symmetries allow us to control the families of periodic orbits near the origin that form a Hopf link at the corresponding energy surfaces.

We start treating the case $|m_1|>m_2=1$, i.e., $\alpha_2$ is a non-trivial multiple of $\alpha_1$. Under this condition, $\gamma_{1,0}\subset H_2^{-1}(1)$ degenerates while $\gamma_{2,0}$ is nondegenerate. Hence, we need conditions for the existence of  $\gamma_{1,E}\subset H^{-1}(E)$, for $E>0$ sufficiently small.

Fix $N>0$ and assume that $H=H_N+R$ is in normal form up to order $N$, where $R= \mathcal{O}_{N+1}$. We assume that both $\partial_{y_2}H$ and $\partial_{x_2}H$ vanish on the $y_1x_1$-plane ($y_2=x_2=0$) for $(y_1,x_1)$ sufficiently close to $(0,0)$. This assumption on the normal form is supported by Proposition \ref{lem: normal form 2} in the case the Hamiltonian $H = H_2 + \mathcal{O}_3$ satisfies the same conditions.

The first consequence of the symmetry conditions is that the origin $(0,0)$ of the $y_1x_1$-plane is surrounded by a family of periodic orbits in distinct energy surfaces. We denote these periodic orbits by $\gamma_{1, E}$, for every $E>0$ sufficiently small.

Let $\tilde H=\tilde H_N+\epsilon^{N-1}\tilde R$ be the rescaled Hamiltonian, where $\tilde H_N=\epsilon^{-2}H_N(\epsilon\cdot,\epsilon\cdot)$ and $\tilde R=\epsilon^{-(N+1)}R_{N+1}(\epsilon\cdot,\epsilon\cdot)$ are smooth $\epsilon$-dependent functions, and $\epsilon$ lies in a small neighborhood of $0$. The periodic orbits $\gamma_{1,E}\subset H^{-1}(E)$ rescale to periodic orbits $\gamma_{1,\epsilon}^\epsilon= \epsilon^{-1}\gamma_{1,E}\subset \tilde H^{-1}(1)$, where $E=\epsilon^2.$ Moreover, $\gamma_{1,\epsilon}^\epsilon \to \gamma_{1,0} \subset H_2^{-1}(1)$ in $C^\infty$ as $\epsilon \to 0^+$.

The notation $\gamma_\epsilon^\epsilon$ is justified by the following lemma, which is an adaptation of Lemma \ref{lem: perturbation} to the current situation.

\begin{lem}\label{lem: perturbation m2=1}
Assume that $|m_1|>m_2=1$ (i.e., $\alpha_2$ is a non-trivial multiple of $\alpha_1$) and that $H=H_N+R$ satisfies
$\partial_{y_2}H(y_1,0,x_1,0)=\partial_{x_2}H(y_1,0,x_1,0)=0,$
for every $(y_1,x_1)$ near $(0,0)$.
Consider the family $\tilde H_\delta:=\tilde H_N+\delta^{N-1} \tilde R$, defined for $(\epsilon,\delta)$ is close to $(0,0)$, where $\tilde H_N:=\epsilon^{-2}H_N(\epsilon \cdot,\epsilon \cdot)$ and $\tilde R:=\epsilon^{-(N+1)}R(\epsilon  \cdot,\epsilon \cdot)$. There exist $M>0$ and smooth families of periodic orbits $\gamma^\delta_{j,\epsilon}\in \tilde H_{\delta}^{-1}(1), j=1,2,$ with period $T^\delta_{j,\epsilon}$, defined for every $(\epsilon,\delta)$ sufficiently close to $(0,0)$, such that
\begin{equation}\label{estimates_gamma1}
{\rm dist}(\gamma^\delta_{j,\epsilon}, \gamma^0_{j,\epsilon})\leq M \delta^{N-1}
,\quad
|T^\delta_{j,\epsilon}-T^{0}_{j,\epsilon}|\leq M\delta^{N-1},\quad |\rho^\delta_{j,\epsilon}-\rho^{0}_{j,\epsilon}|\leq M \delta^{(N-1)/2},
\end{equation}
where $\rho^\delta_{j,\epsilon}$ is the rotation number of $\gamma^\delta_{j,\epsilon}$. Here, $\gamma_{1,\epsilon}^0, \gamma_{2,\epsilon}^0\subset \tilde H_N^{-1}(1)$ are periodic orbits of $\tilde H_N$, contained in $\{z_2=0\}$ and $\{z_1=0\},$ respectively, and forming a Hopf link. Moreover, $\gamma_{j,\epsilon}^\delta \to \gamma_{j,0}$ in $C^\infty$ as $(\epsilon,\delta) \to (0,0).$
\end{lem}

\begin{proof}
Let $\tilde \Sigma^\delta_{\epsilon}\subset \tilde H_\delta^{-1}(1)$ be the sphere-like component, which is $C^\infty$-close to $H_2^{-1}(1)$ for $(\epsilon,\delta)$ close to $(0,0)$. Let $z_j=x_j+iy_j,j=1,2$. From the assumptions on $H$, we obtain  $\partial_{z_2}\tilde H_N(z_1,0)=\partial_{z_2}\tilde R(z_1,0)=0$ for every $z_1$ near $0$. This implies that $\partial_{z_2}\tilde H_{\delta}(z_1,0)=0$ for every $(\epsilon,\delta)$ near $(0,0)$ and every $z_1$ near $0$. Hence, it follows from Hamilton's equations for $\tilde H_\delta$ that there exists a smooth family of periodic orbits $\gamma^\delta_{1,\epsilon}\subset \tilde \Sigma^\delta_\epsilon\cap \{z_2=0\}$ near $\gamma_{1,0}\subset H_2^{-1}(1)$ with period $T^\delta_{1,\epsilon}\approx T_{1,0}$,
solving the equation $\dot{\bar z}_1=2i\partial_{z_1}\tilde H_\delta(z_1,0)$, defined for $(\epsilon,\delta)$ sufficiently close to $(0,0)$. In particular, $\gamma^\delta_{1,\epsilon}$ and $T^\delta_{1,\epsilon}$ smoothly depend on $(\epsilon,\delta^{N-1})$. Indeed, $\gamma_{1,\epsilon}^\delta$ coincides with the regular circle-like level curve of  $\tilde H_\delta|_{z_2=0}$ at level $1$.  The regularity of such curves, i.e. $\nabla \tilde H_\delta|_{\gamma_{1,\epsilon}^\delta} \pitchfork \dot \gamma_{1,\epsilon}^\delta,$ implies that $\gamma_{1,\epsilon}^\delta$  smoothly depends on $(\epsilon,\delta^{N-1})$ for $(\epsilon,\delta)$ sufficiently close to $(0,0)$. By uniqueness, $\gamma_{1,\epsilon}^\delta$ coincides with the previously defined periodic orbit $\gamma_{1,\epsilon}^\epsilon \subset \tilde H^{-1}(1)$ for $\delta=\epsilon$ small, and  $\gamma_{1,\epsilon}^0 \subset \tilde H_N^{-1}(1).$ Moreover, $\gamma_{1,\epsilon}^\delta \to \gamma_{1,0}\subset H_2^{-1}(1)\cap \{y_2=x_2=0\}$ in $C^\infty$ as $(\epsilon,\delta) \to (0,0).$

Since $\gamma_{2,0}\subset H_2^{-1}(1)$ is nondegenerate, we obtain another smooth family of periodic orbits $\gamma_{2,\epsilon}^\delta\subset \tilde H_\delta^{-1}(1)$ smoothly depending on $(\epsilon,\delta^{N-1})$ near $(0,0)$. Given the smoothness of $\gamma_{j,\epsilon}^\delta$ and $T_{j,\epsilon}^\delta, j=1,2$, the estimates in \eqref{estimates_gamma1} are proved as in Lemma \ref{lem: perturbation}.
\end{proof}

Under the assumptions of Lemma \ref{lem: perturbation m2=1}, the following theorem gives conditions on the Birkhoff-Gustavson normal form that imply the existence of a non-resonant Hopf link for $H$.

\begin{thm}\label{thm: m2=1, m1>1}
Assume that  $|m_1|> m_2=1$  and that $H=H_N+\mathcal{O}_{N+1}$ satisfies $\partial_{y_2}H(y_1,0,x_1,0)=\partial_{x_2}H(y_1,0,x_1,0)=0
$ for every $(y_1,x_1)$ near $(0,0)$. Assume also that one of the following conditions is satisfied:
\begin{itemize}
\item[(i)] $|m_1|> 2$, $|m_1|> \nu -1$ and
$\Omega_{\nu,1},\Omega_\nu\neq 0$.
\item[(ii)]  $|m_1|> 2$, $|m_1|=\nu -1$, $\Omega_{\nu,1},\Omega_{\nu}\neq 0$ and $a_{0,2,2|m_1|,0}=0$.
\item[(iii)]  $|m_1|> 2$, $|m_1|<\nu -1$, $\Omega_{\nu,2}\neq 0$ and $a_{0,2,2|m_1|,0}\neq 0$.
\item[(iv)] $|m_1|=2,$ $\nu=2$, $\Omega_{2,1}\neq 0$ and $a_{0,1,2,0} \neq 0.$
\item[(v)]  $|m_1|= 2$, $\nu=2$, $\Omega_{2,1},\Omega_2\neq 0$.
\item[(vi)]  $|m_1|= 2$, $\nu=3$, $\Omega_{3,1},a_{0,1,2,0}\neq 0$ and $a_{0,2,4,0}=0$.
\end{itemize}
Here,  $2\leq \nu\leq \lfloor N/2\rfloor$ is the smallest integer so that one of the numbers in \eqref{condition of nv} does not vanish. Then the sphere-like component $\Sigma_E\subset H^{-1}(E)$ carries a non-resonant Hopf-link for every $E>0$ sufficiently small.
\end{thm}

\begin{proof}
We know from Lemma \ref{lem: perturbation m2=1} that the energy surfaces $H_N^{-1}(E)$ and $H^{-1}(E)$ admit respective Hopf links $\gamma_1\cup \gamma_2$ and $\gamma_{1,E}\cup \gamma_{2,E}$, for every $E>0$ sufficiently small. Let $\rho_i$ and $\rho_{i,E}$ denote the rotation number of $\gamma_i$ and $\gamma_{i,E},i=1,2$. Under the conditions in (i)-(ii), we obtain from Theorem~\ref{thm: non-resonant theorem_A} that
$(\rho_1-1)(\rho_2-1)=1+2^\nu\Omega_\nu E^{\nu-1}+O(E^\nu).$
Under the conditions in (iii), we obtain from Theorem \ref{thm: non-resonant theorem_A} that $(\rho_1-1)(\rho_2-1)=1+|m_1|(\frac{2}{\alpha_2})^\nu\Omega_{\nu,2} E^{\nu-1}+O(E^\nu).$
In cases (iv), (v), and (vi), Theorem \ref{thm: non-resonant theorem B} gives either
$(\rho_1-1)(\rho_2-1)=1+\frac{1}{2}(\frac{2}{\alpha_1})^\nu\Omega_{\nu,1} E^{\nu-1}+O(E^\nu)$ or $(\rho_1-1)(\rho_2-1)=1+2^\nu\Omega_{\nu} E^{\nu-1}+O(E^\nu), \nu=2,3.$

By Lemma \ref{lem: perturbation m2=1}, $|\rho_{i,E}-\rho_i|<ME^{(N-1)/4}$ for some constant $M$. Since $N$ can be taken arbitrarily large, we conclude in all cases that $(\rho_{1,E}-1)(\rho_{2,E}-1)\neq 1$ for every $E>0$ sufficiently small. This implies that the Hopf link $\gamma_{1,E}\cup \gamma_{2,E}\subset H^{-1}(E)$ is non-resonant.
\end{proof}

 Now, we consider the most resonant case $|m_1|=m_2=1,$ that is $\alpha_1=\alpha_2$, with the same type of symmetry.
 The following theorem is a straightforward generalization of Theorem \ref{thm: m2=1, m1>1}.

\begin{thm}\label{thm: m2=m1=1}
Assume that $|m_1|=m_2=1$ and that the Birkhoff-Gustavson normal form $H=H_N  + \mathcal{O}_{N+1}$ satisfies
$$
\partial_{y_1}H(0,y_2,0,x_2)=\partial_{x_1}(0,y_2,0,x_2)=\partial_{y_2}H(y_1,0,x_1,0)=\partial_{x_2}(y_1,0,x_1,0)=0,
$$
for every $(y_2,x_2)$ and $(y_1,x_1)$ sufficiently close to $(0,0)$. Assume also that  one of the following conditions is satisfied:
\begin{itemize}
\item[(i)] $\nu=2$, $\Omega_{2}\neq 0$ and $a_{0,2,2,0}=0$;
\item[(ii)] $\nu=2$, $\Omega_{2,1}=-\Omega_{2,2}\neq 0, a_{0,2,2,0}=0$ and $\beta_1+\beta_2+2\Omega_{2,1}\Omega_{2,2}\neq 0$,
\end{itemize}
where $2\leq \nu\leq \lfloor N/2\rfloor$ is the smallest integer so that one of the numbers in \eqref{condition of nv} does not vanish and $\beta_i$ is given by \eqref{beta1,2}. Then the sphere-like component $\Sigma_E\subset H^{-1}(E)$ carries a non-resonant Hopf-link for every $E>0$ sufficiently small.
\end{thm}

\begin{proof}
We know from Lemma \ref{lem: perturbation m2=1} that the energy surfaces $H_N^{-1}(E)$ and $H^{-1}(E)$ admit Hopf links $\gamma_1\cup \gamma_2$ and $\gamma_{1,E}\cup \gamma_{2,E}$, respectively. Let $\rho_j$ and $\rho_{j,E}$ denote the rotation number of $\gamma_j$ and $\gamma_{j,E},j=1,2$.
Lemma \ref{lem: perturbation m2=1} also gives $|\rho_{j,E}-\rho_j|\leq ME^{(N-1)/4}, j=1,2$, for some $M>0$.
By Theorem \ref{thm: non-resonant theorem C}, we have
$(\rho_1-1)(\rho_2-1)=1+4\Omega_2 E+8(\beta_1+\beta_2+2\Omega_{2,1}\Omega_{2,2})E^2+O(E^3)$, where $\Omega_2=\Omega_{2,1}+\Omega_{2,2}$.
Since $N$ can be taken arbitrarily large, we conclude that $(\rho_{1,E}-1)(\rho_{2,E}-1)\neq 1$ for every $E>0$ sufficiently small. This implies that the Hopf link $\gamma_{1,E}\cup \gamma_{2,E}$ is non-resonant.
\end{proof}

Theorem \ref{thm: main theorem 2} directly follows from Theorem \ref{thm: m2=1, m1>1}.

\subsection{\texorpdfstring{$\Z_p$}{Lg}-symmetric Hamiltonians} \label{sec: zp symmetry}
Let $p\geq 3$ be an integer and let $\hat R$ be the rotation by $2\pi/p$ on $\R^2$, that is $\hat R y = e^{2\pi i /p} y$, for every $y\in \C \equiv \R^2$. Let
$$
R(y,x) := (\hat Ry,\hat R x)=(e^{2\pi i /p}y, e^{2\pi i /p} x).
$$
be the unitary map on $\R^4$ induced by $\hat R$ on the Lagrangian planes $(y_1,y_2)$ and $(x_1,x_2)$.

We consider $\Z_p$-invariant Hamiltonians $H\circ R = H =H_2 + \mathcal{O}_3$ satisfying  $m_2=|m_1|=1$. Then $R$ induces a $\Z_p$-action on  the sphere-like components $\Sigma_E\subset H^{-1}(E)$ near the origin, and $\Sigma_E / \Z_p$ is diffeomorphic to the lens space $L(p,p-1)$. The contact form $\lambda_0 = \frac{1}{2} \sum_{j=1}^2 y_jdx_j - x_j dy_j|_{\Sigma_E}$ is $R$-invariant and  descends to the standard contact structure on $L(p,p-1)$. Indeed, up to a symplectic change of coordinates, $R$ is equivalent to the map
\begin{equation}\label{equ: mathcal R}
\mathcal{R}:(x_1+iy_1,x_2+iy_2) \mapsto (e^{-2\pi i /p}(x_1+iy_1), e^{2\pi i /p}(x_2+iy_2)),
\end{equation}
which rotates the symplectic planes $(x_1,y_1)$ and $(x_2,y_2)$ by $2\pi/p$ in opposite directions.

In the following we proceed by induction and assume that $H=H_{N-1} + \hat H_{N-1}^{(N)} + \mathcal{O}_{N+1}$ is in normal form up to order $N-1$ and $\hat H_{N-1}^N=K + L$ is a $N$-homogeneous polynomial in coordinates $(y,x)$, with $L\in {\rm im} D$ and $K \in \ker D$.  We assume that $H$ and $H_{N-1}$ are $R$-invariant. In particular, $\hat H_{N-1}^{(N)}$ is $R$-invariant.  We show that the Birkhoff-Gustavson normal form $H=H_N+\mathcal{O}_{N+1}$ and the Hamiltonian $H$ in the new coordinates $(\eta,\xi)$ are also $R$-invariant. To do so, we consider a symplectic change of coordinates $\Phi: (\eta_1,\eta_2,\xi_1,\xi_2) \mapsto (y_1,y_2,x_1,x_2)$ that is $R$-invariant, and  the new Hamiltonian $H \circ \Phi_N = H_N + \mathcal{O}_{N+1}$ is also $R$-invariant. To achieve these conditions, it is crucial to assume that $G_N$ and $K_N$ are $R$-invariant, i.e.,
\begin{equation}\label{GK}
G_N \circ R = G_N \quad \mbox{ and } \quad K_N \circ R = K_N
\end{equation}
where $G_N(\eta,x)$ is the $N$-homogeneous polynomial satisfying $D\cdot G_N = -L\in {\rm im} D$. Recall that $G_N$ determines the generating function $W_N(\eta,x) = \eta x +G_N(\eta,x)$, which defines $\Phi$ by
$$
\left\{\begin{aligned}
\xi & = x + \partial_\eta G(\eta,x),\\
  y & = \eta + \partial_x G(\eta,x).
\end{aligned}\right.
$$
Applying $\hat R$ to each equation above, we obtain
$$
\left\{\begin{aligned}
\hat R\xi & = \hat R x + \partial_\eta G(\hat R \eta,\hat R x),\\
\hat R  y & = \hat R \eta + \partial_x G (\hat R \eta,\hat R x).
\end{aligned}\right.
$$
We have used that $R^T(\nabla G \circ R)= \nabla G,$ which follows from \eqref{GK}, and that $R^T  = R^{-1} = (\hat R^{-1},\hat R^{-1}).$ We conclude that
$$
\Phi \circ R(\eta,\xi)=\Phi(\hat R \eta ,\hat R \xi) = (\hat R y,\hat R x) = R(y,x) = R \circ \Phi(\eta,\xi),$$
i.e., $\Phi$ is $R$-invariant. This also implies that the new Hamiltonian $H \circ \Phi$ is $R$-invariant
$$
H \circ \Phi \circ R = H \circ R \circ \Phi = H \circ \Phi.
$$
Now the order-$N$ polynomial $H_N$ is given by $H_N = H_{N-1} + K_N.$ Since $H_{N-1}$ and $K_N$ are $R$-invariant, we conclude that $H_N$ is also $R$-invariant.

In the following proposition we show that assumptions \eqref{GK} hold true. This means that under the change of coordinates $(\eta,\xi) \mapsto (y,x)$, $H_N$ can always be taken $R$-invariant.

\begin{prop}[$\Z_p$-symmetry]\label{prop: zp symmetry}
Assume that $H=H_2 + \mathcal{O}_3$ is $R$-invariant in coordinates $(y,x)$ and that $\alpha_1=\alpha_2=1$. Given $N>0$, there exist an $R$-invariant symplectic change of coordinates $\Phi:(\eta,\xi) \mapsto (y,x)$ so that the Birkhoff-Gustavson normal form  $H\circ \Phi=H_N+\mathcal{O}_{N+1}$ is $R$-invariant in coordinates $(\eta,\xi)$.
\end{prop}

\begin{proof}
In complex coordinates $z_j=x_j+iy_j,j=1,2,$ we associate to the monomial $Q=z_1^{k_1}z_2^{k_2}\bar z_1^{l_1}\bar z_2^{l_2}$ the $R$-invariant polynomial
$$
Q_R:=\sum_{j=0}^{p-1}Q_j=\sum_{j=0}^{p-1}z_{1,j}^{k_1}z_{2,j}^{k_2}\bar z_{1,j}^{l_1}\bar z_{2,j}^{l_2},\quad Q_j=Q\circ R^j.
$$
where $z_{1,j}=\cos(2\pi j/p) z_1-\sin(2\pi j/p) z_2$ and $z_{2,j}=\sin(2\pi j/p) z_1 + \cos(2\pi j/p) z_2$ for every $j$. Since $H$ is $R$-invariant, the $s$-homogeneous part $H^{(s)}$ of $H$ is $R$-invariant for every $s$. If $p$ is even, then $Q_{j}=(-1)^{k_1+k_2+l_1+l_2}Q_{j+p/2}$, and thus $H^{(s)}$ vanishes for every odd degree $s$. Moreover, the relation $H\circ R=H$ implies that $\nabla H\circ R=(R^T)^{-1}\circ \nabla H=R\circ \nabla H$ since $R^T = R^{-1}$.

Consider the generating function $W_3(\eta,x)=x_1\eta_1+x_2\eta_2+G_3(\eta,x)$ and the induced canonical transformation $(y,x)=\Phi_3(\eta,\xi)$, where $G_3$ is a $3$-homogeneous polynomial in coordinates $(\eta,x)$, which formally replace coordinates $(y,x)$. Recall from \eqref{D on monomials} that $D\cdot Q=i^{-1}(k_1-l_1+k_2-l_2)Q$ for every complex monomial $Q$. Recall also from section \ref{sec: construction} that $G_3$ satisfies the relation $-D\cdot G_3+\Gamma^{(3)}=H^{(3)}$, where $H^{(3)}$ is the $3$-homogeneous part of $H$ and $\Gamma^{(3)}\in \ker D$. Notice that $G_3$ is uniquely determined up to elements in $\ker D$. We claim that $\Gamma^{(3)}$ is $R$-invariant. In fact, if $Q\in\ker D$ is a $3$-homogeneous polynomial, then  $(D\cdot Q)(w)=\{H_2,Q\}(w)=(J\nabla H_2(w),\nabla Q(w))=0,\forall w=(y_1,y_2,x_1,x_2)\in \mathbb{R}^4$. Since $R^T\nabla Q\circ R=\nabla (Q\circ R)$ and $R^T\nabla H_2\circ R=\nabla H_2$ due to $H_2\circ R=H_2$, we compute
\begin{equation}\label{equ: D and R}
\begin{aligned}
(D\cdot (Q\circ R))(w)&=(J\nabla H_2(w), \nabla (Q\circ R)(w))=(J R^T\nabla H_2(Rw), R^T\nabla Q(Rw))\\
&=(J \nabla H_2(Rw), \nabla Q(Rw))=(D\cdot Q)(Rw)=0,\quad \forall w\in\mathbb{R}^4.
\end{aligned}
\end{equation}
We have used that $RJR^T=J$ since $R$ is symplectic. We conclude that $Q_j:=Q\circ R^j$ also lies in $\ker D$ for every $j$. We also see from \eqref{equ: D and R} that if $Q\not \in \ker D$, then $Q \circ R\notin \ker D.$ Since $H^{(3)}$ is $R$-invariant and $\Gamma^{(3)}$ is the component of $H^{(3)}$ in $\ker D$, $\Gamma^{(3)}$ is also  $R$-invariant. Now we assume $Q$ is a complex monomial in $H^{(3)}-\Gamma^{(3)}$. Since $H^{(3)}-\Gamma^{(3)}$ is $R$-invariant, the associated $R$-invariant polynomial $Q_R$ generated by $Q$ is also in $H^{(3)}-\Gamma^{(3)}$. According to the computation above, we have
$$
\begin{aligned}
D\cdot Q_R&=\sum_{j=0}^{p-1}D\cdot(Q\circ R^j)=\sum_{j=0}^{p-1}(D\cdot Q)\circ R^j=\sum_{j=0}^{p-1}i^{-1}(k_1-l_1+k_2-l_2)Q\circ R^j\\
&=i^{-1}(k_1-l_1+k_2-l_2)Q_R.
\end{aligned}
$$
It follows that $\frac{i}{k_1-l_1+k_2-l_2}Q_R$ is part of $G_3$. Since $H^{(3)}-\Gamma^{(3)}$ is a finite sum of such $R$-invariant polynomials $Q_{R,j},j=1\ldots k$, $G_3$ is given by $\sum_{j=1}^k \frac{i}{k_{1,j}-l_{1,j}+k_{2,j}-l_{2,j}}Q_{R,j}$, which is also $R$-invariant.  From the previous discussion, we have $\nabla G_3\circ R=R\circ \nabla G_3$. Now we show that $\Phi_3\circ R=R\circ \Phi_3$. Notice that the relation $\Phi_3(\eta,\xi)=(y,x)$ is equivalent to the identity
\begin{equation}\label{equ: relation of Phi_3}
(\xi,y)=(x,\eta)+(\partial_{\eta}G_3,\partial_{x} G_3).
\end{equation}
Let $R(y,x):=(\hat R y,\hat R x)$. Applying $R$ to both sides of \eqref{equ: relation of Phi_3}, we obtain
$$
(\hat R \xi,\hat R y)=(\hat R x,\hat R \eta)+R\nabla G_3(\eta,x)=(\hat R x,\hat R \eta)+\nabla G_3(\hat R \eta,\hat R x).
$$
This implies that $\Phi_3$ also sends $(\hat R \eta,\hat R \xi)$ to $(\hat R y,\hat R x)$. Therefore, $\Phi_3\circ R(\eta,\xi)=R(y,x)=R\circ \Phi_3(\eta,\xi)$. This implies that $\hat H_3:=H\circ \Phi_3$ is $R$-invariant.

Now, by induction, we assume $G_{s-1}\circ R=G_{s-1}$ and $\hat H_{s-1}\circ R=\hat H_{s-1}$. Consider the generating function $W_s=x_1\eta_1+x_2\eta_2+G_s$ and the induced canonical transformation $\Phi_s$, where $G_s$ is a $s$-homogeneous polynomial determined by $-D\cdot G_{s}+\Gamma^{(s)}=\hat H^{(s)}_{s-1}$ and $\hat H^{(s)}_{s-1}$ is the $s$-homogeneous part of $\hat H_{s-1}$, see section \ref{sec: construction}. From the previous discussion, we know that $\hat H^{(s)}_{s-1}$ is $R$-invariant. Again, since $\Gamma^{(s)}$ is formed by the monomials of $\hat H^{(s)}_{s-1}$ in $\ker D$ and $G_s$ is uniquely determined by the identity $D \cdot G_s = \Gamma^{(s)} - \hat H^{(s)}_{s-1}$ up to elements in $\ker D$,  both $\Gamma^{(s)}$ and $G_s$ are $R$-invariant. Therefore, we obtain $\Phi_{s}\circ R=R\circ \Phi_{s}$ as before and $\hat H_s:=\hat H_{s-1}\circ \Phi_s=\Gamma^{(2)}+\cdots+\Gamma^{(s)}+R_{s+1}$ is $R$-invariant. This also implies that the normal form $H_s=\Gamma^{(2)}+\cdots+\Gamma^{(s)}$ and the remainder term $R_{s+1}$ are $R$-invariant.  \end{proof}

In the following, we consider the Birkhoff-Gustavson normal form $H = H_N + R_{N+1}$, where $\alpha_1=\alpha_2=1$ and $R_{N+1} = \mathcal{O}_{N+1}$. As before, we define $\tilde H_N (y,x):= \epsilon^{-2}H_N(\epsilon y, \epsilon x)$ and $\tilde R(y,x) = \epsilon^{-(N+1)}R_{N+1}(\epsilon y, \epsilon x),$ for every $\epsilon>0$ small.

\begin{prop}\label{prop: smooth family}
Assume that $\alpha_1=\alpha_2=1$ and let $H=H_N + \mathcal{O}_{N+1}$ be in normal form up to order $N$ as in Proposition \ref{prop: zp symmetry}. Moreover, $H$ is invariant under the $\mathbb{Z}_p$-action generated by the unitary map $R:(y,x)\mapsto (e^{2\pi i/p}y,e^{2\pi i/p}x)$, $p\geq 3$.  Let $\tilde H_\delta:=\tilde H_N+\delta^{N-1}\tilde R$, where $(\epsilon,\delta)\in \R^2$ is close to $(0,0)$. Then
\begin{itemize}

\item[(i)] There exist smooth families of periodic orbits $\gamma^\delta_{j,\epsilon}\subset \tilde H_\delta^{-1}(1), j=1,2$, forming a Hopf link, where $(\epsilon,\delta)$ is sufficiently close to $(0,0)$. These periodic orbits are $R$-invariant and $p$-cover nondegenerate periodic orbits $\hat \gamma_{1,\epsilon}^\delta$ and $\hat \gamma_{2,\epsilon}^\delta$ on $\tilde H_\delta^{-1}(1)/\Z_p \equiv L(p,p-1).$

\item[(ii)] Let $\Psi(v_1,v_2,u_1,u_2):=2^{-1/2}(v_1+v_2,u_1-u_2,u_1+u_2,v_2-v_1)$. The Hamiltonian $K(v,u):=H \circ \Phi \circ \Psi = H_N \circ \Psi + \mathcal{O}_{N+1}$ is in normal form up to order $N$. It is $\mathcal{R}$-invariant, where $\mathcal{R}(z_1,z_2)=(e^{-2\pi i/p}z_1,e^{2\pi i /p} z_2)$ for every $  (z_1,z_2)=(u_1+iv_1,u_2+iv_2).$ The monomials $Q=cz_1^{k_1}z_2^{k_2}\bar z_1^{l_1}\bar z_2^{l_2}$ of $K_N:=H_N\circ \Psi$ satisfy $|l_1-k_1|=|k_2-l_2|\neq 1$. Moreover, $\gamma^0_{1,\epsilon} \subset \{|z_1|= c_1,z_2=0\}$ and $\gamma^0_{2,\epsilon} \subset \{z_1=0,|z_2| = c_2\}$, where $c_j=c_j(\epsilon)>0$ for every $\epsilon>0$ sufficiently small.

\item[(iii)] The rotation numbers of $\gamma_{j,\epsilon}^\delta$ satisfy $|\rho^\delta_{j,\epsilon}-\rho^0_{j,\epsilon}|<M\delta^{(N-1)/2}$ for every $(\epsilon,\delta)$ sufficiently close to $(0,0)$.
\end{itemize}
\end{prop}

\begin{proof}
Since $H$ and $H_N$ are $R$-invariant, see Lemma \ref{prop: zp symmetry}, the family $\tilde H_\delta$ is also $R$-invariant for every $(\delta,\epsilon)$ near $(0,0)$. Notice that $\mathcal{R}:=\Psi^{-1}\circ R\circ \Psi$ is given by $\mathcal{R}(z_1,z_2)= (e^{-2\pi i/p}z_1,e^{2\pi i/p}z_2)$, where $z_j=u_j+iv_j,j=1,2$. Also,  $H_2\circ\Psi=H_2$ and the sphere-like component $\Sigma_\epsilon^\delta\subset (\tilde H_\delta\circ \Psi)^{-1}(1)$ is $\mathcal{R}$-invariant.
Taking the quotient of $\Sigma^\delta_{\epsilon}$ under the $\mathbb{Z}_p$-action generated by $\mathcal{R}$, we obtain a smooth family of three-manifolds $\hat \Sigma^\delta_\epsilon$ diffeomorphic to the lens space $L(p,p-1)=S^3/\mathbb{Z}_p$, and converging in $C^\infty$ to $S^3/\Z_p$ as $(\epsilon,\delta)\to (0,0)$. Moreover, since the Liouville form $\lambda_0=\frac{1}{2}\sum_{i=1}^2(v_idu_i-u_idv_i)=\frac{1}{4i}\sum_{i=1}^2(z_id\bar z_i-\bar z_i dz_i)$ is  $\mathcal{R}$-invariant, $\hat \Sigma^\delta_\epsilon$ is equipped with the standard contact form $\hat \lambda^\delta_\epsilon$, whose contact structure is contactomorphic to the standard contact structure $\xi_0=\ker \lambda_0$ on $L(p,p-1)=S^3 / \Z_p$. In summary, $(\hat \Sigma^\delta_\epsilon,  \xi^\delta_\epsilon:=\ker \hat \lambda^\delta_\epsilon)$ is a family of contact lens spaces all contactomorphic to $(L(p,p-1), \xi_0)$. Recall that the contact structure $\xi_0 \to L(p,p-1)$ is trivial. We fix one such trivialization induced by a $\mathcal{R}$-invariant trivialization of $(S^3,\xi_0)$.

In coordinates $(y_1,y_2,x_1,x_2)$, consider the $\Z_p$-invariant periodic orbits of $H_2^{-1}(1)$ given by
$$
\gamma_{1,0}(t) = (\cos t ,\sin t,\sin t ,-\cos t) \quad \mbox{ and } \quad \gamma_{2,0}(t) = (\cos t,-\sin t,\sin t,\cos t).
$$
Their rotation numbers with respect to the global trivialization of $\xi_0$ are both equal to $2$. They $p$-cover periodic orbits $\hat \gamma_{1,0}$ and $\hat \gamma_{2,0}$ on $H_2^{-1}(1)/\Z_p$, respectively, whose rotation numbers in the global frame, are $\hat \rho_{j,0} = \frac{2}{p}, j=1,2$. Since $p\geq 3$,  $\hat \gamma_{1,0}$ and $\hat \gamma_{2,0}$ are nondegenerate.
In coordinates $(v_1,v_2,u_1,u_2)$, these periodic orbits satisfy
$$
\gamma_{1,0} \subset \left\{z_2=0, |z_1|=2^{1/2}\right\} \quad \mbox{ and } \quad \gamma_{2,0} \subset \left\{z_1=0, |z_2|=2^{1/2}\right\}.
$$

As $(\epsilon,\delta) \to (0,0)$, $\Sigma^\delta_\epsilon$ converges in $C^\infty$ to the standard three-sphere $\Sigma_0:=H_2^{-1}(1)$. Since $\hat \gamma_{1,0}$ and $\hat \gamma_{2,0}$ are nondegenerate, we obtain smooth $\Z_p$-invariant families of periodic orbits $\hat \gamma_{1,\epsilon}^\delta,\hat \gamma_{2,\epsilon}^\delta$  on $\hat \Sigma^\delta_\epsilon=\Sigma^\delta_\epsilon / \Z_p$, which are continuations of $\hat \gamma_{1,0}$ and $\hat \gamma_{2,0}$, respectively. If $\delta=0$, then $\tilde H_\delta=\tilde H_N$, and  $\gamma_i=\epsilon \gamma_{j,\epsilon}^0,j=1,2$ are periodic orbits on $H_N^{-1}(\epsilon^2)$.

Denote by $K_N=H_N\circ \Psi$ the $\mathcal{R}$-invariant Hamiltonian in coordinates $(v,u)$.  It is straightforward to check that in coordinates $(v,u)$ the operator $D$ has the same form as in coordinates $(y,x)$, i.e.,  $D= \sum_{j=1}^2 v_j\partial_{u_j} - u_j \partial_{v_j}=\sum_{j=1}^2 y_j\partial_{x_j}-x_j \partial_{y_j}$.  In particular, $D \cdot K_N=0$ and thus $H \circ \Psi = K_N + \mathcal{O}_{N+1}$ is in normal form up to order $N$. Due to the previous discussion, the $j$-homogeneous part $K^{(j)}_N$ of $K_N$ is $\mathcal{R}$-invariant and $D \cdot K^{(j)}_N=0$ for every $j=1,\ldots N$. Let $Q=z_1^{k_1}z_2^{k_2}\bar z_1^{l_1}\bar z_2^{l_2}$ be a monomial of $K^{(j)}_N$. From \eqref{D on monomials}, we have $D\cdot Q =0$, and thus $(1,1)\cdot (k-l) =0\Rightarrow k_1-l_1+k_2-l_2=0$. Moreover, $Q\circ \mathcal{R}=e^{2\pi(l_1-k_1+k_2-l_2)i/p}Q$. Since $K_N^{(j)}$ is $\mathcal{R}$-invariant, we conclude that $Q\circ \mathcal{R} = Q$ and $l_1-k_1+k_2-l_2=mp$ for some $m\in \Z$. It follows that $k_2-l_2=l_1-k_1=\frac{mp}{2}$ for some $m\in \Z$. Since $p\geq 3,$ we have
\begin{equation}\label{condkl}
|k_2-l_2|=|l_1-k_1|=|\frac{mp}{2}|>1, \quad \forall m\in \Z\setminus \{0\}.
\end{equation}
This shows that $K_N$ has no monomials of the form $z_1 z_2^{k_2}\bar z_2^{l_2}$ or $\bar z_1z_2^{k_2}\bar z_2^{l_2}.$ In particular, $\partial_{z_1} K_N|_{z_1=0}=0$. Similarly, $\partial_{z_2} K_N|_{z_2=0}=0$. These identities imply that $K_N$ admits periodic orbits $\gamma_1 \subset \{|z_1|=c_1>0, z_2=0\}$ and $\gamma_2 \subset \{z_1=0, |z_2|=c_2>0\}$ forming a Hopf link for every positive energy sufficiently small. This follows from  Hamilton's equations \eqref{H system 1}. Indeed, conditions \eqref{cond: gamma2 m2=1} and \eqref{cond: gamma1 m2=1} follow from \eqref{condkl}. By uniqueness of continuations,  $\gamma_1$ and $\gamma_2$ correspond to $\gamma_{1,\epsilon}^0$ and $\gamma_{2,\epsilon}^0$ after re-scaling.

Finally, since $\hat \gamma_{1,\epsilon}^\delta\subset \hat \Sigma_\epsilon^\delta\equiv L(p,p-1)$ is nondegenerate, we can argue as in the proof of Lemma \ref{lem: perturbation} to conclude that there exists a uniform $M>0$ such that $|\rho^\delta_{i,\epsilon}-\rho^0_{i,\epsilon}|\leq M\delta^{(N-1)/2}$ for every $(\epsilon,\delta)$ sufficiently close to $(0,0)$. Moreover, if $p\neq 4$, then the Floquet multipliers of $\hat \gamma_{j,\epsilon}^\delta, j=1,2,$ do not coincide with $\pm 1$ and the previous estimate improves to $|\rho^\delta_{i,\epsilon}-\rho^0_{i,\epsilon}|\leq M\delta^{N-1}$ for every $(\epsilon,\delta)$ sufficiently close to $(0,0)$.
\end{proof}

\begin{thm}\label{thm: m2=m1=1 3}
Under the conditions of Proposition \ref{prop: smooth family}, assume that $K = K_N + \mathcal{O}_{N+1},$ and one of the following conditions is satisfied:
\begin{itemize}
\item[(i)] $\nu=2$, $\Omega_{2}\neq 0$, $a_{0,2,2,0}=0$;
\item[(ii)] $N\geq 6$, $\nu=2$, $\Omega_{2,1}=-\Omega_{2,2}\neq 0, a_{0,2,2,0}=0$, $\beta_1+\beta_2+2\Omega_{2,1}\Omega_{2,2}\neq 0$,
\end{itemize}
where $2\leq \nu\leq \lfloor N/2\rfloor$ is the smallest integer so that one of the numbers in \eqref{condition of nv} does not vanish and $\beta_i$ is given by \eqref{beta1,2}. Then the sphere-like component of $K^{-1}(E)$ carries a non-resonant Hopf link $\gamma_{1,E}\cup\gamma_{2,E}$ for every $E>0$ sufficiently small.
\end{thm}

\begin{proof}
Recall from Proposition \ref{prop: smooth family}-(ii) that $K_N$ admits two periodic orbits $\gamma^0_{i,\epsilon}$ on $(v_i,u_i)$-plane, $i=1,2$. Following the notation in Proposition \ref{prop: smooth family}, we first estimate the rotation numbers $\rho_{1,\epsilon}^0,\rho_{2,\epsilon}^0$ of $\gamma_{1,\epsilon}^0, \gamma_{2,\epsilon}^0\subset K_N^{-1}(\epsilon^2),$ respectively.
From the proof of Theorem \ref{thm: non-resonant theorem C}, we have $(\rho_{1,\epsilon}^0-1)(\rho_{2,\epsilon}^0-1)=1+4\Omega_2 E+8(\beta_1+\beta_2+2\Omega_{2,1}\Omega_{2,2})E^2+O(E^3)$, where $\Omega_2=\Omega_{2,1}+\Omega_{2,2}$ and $E=\epsilon^2$. By Proposition \ref{prop: smooth family}-(iii), we further obtain $|\rho_{i,\epsilon}^0-\rho^\delta_{i,\epsilon}|<M\delta^{(N-1)/2}$. We may assume that $N$ is arbitrarily large. Taking $\delta=\epsilon>0$ small, we see in both cases (i) and (ii) that $(\rho_{1,\epsilon}^\epsilon-1)(\rho_{2,\epsilon}^\epsilon-1)\neq 1$ for every $\epsilon>0$ sufficiently small. Hence, the Hopf link $\gamma_{1,\epsilon}^\epsilon\cup \gamma_{2,\epsilon}^\epsilon \subset K^{-1}(\epsilon^2)$ is non-resonant for every $\epsilon>0$ sufficiently small.
\end{proof}

Theorem \ref{thm: main theorem 3} directly follows from Proposition \ref{prop: smooth family} and Theorem \ref{thm: m2=m1=1 3}.

\begin{rem}\label{rem: normal form in two coordinates}
Let the Hamiltonians $H=H_2+\mathcal{O}_2$ and $K=H\circ \Psi=H_2+\mathcal{O}_2$ be as above. As explained in section \ref{sec: construction}, there exist symplectic changes of coordinates $\Phi_{N,1}$ and $\Phi_{N,2}$ such that $H\circ \Phi_{N,1}=H_N+\mathcal{O}_{N+1}$ and $K\circ \Phi_{N,2}=K_N+\mathcal{O}_{N+1}$ for every $N\in\mathbb{N}$. We remark that  $H_N\circ \Psi$ may not coincide with $K_N$ since $\Phi_{N,1}\circ \Psi$ and $\Psi\circ \Phi_{N,2}$ may not coincide. This is illustrated in Hill's Lunar Problem, see section \ref{Hill}, where $\Phi_{N,1}\circ \Psi$ and $\Psi\circ \Phi_{N,2}$ do not coincide for $N=4$. In this case, $H\circ \Phi_{N,1}\circ \Psi-K\circ \Phi_{N,2}=\mathcal{O}_{6}$. In particular, if $H$ is $R$-invariant and $K$ is $\mathcal{R}$-invariant, then $H_N$ and $K_N$ may not keep the corresponding symmetry at the same time.
\end{rem}

\section{The spatial isosceles three-body problem}\label{SI3BP}

The spatial isosceles three-body problem is the study of the motion of three-point bodies in $\mathbb{R}^3$ subject to Newton's gravitational law, in which two of them have equal masses and move symmetrically with respect to a fixed axis where the third body moves. Following \cite{HLOSY, Mk84}, this problem can be reduced to the dynamics of the mechanical Hamiltonian with two degrees of freedom
$$
H(p_r,p_z,r,z)=\frac{p_r^2+p_z^2}{2}+\frac{\varpi^2}{2r^2}-\frac{1}{r}-\frac{4\alpha^{-1}}{\sqrt{r^2+(1+2\alpha)z^2}},
$$
where $r>0,z\in \R$ are suitable cylindrical coordinates on the plane determined by the bodies, $p_r,p_z$ are the corresponding momenta, and $\alpha,\varpi>0$ are parameters, called mass ratio and angular momentum, respectively.

We consider the dynamics on the energy surface $\mathfrak{M}:= H^{-1}(h)$ for negative energies $h$. Indeed, there exist only two essential parameters: $\alpha$ and $\varpi^2h$. Without loss of generality, we may fix the energy $h=-1$ and replace $\varpi^2h$ with the new parameter $\mathfrak{e}$, called  eccentricity, given by
$$
\mathfrak{e}:=\sqrt{1-\frac{2\varpi^2}{(1+4/\alpha)^2}}\in(0,1),
$$
We mainly study the dynamics of $\mathfrak{M}$ for $\mathfrak{e}$ sufficiently small. In this case, $\mathfrak{M}$ is a sphere-like energy surface close to the  minimum $(0,0,\frac{\varpi^2}{1+4/\alpha},0)$ corresponding to $\mathfrak{e}=0$. If $\mathfrak{M}$ has a non-resonant Hopf link, then the results in \cite{HLOSY} imply the existence of infinitely many periodic orbits in $\mathfrak{M}$ whose projections to the $(r,z)$-plane admit various symmetries.

The following theorem shows the existence of non-resonant Hopf links for every $\alpha>0$, and $\mathfrak{e}>0$ sufficiently small.

\begin{thm}\label{thm: isosceles 3bp}
Fix $\alpha>0$. There exists $\mathfrak{e}_0>0$, such that for every $0<\mathfrak{e}<\mathfrak{e}_0$, the energy surface $\mathfrak{M}=H^{-1}(-1)$ carries a non-resonant Hopf link $\gamma_{1,\mathfrak{e}} \cup \gamma_{2,\mathfrak{e}}$. In particular, $\mathfrak{M}$ admits infinitely many periodic orbits.
\end{thm}

\begin{proof}[Proof of Theorem \ref{thm: isosceles 3bp}] The proof involves a careful change of coordinates that puts $H$ into normal form, see section \ref{sec: construction}. Then we apply the results in section \ref{sec 4}.

Denote $\zeta=(p_r,p_z,r,z)$, replace $r$ with $r+\frac{\alpha\varpi^2}{4+\alpha}$, and add $\frac{(4+\alpha)^2}{2\alpha^2\varpi^2}$ to $H$. We obtain the new Hamiltonian
\begin{eqnarray*}
\bar H(\zeta)=\frac{p_r^2+p_z^2}{2}+\frac{(4+\alpha)^2}{2\alpha^2\varpi^2}+\frac{\varpi^2}{2(r+\frac{\alpha\varpi^2}{4+\alpha})^2}
-\frac{1}{r+\frac{\alpha\varpi^2}{4+\alpha}}-\frac{4\alpha^{-1}}{\sqrt{(r+\frac{\alpha\varpi^2}{4+\alpha})^2+(1+2\alpha)z^2}}.
\end{eqnarray*}
The origin $0\in \bar H^{-1}(0)$ is now the global minimum of $\bar H$ and $\bar H(0)=0$.

After a linear symplectic change of coordinates determined by
$$
S=\mathrm{diag}\bigg(\frac{4+\alpha}{\alpha\varpi^{3/2}},\frac{\sqrt{2}(4+\alpha)^{3/4}(1+2\alpha)^{1/4}}{\alpha\varpi^{3/2}},\frac{\alpha\varpi^{3/2}}{4+\alpha},\frac{\alpha\varpi^{3/2}}{\sqrt{2}(4+\alpha)^{3/4}(1+2\alpha)^{1/4}}\bigg),
$$
we obtain the Hamiltonian
\begin{equation*}
\begin{aligned}
\hat H(\zeta)=\bar H(S\zeta)&=\frac{(4+\alpha)^2p_r^2}{2\alpha^2\varpi^{3}}+\frac{(4+\alpha)^{3/2}(1+2\alpha)^{1/2}p_z^2}{\alpha^2\varpi^{3}}+\frac{(4+\alpha)^2}{2\alpha^2\varpi^2}
+\frac{(4+\alpha)^2}{2\alpha^2\varpi(r+\sqrt{\varpi})^2}\\
&\quad\ -\frac{4+\alpha}{\alpha\varpi^{3/2}(r+\sqrt{\varpi})}-\frac{4(4+\alpha)}{\alpha^2\varpi^{3/2}((r+\sqrt{\varpi})^2+\frac{1}{2}\sqrt{(4+\alpha)(1+2\alpha)}z^2)^{1/2}}.
\end{aligned}
\end{equation*}

Now let $K(\zeta):=\hat H(\zeta)\cdot \frac{2\alpha^2\varpi^3}{(4+\alpha)^2}$. Then $\mathfrak{M}=H^{-1}(-1)$ corresponds to $SK^{-1}(\varpi\mathfrak{e}^2)$.
Replace the coordinates $(p_r,p_z,r,z)$ by  $(y_1,y_2,x_1,x_2)$. The Taylor expansion of $K$ at $0$ gives
\begin{equation}\label{tilde H}
\begin{aligned}
K(y,x)&=2\cdot\frac{y_1^2+x_1^2}{2}+2\sqrt{\frac{4+8\alpha}{4+\alpha}}\cdot\frac{y_2^2+x_2^2}{2}
-\frac{2x_1^3}{\sqrt{\varpi}}-\frac{3}{\sqrt{\varpi}}\sqrt{\frac{4+8\alpha}{4+\alpha}}x_1x_2^2\\
&\quad\  +\frac{3x_1^4}{\varpi}+\frac{6}{\varpi}\sqrt{\frac{4+8\alpha}{4+\alpha}}x_1^2x_2^2-\frac{3(1+2\alpha)x_2^4}{4\varpi}+\mathcal{O}_5.
\end{aligned}
\end{equation}
At this point, we already have the second-order term in normal form
$$
\Gamma^{(2)}=\frac{\alpha_1}{2}(y_1^2+x_1^2)+\frac{\alpha_2}{2}(y_2^2+x_2^2),
$$
where $\alpha_1=2$ and $\alpha_2=2\sqrt{(4+8\alpha)/(4+\alpha)}>2, \forall \alpha>0$. In particular, $\alpha_1,\alpha_2$ are rationally independent if
$\alpha\neq 4(t^2-1)/(8-t^2), \forall t\in \Q.$ The Hamiltonian $K$ is weakly non-resonant at $0$ if $\alpha_2/\alpha_1\notin\mathbb{N}$, that is $\alpha\neq 4(c^2-1)/(8-c^2), \forall c\in\mathbb{N}$. Here, the resonant sequence of $K$ is
$$
\frac{\alpha_2}{\alpha_1}=\sqrt{\frac{4+8\alpha}{4+\alpha}}=2,\ \frac{3}{2},\ \frac{5}{2},\ \frac{4}{3},\cdots\in (1,2\sqrt{2}),
$$
and $\alpha_2/\alpha_1=2$ ($\alpha=3$) is the unique value of $\alpha$ for which $K$ is not weakly non-resonant.

Now we compute the Birkhoff-Gustavson normal form $K_4$ of $K$. By following the construction in section \ref{sec: construction} closely, we first compute the generating function
$W_3(\eta,x)=x_1\eta_1+x_2\eta_2+G_3(\eta,x)$ so that the induced canonical transformation $\Phi_3:(\eta_1,\eta_2,\xi_1,\xi_2)\rightarrow (y_1,y_2,x_1,x_2)$ satisfies
\begin{equation}\label{equation of Wb}
\hat K_3=K\circ \Phi_3=\Gamma^{(2)}+\Gamma^{(3)}+\hat K^{(4)}_3+\mathcal{O}_5,\quad D\cdot \Gamma^{(3)}=0,
\end{equation}
where $D=\sum_{i=1}^2\alpha_i(\eta_i\partial_{\xi_i}-\xi_i\partial_{\eta_i})$ and $\hat K^{(4)}_3$ is the $4$-order homogeneous part of $\hat K_3$. In particular, we have
\begin{equation}\label{equation of W^3}
-D\cdot G_3+\Gamma^{(3)}=K^{(3)}=
-\frac{2x_1^3}{\sqrt{\varpi}}-\frac{3}{\sqrt{\varpi}}\sqrt{\frac{4+8\alpha}{4+\alpha}}x_1x_2^2,
\end{equation}
where $K^{(3)}$ is the $i$-th order part of $K$. In the equation above, we formally replace $(y,x)$ with $(\eta,x)$.

We can solve for $G_3$ as follows. Using complex variables $z_j=x_j+iy_j,j=1,2$, we have
$$
\begin{aligned}
K^{(3)}(z,\bar z)&=-\frac{z_1^3}{4 \sqrt\varpi} - \frac{3|z_1|^2 z_1}{4 \sqrt\varpi} - \frac{3 |z_1|^2 \bar z_1}{4 \sqrt\varpi} - \frac{\bar z_1^3}{4 \sqrt\varpi} \\
&\quad\ -\sqrt{\frac{4 + 8 \alpha}{4 + \alpha}}\bigg(\frac{3z_1 z_2^2}{8\sqrt\varpi} + \frac{3\bar z_1 z_2^2}{8\sqrt\varpi} + \frac{
 3 |z_2|^2 z_1}{4 \sqrt\varpi} + \frac{
 3 |z_2|^2 \bar z_1}{4\sqrt\varpi} + \frac{
 3 z_1 \bar z_2^2}{8\sqrt\varpi} + \frac{
 3 \bar z_1 \bar z_2^2}{8\sqrt\varpi}\bigg).
\end{aligned}
$$
If $\alpha \neq 3$, then $|m_1|>m_2\geq 2$, which means that $\sigma=z_2^{m_2}\bar z_1^{|m_1|}$ has degree $\geq 5$. Hence, from \eqref{hatH}, we obtain $\Gamma^{(3)}(\eta,x)=0$. If $\alpha=3$, then $m_2=1$ and $|m_1|=2$. In this case, $\sigma=z_2\bar z_1^2$, $\bar \sigma = z_1^2 \bar z_2$, and $\Gamma^{(3)}$ is a combination of $\sigma$ and $\bar \sigma$. Since $K^{(3)}$ has no such terms, we conclude that $\Gamma^{(3)}(\eta,x)=0$, as in the previous cases. Hence $-D \cdot G_3 = K^{(3)}$. Using \eqref{D on monomials}, we solve this equation for $G_3$ and obtain
\begin{equation}\label{expression of W^3}
\begin{aligned}
G_3(\eta,x)&=-\frac{\eta_1(2\eta_1^2+3x_1^2)}{3\sqrt{\varpi}}
-\frac{12(1+2\alpha)\cdot\eta_2x_1x_2}{(12+31\alpha)\sqrt{\varpi}}\\
&\quad\ -\frac{3\sqrt{1+2\alpha}\cdot \eta_1(8(1+2\alpha)\eta_2^2+(4+15\alpha)x_2^2)}{(12+31\alpha)\sqrt{\varpi(4+\alpha)}}.
\end{aligned}
\end{equation}
The $4$-homogeneous part $\hat K_3^{(4)}$ of $K_3$ is given by
$$
\hat K_3^{(4)}(\eta,x) = \sum^2_{i=1}\frac{\alpha_i}{2}\big((\partial_{x_i}G_3)^2-(\partial_{\eta_i}G_3)^2\big)+K^{(4)},
$$
where
$$
K^{(4)}=\frac{3x_1^4}{\varpi}+\frac{6}{\varpi}\sqrt{\frac{4+8\alpha}{4+\alpha}}x_1^2x_2^2-\frac{3(1+2\alpha)x_2^4}{4\varpi},
$$
\begin{eqnarray*}
&\partial_{x_1}G_3=-\frac{2\eta_1 x_1}{\sqrt{\varpi}}-\frac{12(1+2\alpha)\eta_2 x_2}{(12+31\alpha)\sqrt{\varpi}},\ \partial_{\eta_1}G_3=-\frac{2\eta_1^2+x_1^2}{\sqrt{\varpi}}-\frac{3\sqrt{1+2\alpha}(8(1+2\alpha)\eta_2^2+(4+15\alpha)x_2^2)}{(12+31\alpha)\sqrt{\varpi(4+\alpha)}},\\
&\partial_{x_2}G_3=-\frac{12(1+2\alpha)\cdot\eta_2 x_1}{(12+31\alpha)\sqrt{\varpi}}-\frac{6\sqrt{1+2\alpha}(4+15\alpha)\eta_1x_2}{(12+31\alpha)\sqrt{\varpi(4+\alpha)}},\
 \partial_{\eta_2}G_3=-\frac{12(1+2\alpha)(\sqrt{4+\alpha}\cdot x_1x_2+4\sqrt{1+2\alpha}\eta_1\eta_2)}{(12+31\alpha)\sqrt{\varpi(4+\alpha)}}.
\end{eqnarray*}
The new Hamiltonian $\hat K_3= K \circ \Phi_3$ in coordinates $(\eta,\xi)$ becomes
$$
\hat K_3(\eta,\xi)= \Gamma^{(2)}(\eta,\xi) +\hat K_3^{(4)}(\eta,\xi) + \mathcal{O}_{5}.
$$

Next, we rename $(\eta,\xi)$ back to $(y,x)$ and write $\hat K_3=\hat K_3(y,x)$. We aim to find a canonical transformation $\Phi_4:(\eta_1,\eta_2,\xi_1,\xi_2)\rightarrow (y_1,y_2,x_1,x_2)$ induced by the generating function
$W_4(\eta,x)=\eta_1 x_1+\eta_2 x_2+G_4(\eta,x)$,  so that the new Hamiltonian is in the normal form up to order $4$
$$
\hat K_4=\hat K_3\circ \Phi_4=\Gamma^{(2)}+\Gamma^{(4)}+\mathcal{O}_5,\quad D\cdot \Gamma^{(4)}=0.
$$
The $4$-homogeneous polynomials $\Gamma^{(4)}$ and $G_4$ satisfy
\begin{equation}\label{equation of W^4}
-D\cdot G_4+\Gamma^{(4)}(\eta,x)=\hat K^{(4)}_3(\eta,x).
\end{equation}
In fact, we only need to see $\Gamma^{(4)}$, which consists of the part of $\hat K_3^{(4)}$ in $\ker D$. In complex coordinates $z_j=x_j+iy_j,j=1,2,$ the expression for $\hat K_3^{(4)}(y,x)$ becomes
$$
\begin{aligned}
\hat K_3^{(4)}(z,\bar z) &= -\frac{3}{4\varpi}|z_1|^4-\frac{3(24+55\alpha)}{2(12+31\alpha)\varpi}\sqrt{\frac{1+2\alpha}{4+\alpha}}|z_1|^2|z_2|^2\\
&\quad -\frac{9(1+2\alpha)(128+380\alpha+31\alpha^2)}{32(4+\alpha)(12+31\alpha)\varpi}|z_2|^2 + Q,
\end{aligned}
$$
where $Q$ consists of the part of $\hat K_3^{(4)}$ in $\mathrm{Im}D$.

As before, if $\alpha \neq 3$, $\sigma$ has order $\geq 5$ and thus $\Gamma^{(4)}$ is formed by the monomials of $\hat K_3^{(4)}$ of the form $|z_1|^4,$ $|z_1|^2|z_2|^2$ and $|z_2|^4$. If $\alpha= 3$ ($|m_1|=2$ and $m_2=1$), then $\sigma = z_2\bar z_1^2$ and thus the same conclusions hold for $\Gamma^{(4)}$.  It follows that
\begin{equation}\label{Gamma 4}
\begin{aligned}
\Gamma^{(4)}(\eta,x)&=-\frac{3}{4\varpi}(\eta_1^2+x_1^2)^2-\frac{3(24+55\alpha)}{2(12+31\alpha)\varpi}\sqrt{\frac{1+2\alpha}{4+\alpha}}(\eta_1^2+x_1^2)(\eta_2^2+x_2^2)\\
&\quad -\frac{9(1+2\alpha)(128+380\alpha+31\alpha^2)}{32(4+\alpha)(12+31\alpha)\varpi}(\eta_2^2+x_2^2)^2.
\end{aligned}
\end{equation}
This implies that $\nu=2$, where $\nu\geq 2$ is the least integer so that at least one of the terms in \eqref{condition of nv} is non-zero.  Indeed, we have
$$
a_{2,0,2,0}=-\frac{3}{4\varpi},\ a_{1,1,1,1}=-\frac{3(24+55\alpha)}{2(12+31\alpha)\varpi}\sqrt{\frac{1+2\alpha}{4+\alpha}},\ a_{0,2,0,2}=-\frac{9(1+2\alpha)(128+380\alpha+31\alpha^2)}{32(4+\alpha)(12+31\alpha)\varpi},
$$
and for every $\alpha>0$, we have

\begin{equation}\label{nonresonant formula}
\begin{aligned}
\frac{\Omega_{\nu,1}}{\alpha_1\alpha_2}&=\frac{a_{1,1,1,1}\alpha_1-\nu a_{2,0,2,0}\alpha_2}{\alpha_1^2\alpha_2}=\frac{21\alpha}{16(12+31\alpha)\varpi}\neq 0,\\
\frac{\Omega_{\nu,2}}{\alpha_1\alpha_2}&=\frac{a_{1,1,1,1}\alpha_2-\nu a_{0,2,0,2}\alpha_1}{\alpha_1\alpha_2^2}=\frac{3\alpha(260+93\alpha)}{256(12+31\alpha)\varpi}\neq 0,\\
\Omega_\nu&=\frac{\Omega_{\nu,1}+\Omega_{\nu,2}}{\alpha_1\alpha_2}
=\frac{279\alpha(4+\alpha)}{256(12+31\alpha)\varpi}\neq 0.
\end{aligned}
\end{equation}

\begin{rem}
The normal form computed in this section coincides with the one in \cite{Shibayama2009} after replacing $\alpha$ with $\lambda=\sqrt{\frac{4+8\alpha}{4+\alpha}}$.
\end{rem}

If $\alpha\neq 3$ (that is $|m_1|>m_2\geq 2$), 
since $\Omega_{\nu,1},\Omega_{\nu,2},\Omega_{\nu}\neq 0$, it follows from Theorem \ref{thm: main theorem}-(i),(ii) that the sphere-like component $\Sigma_E$ of the energy surface carries a pair of periodic orbits $\gamma_{1,E},\gamma_{2,E}$ forming a non-resonant Hopf link for every $E>0$ sufficiently small. In particular, from the proof of Theorem \ref{thm: m2>2} and \ref{thm: m2=2}, we obtain that
$$
(\rho_{1,E}-1)(\rho_{2,E}-1)=1+\frac{279\alpha(4+\alpha)E}{64(12+31\alpha)\varpi}+\mathcal{O}(E^{2})\neq 1,\quad \forall E\ll 1.
$$
Moreover, if $\alpha=3$ (that is $m_2=1$ and $|m_1|=2$), then we have
\begin{eqnarray*}
K(\zeta)=p_r^2+r^2+2(p_z^2+z^2)
-\frac{2r^3}{\sqrt{\varpi}}-\frac{6rz^2}{\sqrt{\varpi}}+\frac{3r^4}{\varpi}+\frac{12r^2z^2}{\varpi}-\frac{21z^4}{4\varpi}+\mathcal{O}_5.
\end{eqnarray*}
The normal form becomes
$$
\begin{aligned}
\Gamma_4=\Gamma^{(2)}+\Gamma^{(4)}&=2\cdot \frac{\eta_1^2+\xi_1^2}{2}+4\cdot \frac{\eta_2^2+\xi_2^2}{2}\\
&\quad\ -\frac{3}{4\varpi}(\eta_1^2+\xi_1^2)^2-\frac{27}{10\varpi}(\eta_1^2+\xi_1^2)(\eta_2^2+\xi_2^2)+\frac{663}{160\varpi}(\eta_2^2+\xi_2^2)^2.
\end{aligned}
$$
Since $K$ is a smooth function of $p_z^2,z^2$, see \eqref{tilde H}, the condition \eqref{cond for y2,x2} is naturally satisfied. Then due to $\Omega_{2,1},\Omega_{2,2},\Omega_2\neq 0$,  Theorem \ref{thm: main theorem 2}-(v), the sphere-like component $\Sigma_E$ of the energy surface also carries a pair of periodic orbits $\gamma_{1,E},\gamma_{2,E}$ forming a non-resonant Hopf link for every $E>0$ sufficiently small.
In particular, from the proof of Theorem \ref{thm: m2=1, m1>1}-(v), we obtain that
$$
(\rho_{1,E}-1)(\rho_{2,E}-1)=1+\frac{279E}{320\varpi}+\mathcal{O}(E^{2})\neq 1,\quad \forall E\ll 1.
$$
The proof is now completed.
\end{proof}


\begin{rem}\label{rem: alpha=0}
If $\alpha=0$, then $|m_1|=m_2=1$ and
$$
\tilde H(\zeta)=p_r^2+p_z^2+r^2+z^2-\frac{2r^3}{\sqrt{\varpi}}-\frac{3rz^2}{\sqrt{\varpi}}+\frac{3r^4}{\varpi}+\frac{6r^2z^2}{\varpi}-\frac{3z^4}{4\varpi}+\mathcal{O}_5.
$$
The Birkhoff-Gustavson normal form becomes
$\Gamma^{(2)}+\Gamma^{(4)}=\eta_1^2+\eta_2^2+\xi_1^2+\xi_2^2-\frac{3}{4\varpi}(\eta_1^2+\eta_2^2+\xi_1^2+\xi_2^2)^2.$
We obtain $\Omega_{\nu,1}=\Omega_{\nu,2}=\Omega_{\nu}=0$ and thus Theorem \ref{thm: m2=m1=1} does not apply.
\end{rem}

\section{Hill's lunar problem}\label{Hill}
The Hamiltonian of Hill's lunar problem is given in canonical coordinates $(p,q)$ by
\begin{equation}\label{Hill's problem}
H(p,q)=\frac{1}{2}(p_1^2+p_2^2)+p_1q_2-p_2q_1-\frac{3}{|q|}-\frac{3q_1^2}{2}+\frac{1}{2}(q_1^2+q_2^2).
\end{equation}
We observe that $H$ is invariant under the $\mathbb{Z}_4$-action generated by linear map $R:(y_1,y_2,x_1,x_2)\mapsto (-y_2,y_1,-x_2,x_1)$.

In Levi-Civita coordinates $q=v^2,p=2u/\bar v$, for which $dp\wedge dq=4du\wedge  dv$, the Hamiltonian becomes
$$
\begin{aligned}
H(u,v)&=\frac{2|u|^2}{|v|^2}+2(u_1v_2-u_2v_1)-\frac{3}{|v|^2}-(v_1^2-v_2^2)^2+2v_1^2v_2^2,
\end{aligned}$$
Let $K(u,v)=\frac{|v|^2}{4}(H+c)$ be the regularized Hamiltonian. After re-scaling coordinates $(u,v)=(2^{1/4}c^{3/4}y,2^{3/4}c^{1/4}x)$, we define the new Hamiltonian
\begin{equation}\label{equ: regularized Hill problem}
\begin{aligned}
\hat H(y,x)&=K(2^{1/4}c^{3/4}y,2^{3/4}c^{1/4}x)/(2^{1/2}c^{3/2})+3/(4\cdot 2^{1/2}c^{3/2})\\
&=\frac{1}{2}(|y|^2+|x|^2)+2|x|^2(iy\cdot x)-4|x|^6+24|x|^2x_1^2x_2^2\\
&=\frac{1}{2}\big((y_1+f_1)^2+(y_2+f_2)^2\big)+V(x),
\end{aligned}
\end{equation}
where $f_1=2x_2|x|^2,f_2=-2x_1|x|^2$ and $V=|x|^2(1/2-6(x_1^2-x_2^2)^2)$. Notice that the last constant term in the definition of $\hat H$ is added so that $\hat H$ does not depend on the parameter $c$. In this way, the regularized dynamics on $H^{-1}(-c), c\gg 0,$ corresponds to
the dynamics on the sphere-like component of $\hat H^{-1}(3/( 2^{5/2}c^{3/2})).$ Thus, the dynamics of $H$ for large negative energies corresponds to the dynamics of $\hat H$ for positive energies close to $0$.

For simplicity, $\hat H$ will be denoted by $H$. Then $H$ splits into homogeneous terms as
\begin{equation}
H (y_1 , y_2, x_1 , x_2) = H^{(2)} + H^{(4)} + H^{(6)} ,
\label{HH}
\end{equation}
where
$H^{(2)} = \frac{1}{2} (|y|^2 + |x|^2),$
$H^{(4)} = 2|x|^2(iy\cdot x)$ and $H^{(6)} = -4|x|^6+24|x|^2x_1^2x_2^2.$
The Hamiltonian equations of $H$ become
$$
\left\{
\begin{aligned}
\dot y_1&=-\partial_{x_1}H=-(y_1+f_1)\partial_{x_1}f_1-(y_2+f_2)\partial_{x_1}f_2-\partial_{x_1}V,\\
\dot y_2&=-\partial_{x_2}H=-(y_1+f_1)\partial_{x_2}f_1-(y_2+f_2)\partial_{x_2}f_2-\partial_{x_2}V,\\
\dot x_1&=\partial_{y_1}H=y_1+f_1,\\
\dot x_2&=\partial_{y_2}H=y_2+f_2.
\end{aligned}
\right.
$$

Hill's lunar problem can be seen as an approximation to the circular planar restricted three-body problem in the case that the mass of the primary body is much larger than the mass of the secondary body, like the sun (primary) and the Earth (secondary). Hill's approximation aims to describe the motion of a body with a negligible mass close to the secondary body (Earth satellite). The Hamiltonian in \eqref{HH} was taken from \cite{ST}, where an interesting numerical study of Hill's problem is given. In \cite{ST}, the reader finds details on the spatial and time scales that were chosen to obtain a Hamiltonian function that does not depend
on any physical parameter. The energy $E$ of $H$ is related to Hill's Jacobi constant $c_H=-2c$ through $E=\frac{3}{2}|c_H|^{-3/2}$.

Since our analysis restricts to a neighborhood of the origin, it is interesting to say that $(q_1,q_2)=(0,0)$ corresponds to the collision of the negligible-mass body (like the satellite)
with the secondary body (like the Earth).
The usual Levi-Civita regularization of the collisions allows for motion to
be smoothly continued through and after the collision.

The three terms $H^{(2)}$, $H^{(4)}$, and $H^{(6)}$
have a physical simple interpretation that helps in the mathematical analysis. The $H^{(2)}$ part corresponds to the Keplerian motion of the satellite around the Earth.
Notice that all orbits of $H^{(2)}$ are periodic and have the same period. The term $H^{(4)}$ corresponds to a ``fictitious'' Coriolis force that appears due to the following. In the derivation of the restricted circular planar three-body equations, it is supposed that the Earth rotates around the Sun along a circular orbit, therefore, with constant
angular velocity. Then, we can choose a spatial reference frame that rotates with the Earth, such that in this frame, both the Sun and the Earth are at rest. The choice of this reference frame allows for the elimination of an explicit time dependence on the equations of motion of the satellite. Nevertheless, the term $H^{(4)}$ appears as the result of this non-inertial reference frame. Therefore, $H^{(2)}+H^{(4)}$ still represents the Kepler motion of the satellite around the Earth but in a rotating reference frame. Finally, the term $H^{(6)}$ accounts for the perturbation of the motion due to the presence of the Sun. This is the term that breaks the integrability of the system and is responsible for the complicated dynamics. The influence of these terms vanishes in the limit as $|q|\to 0$. This explains why, in the following normal form analysis, we are forced to go to the sixth order.

Same as the steps in section \ref{sec: construction}, one can also establish the Birkhoff-Gustavson normal form of $H$ given in equation (\ref{HH}). Firstly, we see that $H_2=H^{(2)}$, since $H^{(2)}$ is already in the normal form. Due to the specific properties of Lunar problem, we introduce a different way to obtain the Birkhoff-Gustavson normal form up to order six.

The first step in establishing the Birkhoff-Gustavson normal form of $H$ given in equation (\ref{HH}) is to find a canonical transformation that normalizes the Hamiltonian up to order four.
Using the procedure described in \cite{moser68} we are led to the following generating function
\[
W(y,\xi)=y_1\xi_1+y_2\xi_2+\underbrace{(\xi_2y_1-\xi_1y_2)(-y_1\xi_1-y_2\xi_2)}_{\defi G(y,\xi)}=y\cdot \xi-(iy\cdot \xi)(y\cdot \xi),
\]
where $(\eta_1,\eta_2,\xi_1,\xi_2)$ are the new variables. These variables are implicitly determined by equations
\[
x_j=\xi_j+\partial_{y_j}G,\qquad y_j=\eta_j-\partial_{\xi_j}G, \qquad j=1,2.
\]
which can be solved up to fifth order as:
\begin{eqnarray*}
y&=&\eta+i\eta A+\eta F+2\eta(F^2-A^2)+4i\eta FA,\\
x&=&\xi+i\xi A-\xi F-\xi(F^2-A^2)+2i\xi FA,
\end{eqnarray*}
where
$$
F=i\eta\cdot \xi,\qquad A=\eta\cdot \xi.
$$
In the new coordinates, the Hamiltonian $H$ becomes
\begin{equation}
\hat H_4(\eta,\xi)=H_2+2H_2F-7|\xi|^2 F^2+5H_2F^2-3H_2A^2+3A^2|\xi|^2-4|\xi|^6+24|\xi|^2\xi_1^2\xi_2^2+\mathcal{O}_7,\label{H2}
\end{equation}
In particular, the new Hamiltonian is already in the normal form up to order four.

The next step is to find the Birkhoff-Gustavson normal form up to order six. We do this last step of the normalization in a different but equivalent way: averaging. The procedure is explained, for instance, in \cite{anrMMCM}. We first introduce the canonical coordinates
$$\xi_j =\sqrt{2\tilde I_j}\cos \tilde\theta_j,\quad \eta_j =\sqrt{2\tilde I_j}\sin \tilde\theta_j,\quad j=1,2.$$
In these coordinates $ H=\tilde I_1+\tilde I_2+\mathcal{O}_3$. Then we define the new canonical coordinates
$I,\theta$, as
$$ I_1= \tilde I_1+\tilde I_2,\quad I_2=\tilde I_2,\quad \theta_1=\tilde \theta_1,\quad \theta_2=\tilde \theta_2-\tilde \theta_1.$$
In the new coordinates the Hamiltonian function
has the same form as in (\ref{H2}) with:
\begin{align*}
H_2=I_1,\qquad\quad F=-2\sqrt{I_2(I_1-I_2)}\sin\theta_2,\qquad \ A=(I_1-I_2)\sin (2\theta_1)+I_2\sin(2\theta_1+2\theta_2),\\
\xi_1^2=2(I_1-I_2)\cos^2\theta_1,\ \xi_2^2=2I_2\cos^2(\theta_1+\theta_2),\ |\xi|^2=2(I_1-I_2)\cos^2\theta_1+2I_2\cos^2(\theta_1+\theta_2).
\end{align*}
We consider the order-$6$ part of the Hamiltonian \eqref{H2} in complex coordinates $w_j=\xi_j + i \eta_j, j=1,2$, and identify the terms that are in normal form. To do that, we average each term with respect to $\theta_1$. For a monomial $Q= w_1^{k_1}\bar w_1^{k_2} w_2^{l_1}\bar w_2^{l_2}$, we obtain
$$
\begin{aligned}
\int_0^{2\pi}w_1^{k_1}\bar w_1^{k_2} w_2^{l_1}\bar w_2^{l_2}d\theta_1&=(2I_2)^{\frac{l_1+l_2}{2}}(2I_1-2I_2)^{\frac{k_1+k_2}{2}}\int_0^{2\pi}e^{i(k_1-k_2)\theta_1+i(l_1-l_2)(\theta_1+\theta_2)}d\theta_1.
\end{aligned}
$$
We see that the average of $Q$ vanishes if and only if $k_1-k_2\neq -(l_1-l_2)$, or equivalently, if and only if $Q$ is a normalized term. Hence, after integration, the normalized terms are identified, and $\Gamma^{(6)}$ is determined. Proceeding in this way, the Birkhoff-Gustavson normal form of $H$ up to order six becomes $H=H_6 + \mathcal{O}_7 = \Gamma^{(2)} +\Gamma^{(4)} + \Gamma^{(6)} + \mathcal{O}_7$, where
$$
\Gamma^{(2)}=H_2,\quad \Gamma^{(4)}=2H_2F,\quad \Gamma^{(6)}= -8 H_2 F^2 -10 H_2^3 +60 H_2I_2(I_1-I_2).
$$
In complex coordinates $w_j, j=1,2,$, we have
$$
\begin{aligned}
\Gamma^{(2)}&=\frac{1}{2}(|w_1|^2+|w_2|^2),\quad \Gamma^{(4)}=\frac{1}{2i}(|w_1|^2+|w_2|^2)(w_1\ov w_2-\ov w_1 w_2),\\
\Gamma^{(6)}&= -8 H_2 F^2 -10 H_2^3 +15 H_2 |w_2|^2|w_1|^2.
\end{aligned}
$$

An analysis of Hamilton's equations of $H_6$ shows
the existence of two periodic orbits $\gamma_1$ and $\gamma_2$ (as stated in Theorem \ref{existence})
such that the $w_1$ and $w_2$ components of  $\gamma_1$ satisfy $w_2=iw_1$ and the  $w_1$ and $w_2$ components of $\gamma_2$ satisfy $w_2=-iw_1$. They are the direct and retrograded orbits, respectively. We see that the two orbits are orthogonal with respect to the canonical Hermitian form on $\C^2$.
This suggests the choice of another set of variables $z_1,z_2$, determined by
\begin{equation}\label{equ: new coordinates}
\left(\begin{array}{l}
w_1\\
w_2 \end{array}
\right)=\left(\begin{array}{ll}
1/\sqrt{2}& \ 1/\sqrt{2}\\
i/\sqrt{2}& -i/\sqrt{2}
\end{array} \right)
\left(\begin{array}{l}
z_1\\
z_2 \end{array}
\right).
\end{equation}
This transformation leaves invariant both $H_2$ and the symplectic form. In the new coordinates $(z_1,z_2)$, we see that the $z_2$ component of $\gamma_1$ is zero and the $z_1$ component of $\gamma_2$ is zero. 
In real coordinates, the transformation \eqref{equ: new coordinates} can be written as $(\eta_1,\eta_2,\xi_1,\xi_2)=2^{-1/2}(v_1+v_2,u_1-u_2,u_1+u_2,v_2-v_1)$, where $z_j=u_j+iv_j,j=1,2$. Then in coordinates $(v,u)$, the symplectic map $R$ becomes $\mathcal{R}:(v_1,v_2,u_1,u_2)\rightarrow (-u_1,u_2,v_1,-v_2)$, i.e. rotates $z_1$ to $e^{-2\pi i/4}z_1$ and $z_2$ to $e^{2\pi i/4}z_2$. Therefore, $\gamma_1$ and $\gamma_2$ are both $\mathcal{R}$-invariant.
This change of variables takes the Birkhoff-Gustavson normal form to another one (the Birkhoff-Gustavson normal form is not unique).
Then we have
\begin{equation}\label{equ: hat H}
\begin{aligned}
\Gamma^{(2)}&=\frac{1}{2}(|z_1|^2+|z_2|^2),\quad \Gamma^{(4)}=\frac{1}{2}(-|z_1|^4+|z_2|^4),\\
\Gamma^{(6)}&=-\frac{1}{2}(|z_1|^2+|z_2|^2)\cdot\big(\frac{3}{4}(|z_1|^4+|z_2|^4)+|z_1|^2|z_2|^2+
\frac{15}{4}\big((\ov z_1 z_2)^2+(z_1\ov z_2)^2\big)\big).
\end{aligned}
\end{equation}
This implies $\nu=2$. By Proposition \ref{prop: zp symmetry} with $p=4$, we know that $\hat H$ and the Birkhoff-Gustavson $\hat H_N$ are both $\mathcal{R}$-invariant for every $N\geq 3$.

Denote $\gamma_1$ and $\gamma_2$ as $(z_1,z_2)(t)=(c_1e^{-i\om_1t},0)$ and $(z_1,z_2)(t)=(0,c_2e^{-i\om_2t})$, respectively, where $c_1,c_2>0$ are sufficiently small. Then $\gamma_1$ and $\gamma_2$ solve the equations
$$
\dot{\bar z}_1=i\bar z_1\left(1-2c_1^2-\frac{9}{4}c_1^4\right),\quad
\dot{\bar z}_2=i\bar z_2\left(1+2c_2^2-\frac{9}{4}c_2^4\right),
$$
respectively. Moreover, we have $\omega_1=1-2c_1^2-\frac{9}{4}c_1^4$ and $\omega_2=1+2c_2^2-\frac{9}{4}c_2^4$.

For the periodic solution $\gamma_1$, we have $E=H_6(c_1^2,0)=\frac{c_1^2}{2}(1-c_1^2-\frac{3}{4}c_1^4)$. Since $\partial_{|z_1|^2}H_6(0,0)=1/2$, there
exists a function $c_1^2=c_1(E)^2$ for $E$ sufficiently small such that $E=H_6(c_1(E)^2,0)$. We compute that
$$
c_1(E)^2=2E+4E^2+22E^3+O(E^4),\quad \om_1(E)=1-4E-17E^2+O(E^3).
$$
Moreover, along $\gamma_1$, we have $\hat \omega_2=2\partial_{z_2\bar z_2}H_6(\gamma_1)=1-\frac{7}{4}c_1^4=1-7E^2+O(E^3)$. Then
$$\frac{\hat \omega_2}{\omega_1}=1+4\Omega_{2,1}E+8\beta_1E^2+O(E)^3,$$
where $\Omega_{2,1}=1$ and $\beta_1=13/4$.

For the periodic solution $\gamma_2$,  we have $E=H_6(0,c_2^2)=\frac{c_2^2}{2}(1+c_2^2-\frac{3}{4}c_2^4)$. The computation gives
$$
c_2(E)^2=2E-4E^2+22E^3+O(E^4),\qquad \om_2(E)=1+4E-17E^2+O(E^3).
$$
Moreover, along $\gamma_2$, we have $\hat \omega_1=2\partial_{z_1\bar z_1}H_6(\gamma_2)=1-\frac{7}{4}c_2^4=1-7E^2+O(E^3).$ Then
$$\frac{\hat \omega_1}{\omega_2}=1+4\Omega_{2,2}E+8\beta_2E^2+O(E)^3,$$
where $\Omega_{2,2}=-1$ and $\beta_2=13/4$. Finally, by using \eqref{o2o1/w1w2}, we compute that
$$
\frac{\hat \omega_2\hat \omega_1}{\omega_1\omega_2}=1+8(\frac{13}{2}-2)E^2+O(E^3)=1+36E^2+O(E^3).
$$

Let $\rho_j$ denote the rotation number of the periodic orbits $\gamma_j,j=1,2$. Since $\alpha_1=\alpha_2=1,\nu=2$ and $a_{0,2,2,0}=0$, we apply \eqref{rot of gamma1 p=2 1}, \eqref{rot of gamma2 q=2 1} and Theorem \ref{thm: non-resonant theorem C} to obtain  $\rho_1=2+4E+26E^2+O(E^3),  \rho_2=2-4E+26E^2+O(E^3)$ and $(\rho_{1}-1)(\rho_{2}-1)=1+36E^2+O(E^3)\neq 1$, for every $E>0$ sufficiently small. A direct application of Theorem \ref{thm: main theorem 3}-(ii) gives the following result.

\begin{thm}\label{thm: lunar problem}
For every $E>0$ sufficiently small, the regularized Hill's lunar problem \eqref{equ: regularized Hill problem} admits a non-resonant Hopf link
$\gamma_{1,E}\cup\gamma_{2,E}\subset \Sigma_E \subset \hat H^{-1}(E)$. Their respective rotation numbers satisfy $\rho_{1,E}=2+4E+26E^2+O(E^3)$ and $\rho_{2,E} = 2-4E+26E^2+O(E^3)$. In particular, $H^{-1}(-c)$ contains infinitely many periodic orbits.
\end{thm}

The Hopf link $\gamma_{1,E} \cup \gamma_{2,E},$ for $E>0$ sufficiently small, given in Theorem \ref{thm: lunar problem}, corresponds to a pair of periodic orbits $\gamma_{1,c}\cup \gamma_{2,c}\subset H^{-1}(-c)$ of the Hamiltonian \eqref{Hill's problem}, for $c=\frac{1}{2}(\frac{3}{2E})^{2/3}\gg 0$ sufficiently large. Their projections to the $q$-plane are simple closed curves surrounding the singularity at the origin in opposite directions, thus $\gamma_{1,c}$ and $\gamma_{2,c}$ are direct and retrograde orbits, respectively.  Notice that the non-resonance condition on their rotation numbers is degenerate at first
order in $E$. This is just a consequence of the fact that $H^{(2)}+H^{(4)}$ is the Hamiltonian function of the  Keplerian motion (after regularization of collisions and
time rescaling) in a rotating reference frame.

\section{The He\'non-Heiles Hamiltonian}\label{Henon-Heiles}
In this section, we compute the Birkhoff-Gustavson normal form  of the H\'enon-Heiles Hamiltonian
\begin{equation}\label{equ: Henon-Heiles ham}
H=\frac{y_1^2+x_1^2}{2}+\frac{y_2^2+x_2^2}{2}+x_1^2x_2-\frac{x_2^3}{3},
\end{equation}
and find non-resonant Hopf links for small energies. We refer to \cite{ Gustavson1966, KR} for similar computations. Notice that $H$ is invariant under the $\mathbb{Z}_3$-action generated by the unitary map $R:(y,x)\mapsto (e^{2\pi i/3}y,e^{2\pi i/3}x)$. We aim to apply Theorem \ref{thm: main theorem 3} to obtain non-resonant Hopf links for small energies.

We have $\alpha_1=\alpha_2=1$. The Hamiltonian $H$ is in normal form up to order $2$ with $H_2 = \Gamma^{(2)} = \frac{1}{2}(y_1^2+x_1^2)+\frac{1}{2}(y_2^2+x_2^2)$.
For $N\geq 3$, assume that $H = H_{N-1} + \hat H_{N-1}^{(N)} + \mathcal{O}_{N+1}$ is in normal form up to order $N-1$, and $\hat H_{N-1}^{(N)}$ is a $N$-homogeneous polynomial. We follow section \ref{normal1} to obtain the order-$N$ normal form of $H$. Consider the generating function $W_N=x_1\eta_1+x_2\eta_2+G_N$ and the associated canonical transformation $\Phi_N:(\eta,\xi)\mapsto (y,x)$, where $G_N$ is a $N$-homogeneous polynomial in $(\eta,x)$.  There exist unique $G_N\in {\rm im} D$ and $\Gamma^{(N)}\in \ker D$ satisfying $-D\cdot G_N+\Gamma^{(N)}=\hat H^{(N)}_{N-1}(\eta,x)$. 

Denote $z_j=x_j+iy_j, j=1,2.$  Consider a monomial $z^k\bar z^l$ of odd degree $N=2M+1=k_1+k_2+l_1+l_2$. Since $|m_1|=m_2=1$, we see from \eqref{D on monomials} that $D \cdot z^k\bar z^l = (2(M-l_1-l_2)+1)z^k\bar z^l\neq 0$, which implies that $\Gamma^{(2N+1)}=0$ for every odd $N.$

Hence, for $N=3$, we have $\Gamma^{(3)}=0$, $H_3 = H_2$ and $-D\cdot G_3 = H^{(3)}_2$, where
$$
H^{(3)}_2=\frac{1}{8} (z_1 + \bar z_1)^2 (z_2 + \bar z_2) - \frac{1}{24} (z_2 +\bar z_2)^3.
$$
This gives
$$
\begin{aligned}
G_3=\frac{1}{9} \left(6 \eta_1 (\eta_1 \eta_2 + x_1 x_2) - 2 \eta_2^3 + 3 \eta_2 (x_1^2 - x_2^2)\right).
\end{aligned}
$$
The new Hamiltonian becomes
$\hat H_3 = H\circ \Phi_3 = \Gamma^{(2)} + \hat H_3^{(4)} + \mathcal{O}_5, $
where
$$
\begin{aligned}
\hat H^{(4)}_3(y,x)&=\sum_{i=1}^2\frac{1}{2}\left((\partial_{x_i}G_3(y,x))^2-(\partial_{y_i}G_3(y,x))^2\right).
\end{aligned}
$$
An extensive computation gives $\hat H_3^{(4)} = \Gamma^{(4)}(y,x) + \hat G_4(y,x),$ where
$$
\Gamma^{(4)}=-\frac{5}{48}(|z_1|^4+|z_2|^4) +\frac{1}{12}|z_1|^2|z_2|^2 - \frac{7}{48}(z_1^2\bar z_2^2+\bar z_1^2z_2^2)\in \ker D,
$$
and $\hat G_4 \in {\rm im} D$. The new Hamiltonian becomes $\hat H_4 = H \circ \Phi_3 \circ \Phi_4 = H_4 + \mathcal{O}_5= \Gamma^{(2)} + \Gamma^{(4)} + \mathcal{O}_5$. It is not necessary to consider $N>4.$

Consider the symplectic transformation $\Psi(y,x) = 2^{-1/2}(y_1+y_2,x_1-x_2,x_1+x_2,y_2-y_1)$ given in \eqref{Def_Psi}. Then $\hat H_4 \circ \Psi = H_4 \circ \Psi + \mathcal{O}_5$ is in normal form. We compute
$$
\begin{aligned}
H_4\circ \Psi =\frac{1}{2}(|z_1|^2+|z_2|^2)+ \frac{1}{24}|z_1|^4 - \frac{1}{2}|z_1|^2|z_2|^2 + \frac{1}{24}|z_2|^4.
\end{aligned}
$$
We see that $H_4\circ \Psi$ admits a pair of periodic orbits $\gamma_1=(c_1e^{-i\omega_1 t},0)$ and $\gamma_2=(0,c_2e^{-i\omega_2t})$ forming a Hopf link and solving Hamilton's equations
$$
\dot{\bar z}_1=i\bar z_1\left(1+\frac{1}{6}|z_1|^2-|z_2|^2\right),\quad \dot{\bar z}_2=i\bar z_2\left(1+\frac{1}{6}|z_2|^2-|z_1|^2\right).
$$
From the expression above, we obtain  $\nu=2$, $a_{0,2,2,0}=0$ and $\Omega_{2,1}=\Omega_{2,2}=-7/12$, $\Omega_2=-7/6\neq 0$. Finally, Theorem \ref{thm: main theorem 3}-(i) implies that $H$ admits a non-resonant Hopf link $\gamma_1\cup \gamma_2\subset H^{-1}(E)$ for every $E>0$ sufficiently small.  The rotation numbers of $\gamma_1$ and $\gamma_2$ satisfy $\rho_1=\rho_2=2-7E/3+O(E^2)$ and thus $(\rho_1-1)(\rho_2-1)=1-14E/3+O(E^2)$. Hence, $\gamma_1\cup \gamma_2$ is a non-resonant Hopf link. In the original coordinates, the projection of $\gamma_1$ and $\gamma_2$ to the $x$-plane is an embedded closed curve about the origin and traversed in opposite directions.

\begin{thm}\label{thm: Henon-Heiles}
For every $E>0$ sufficiently small, the H\'enon-Heiles Hamiltonian \eqref{equ: Henon-Heiles ham} admits a pair of periodic orbits $\gamma_{1,E}, \gamma_{2,E}\subset \Sigma_E \subset H^{-1}(E)$ forming a non-resonant Hopf link. Their respective rotation numbers satisfy $\rho_{1,E}=\rho_{2,E} = 2-7E/3+O(E^2)$. In particular, $\Sigma_E$ contains infinitely many periodic orbits.
\end{thm}

\appendix

\section{The quaternion trivialization}\label{sec_quaternion}

In this section, the smooth Hamiltonian function $H:\R^4\to \R$ does not need to be of the form
in equation \eqref{H1}. As before, $\Sg_E$ denotes a level set of $H$. It will be
supposed that $E$ is a regular value of $H$. Recall that in coordinates $(y_1,y_2,x_1,x_2)\in \R^4,$ the standard symplectic form is given by $\omega_0=dy_1\wedge dx_1 + dy_2\wedge dx_2$. Let $I_2$, $J_2$, and $0_2$ be the following
$2\times 2$ matrices:
\begin{equation}
I_2=\left(\begin{array}{ll}
1 & 0\\
0 & 1 \end{array}
\right), \quad J_2=  \left(\begin{array}{cc}
0 & -1\\
1 & 0 \end{array}
\right),\quad
0_2=\left(\begin{array}{ll}
0 & 0\\
0 & 0 \end{array}
\right)
 \label{matrix1}
\end{equation}
Using the blocks $I_2$,  $J_2$, and $0_2$    we define the following $4\times 4$ matrices:
\begin{equation}\label{matrix2}\begin{aligned}
  I & =\left(\begin{array}{ll}
I_2 & 0_2\\
0_2 & I_2 \end{array}
\right),
& \quad
A_1= \left(\begin{array}{ll}
\ 0_2 & -J_2\\
-J_2 & \  0_2 \end{array}
\right),\\
 A_2 & = \left(\begin{array}{cc}
J_2& \ 0_2\\
0_2 &-J_2  \end{array}\right),
& \quad A_3  =J=\left(\begin{array}{cc}
0_2 & -I_2\\
I_2 & 0_2 \end{array}
\right).
\end{aligned} \end{equation}
The matrices $A_1,A_2,A_3$ are anti-symmetric, $A_kA_k=-I$, for $k=1,2,3$, and
they multiply as $A^iA^j=\ep_{ijk}A^k$, where $\ep_{ijk}$ is the usual anti-symmetric tensor of
Levi-Civita. Since $E$ is a regular value of $H$, $\nabla H$
 never vanishes on  $\Sg_E$. Let $V_0= \nabla H/ |\nabla H|$ be a unit vector normal to $\Sigma_E$. Fixing the oriented basis
$\{\partial_{y_1},\partial_{y_2},\partial_{x_1},\partial_{x_2}\}$ on $\R^4$ we define the unit vector fields $V_k=A_k V_0$ on $\Sg_E$, for $k=1,2,3$. More precisely, we have
\begin{equation*}
\begin{aligned}
&V_0=|\nabla H|^{-1}(\partial_{y_1} H,\partial_{y_2} H,\partial_{x_1} H,\partial_{x_2} H),
&V_1=|\nabla H|^{-1}(\partial_{x_2} H,-\partial_{x_1} H,\partial_{y_2} H,-\partial_{y_1} H);\\
&V_2=|\nabla H|^{-1}(-\partial_{y_2} H,\partial_{y_1} H,\partial_{x_2} H,-\partial_{x_1} H),
&V_3=|\nabla H|^{-1}(-\partial_{x_1} H,-\partial_{x_2} H,\partial_{y_1} H,\partial_{y_2} H).
\end{aligned}
\end{equation*}
If $(\cdot,\cdot)$ is the Euclidean  inner product on $\R^4$ then
$(V_k,V_j)=(A_kV_0,A_jV_0)=0,$ for every $k\neq j$, that implies that $V_1,V_2,V_3$ form an orthonormal basis of each tangent space to
$\Sg_E$. Moreover,  $V_3=J\nabla H$ is parallel to
the Hamiltonian vector field $X$ on $\Sg_E$, and $\omega_0(V_1,V_2)=1$.

Now, given any  periodic orbit $\gamma$ on $\Sg_E$ we can compute its rotation number with respect
to $V_1,V_2,V_3$ in the following way. Given a point $p\in\gamma$, let
$\gamma(t)$ be a time parametrization of $\gamma$ such that $\gamma(0)=p$, and let $L(t)$ be the Hessian matrix
of $H$ along $\gamma$. Using the Cartesian frame $\{\partial_{y_1}, \partial_{y_2},\partial_{x_1}, \partial_{x_2}\}$, we write the linearized Hamiltonian equations over $\gamma$ as $\dot W=JL(t)W$.
For a given initial condition of the form  $W(0)=z_1(0)V_1(p)+z_2(0)V_2(p)$,
let  $W(t)$ be the corresponding solution to the above linearized equations.
The $V_1$ and $V_2$ components of $W$ are given by
$z_1(t)=(V_1(\gamma(t)),W(t))$ and $z_2(t)=(V_2(\gamma(t)),W(t))$, respectively.
To simplify the notation, let
$V_j(t)=V_j(\gamma(t))$, for $j=0,1,2,3$.
Then, differentiating
$z_1(t), z_2(t)$ with respect to time
we get, after some computations, that
$$
\begin{pmatrix}
\dot z_1\\
\dot z_2
\end{pmatrix}
=
J_2
\begin{pmatrix}
(V_1,L V_1)+(V_3,L V_3)&  (V_1,L V_2)\\
(V_2,L V_2)& (V_2,L V_2)+(V_3,L V_3)
\end{pmatrix}
\begin{pmatrix}
z_1\\
z_2 \end{pmatrix}.
$$
Writing $\theta(t)$ for the argument of $z_1(t)+ iz_2(t)$
($\dot \theta(t)$ is counter-clockwise oriented)
we obtain from the above equation that
\begin{equation}\label{thdot}
\dot \theta=
\left\{(V_3,L V_3)+
\begin{pmatrix}
\cos \theta\\
\sin\theta  \end{pmatrix}^T
\begin{pmatrix}
(V_1,L V_1)&  (V_1,L V_2)\\
(V_2,L V_1)& (V_2,L V_2)
\end{pmatrix}
\begin{pmatrix}
\cos \theta\\
\sin\theta  \end{pmatrix}
\right\}.
\end{equation}
Now, if $\theta(t)$ is any solution to this equation then  the rotation number $\rho_\gamma$ associated with
$\gamma$, defined in equation \eqref{rho}, is given by
\begin{equation*}
\rho_\gamma=T\lim_{t\to \infty}\frac{\theta(t)}{2\pi t},
\end{equation*}
where $T$ is the period of $\gamma$.
This limit always exists and does not rely on the initial condition $\theta(0)$. %

\hfill\newline
\noindent{\bf Acknowledgement.}
GR is partially supported by FAPESP grant  2023/07076-4. LL is partially support by National Natural Science Foundation of China (Grant number: 12071255, 12401238), Natural Science Foundation of Shandong Province, China (Grant number: ZR2024QA188). LL thanks the support of the School of Mathematics at Shandong University. PS is partially supported by the National Natural Science
Foundation of China (Grant number: W2431007).

\end{document}